\documentclass[10pt]{article}
\usepackage[square,authoryear]{natbib}
\usepackage{graphicx}
\graphicspath{{figures/}}
\usepackage{marsden_article}
\usepackage{multirow}

\newtheorem{example}[theorem]{Examples}
\newtheorem{remark}[theorem]{Remark}

\numberwithin{equation}{section}

\DeclareMathOperator*{\ext}{ext}

\begin{document}



\title{Variational multirate integrators}

\author{Sina Ober-Blöbaum\\ Department of Mathematics\\Paderborn University\\D-33098 Paderborn\\Germany\\
\and
Theresa Wenger\\85057 Ingolstadt \\Germany\\
\and
Tobias Gail\\ Treuchtlingen\\Germany\\
\and
Sigrid Leyendecker\\ Institute of Applied Dynamics\\Friedrich-Alexander-Universität Erlangen-Nürnberg\\D-91058 Erlangen\\Germany}





\maketitle

\begin{abstract}The simulation of systems that act on multiple time scales is challenging. A stable integration of the fast dynamics requires a highly accurate approximation whereas for the simulation of the slow part, a coarser approximation is accurate enough. With regard to the general goals of any numerical method, high accuracy and low computational costs, a popular approach is to treat the slow and the fast part of a system differently. Embedding this approach in a variational framework is the keystone of this work. By paralleling continuous and discrete variational multirate dynamics, integrators are derived on a time grid consisting of macro and micro time nodes that are symplectic, momentum preserving and also exhibit good energy behaviour. The choice of the discrete approximations for the action determines the convergence order of the scheme as well as its implicit or explicit nature for the different parts of the multirate system. The convergence order is proven using the theory of variational error analysis.
The performance of the multirate variational integrators is demonstrated by means of several examples.
\end{abstract}


\section{Introduction}
\label{intro}

Many systems in nature or in technology exhibit dynamics on varying time scales, imposing a challenge for numerical integration. Tiny step sizes are required to guarantee a stable integration of the fast motion, whereas for the slow motion, a larger time step is accurate enough. Typical examples of such multirate systems can be found in planetary dynamics where the different time scales originate from the extreme range of mass and distance of the planets. On the other end of the scale highly oscillatory systems arise in molecular dynamics, where locally extremely high frequencies superpose global folding processes. Further technological systems that are composed of rigid and elastic parts with varying and in particular with high stiffness are widely encountered, e.g.~in vehicle dynamics or in biomechanics. A large variety of integration methods have been developed to efficiently solve such multirate systems comprising fast and slow dynamics. The methods distinguish with respect to the simulation goals and their adaption to the structure of the underlying system.

One way to efficiently simulate multirate systems is to average what happens on the micro scale and to pass this information to the macro scale.
HMM (heterogeneous multiscale methods) \citep{weinan_heterogeneous_2007} rely on a top-down strategy, where the missing information is filled in an incomplete model on the macro scale by estimating what happens on the micro scale through averaging. This provides a general framework for designing and analysing very heterogeneous, multiscale or even multiphysics problems. An overview of these methods is found in \citep{weinan_heterogeneous_2007} and the application to highly oscillatory systems in \citep{brumm_heterogeneous_2014}. FLAVORS (flow averaging integrators \citep{tao_nonintrusive_2010}) are formulated using variational methods and the average of the flow is performed via a splitting and resynchronisation technique. In \citep{ober2013variational} the FLAVOR approach is applied to electric circuits. Another example of averaging methods are projective integration methods. They may have many levels of time scales and are independent of the possible integration method. 
These were introduced by Kevrekidis et al.~in \citep{kevrekidis03,kevrekidis_04} and in \citep{rico2004coarse}. A second order method of this approach is found in \citep{lee_second-order_2007}.
Multiple time stepping methods are developed for systems where the computationally expensive parts of the right-hand side do not contribute much to the dynamics of the solution.                                                                                                         The idea is to evaluate the expensive part of the vector field less often than the rest. One famous representative of this approach is the impulse method, also known as Verlet-I method or r-RESPA \citep{TBM_RESPA,GHWS_RESPA}. 
Extensions of r-RESPA to a mollified version are introduced by \citep{izaguirre1999longer} proposing to evaluate the slow force at an averaged value to reduce the resonance instability problem of the impulse method. 
Asynchronous variational integrators are developed from a Lagrangian, variational point of view enabling to use different time steps for different potential energy contributions and in different parts of the domain \citep{lew03_1,fong08}. The key feature of asynchronous variational integrators over many of the multiple time stepping methods is that time steps in arbitrary ratios can be considered. This is especially useful in solid dynamics simulations, where time steps can be made to vary smoothly throughout the spatial mesh. 
A lot of other multiple time stepping methods have been proposed in computational mechanics, commonly known as subcycling methods, see e.g.~\citep{gravouil_multi-time-step_2003,daniel_partial_2003}.
Instead of different time steps, IMEX methods combine implicit integration for the fast force with explicit integration for the slow force, see \citep{crouzeix,ZHANG1997297}. The goal is to capture the slow dynamics accurately without resolving the fast one. The IMEX method can also be formulated in a variational splitting framework, see \citep{SternGrin}. An overview of numerical methods for oscillatory, multiscale Hamiltonian systems is given e.g.~in \citep{cohen06}.
Often, both, different time stepping schemes and different time steps are combined in the integration method. An overview of heterogeneous (different time stepping schemes) asynchronous (different time steps) time integrators (HATI) for computational structural dynamics is given in \citep{gravouil_heterogeneous_2015}. 
Further multirate methods are the Multirate Partitioned Runge-Kutta Methods \citep{guenther01}, {Multirate generalized additive Runge Kutta methods \citep{gunther2016multirate}}, the extrapolated multirate methods \citep{constantinescu_extrapolated_2013} and the multirate linear multistep methods presented in \citep{gear_multirate_1984} where the unknowns are separated into fast and slow degrees of freedom. Different coupling strategies in multirate linear multistep methods are discussed in \citep{verhoeven_bdf_2008} amongst others. 
There also exist multirate schemes based on backward differentiation formulas (BDF), applied e.g.~to electric circuit systems \citep{verhoeven_bdf_2008,striebel_multirate_2009} or to mechanical problems \citep{arnold07}.

Most of the mentioned methods do not focus on the preservation of the underlying system's structure. An elegant and flexible approach to derive structure preserving integrators is based on a discrete version of Hamilton's principle. Approximating the continuous action integral and requesting the discrete counterpart to become stationary yields the integrator, i.e.~the discrete Euler-Lagrange equations. A detailed introduction as well as a survey on the history and literature on the variational view of discrete mechanics is given in \citep{MaWe01}. Due to the variational derivation, the integrator is symplectic. Via a discrete analogue to Noethers Theorem it can be shown that variational integrators preserve momentum maps. The variational approach has been successfully applied in different fields. Variational integrators have been developed for classical conservative mechanical systems (for an overview see \citep{lew04,lew03_2}), constrained systems (holonomic \citep{leyendecker08_1,leyendecker09} and non-holonomic systems \citep{sukhatme_geometric_2009}), nonsmooth systems \citep{fetecau03}, stochastic systems \citep{bourabee08}, forced \citep{kane00}
 and controlled \citep{oberbloebaum09} systems. Multiscale systems are considered in \citep{tao_nonintrusive_2010}.

In this work, the integrators are derived in a variational framework and designed to efficiently approximate the dynamics of multirate systems. The methods combine the separation of variables into fast and slow degrees of freedom with the use of different time steps in their discretization and different quadrature rules to approximate the contributions of the action. One compound step is used to solve for all unknowns enabling the method to be extended to solve index-3 DAEs with constraints coupling the slow and fast variables, see \citep{leyendecker13_1}. The presented approach includes the IMEX method of \citep{SternGrin} and a special case of the asynchronous variational integrators, however in general the schemes are different.

The partitioning into fast and slow components a priori may be difficult and may even vary with time. There are strategies how to detect slow and fast components dynamically during the simulation process based on local error approximations, see e.g.~\citep{engstler1997multirate}, \citep{gear_multirate_1984}, and to adapt the orders or the time step lengths accordingly. 
In this work, the split is given a priori, since, the keystone of this work is the fundamental derivation of multirate integrators from a variational Lagrangian point of view and its analysis. {We focus on low order methods to capture qualitative behaviour exactly as it is often done for molecular dynamics.}

The main contributions of this paper are the following. The variational multirate integrator is derived and its structure preserving properties {as symplecticity, momentum preservation in terms of symmetries} are shown. Depending on the choices of the quadrature rules, different integration schemes are derived in a $(p,q)$-formulation and related to existing methods. In Section \ref{Con_analysis}, the convergence orders of the multirate variational integrators are proven. The linear stability is discussed in Section \ref{sec:stability}. 
Section \ref{sec:numerical_Examples} illustrates the theoretical results by means of two dynamical example systems, followed by computing time investigations demonstrating the efficiency of the presented approach {where we observe decreasing computing times limited by the implicit nature of the integrator.}

\section{Lagrangian multirate dynamics}

Consider an $n$-dimensional mechanical system defined in $Q$, where $Q$ is an open subset of $\mathbb{R}^n$, with configuration vector ${q}(t) \in
Q$ and velocity vector $\dot{{q}}(t) \in T_{{q}(t)}Q$,
where $t$ denotes the time variable in the bounded interval $[t_0, t_N]\subset
\mathbb{R}$. The Lagrangian $L \colon TQ \rightarrow \mathbb{R}$ of a mechanical system is given by the difference of the kinetic energy $T(\dot{{q}})$ and a potential {energy} $U({q})$. Let the fact that the Lagrangian contains slow and fast dynamics be characterized by the possibility to additively split the potential energy $U(q)=V(q)+W(q)$ into a slow potential $V$ and a fast potential $W$, where $W(q)$ can be, for example, of the form $W(q) = 1/\varepsilon \, \bar{W}(q),\, \varepsilon \ll 1$. In general, consider possible potential splits according to, for example, stiffness, nonlinearity, dynamical behavior, and evaluation cost \citep{gunther2016multirate}. 
Let $\mathcal{C}(Q)=\mathcal{C}([t_0, t_N], Q, {q}_{0}, {q}_{N})$ denote the space of
smooth curves $q \colon [t_0, t_N]\rightarrow Q$ satisfying ${q}(t_0)={q}_{0}$ and ${q}(t_N)={q}_{N}$, where 
${q}_{0}, {q}_{N} \in Q$ are fixed endpoints. For $q\in \mathcal{C}(Q)$, the action integral is defined as
\begin{equation}\label{eq:ac}
\displaystyle \mathfrak{S}({q})=\int_{t_0}^{t_N} L({q}, \dot{{q}}) \, dt.
\end{equation}
Requiring that the first variation of this action vanishes, i.e.~$\delta \mathfrak{S}=0$, Hamilton's principle of stationary action yields the Euler-Lagrange equations of motion of a conservative mechanical system.
We further assume that the $n$-dimensional configuration variable $q$ can be divided into $n^s$ slow variables $q^s\in Q^s$ and $n^f$ fast variables $q^f\in Q^f$ such that $Q^s \times Q^f=Q$ and $q=(q^s,q^f)^T$ with $n^s+n^f=n$. Let the fast potential depend on the fast degrees of freedom only, i.e.~$W=W(q^f)$, while the slow potential $V=V(q)$ depends on the complete configuration variable. 
The kinetic energy is assumed to be of the form $T(\dot{q}) = \frac{1}{2} \dot{q}^T M \dot{q}$ with a block-diagonal mass matrix $M=\left( \begin{matrix}  M^s & 0 \\ 0 &  M^f\end{matrix}\right)$.
With these assumptions, the regular multirate Lagrangian $L \colon TQ^s \times TQ^f \rightarrow \mathbb{R}$ (which in slight abuse of notation is also denoted by $L$) is given by
\begin{equation}\label{eq:multirateL}
L(q^s, q^f,\dot{q}^s, \dot{q}^f) = T(\dot{q}^s, \dot{q}^f) - V(q^s, q^f) - W(q^f)
\end{equation}
and the Euler-Lagrange equations take the form
\begin{equation} \label{eq:EL_TVW_sf}
\begin{split}
 \frac{\partial V}{\partial q^s} + \frac{d}{dt} \frac{\partial T}{\partial \dot{q}^s}  & = 0\\
 \nabla W +  \frac{\partial V}{\partial q^f} + \frac{d}{dt} \frac{\partial T}{\partial \dot{q}^f}   & = 0
\end{split}
\end{equation}
and can equivalently be written in Hamiltonian form as
\begin{align}\label{eq:Ham}
\begin{split}
\dot{q}^s &=  (M^{s} ) ^{-1} p^s,\\
\dot{p}^s & = -\frac{\partial V}{\partial q^s},\\
\dot{q}^f &=  ( M^{f} )^{-1} p^f,\\
\dot{p}^f & = -\frac{\partial V}{\partial q^f} - \nabla W
\end{split}
\end{align}
where $  p^s =   \frac{\partial T}{\partial \dot{q}^s} $, $  p^f =   \frac{\partial T}{\partial \dot{q}^f}$ and with $\nabla$ being the nabla operator.
 
\begin{remark}
If in addition, the slow potential depends on the slow variables only and on top of that the kinetic energy does not contain any entries coupling $\dot{q}^s$ and $\dot{q}^f$, then the system is completely decoupled and simulation can be performed independently in parallel, without any exchange of information. However, we focus on the coupled scenario described above. Note that the inclusion of additional potentials depending on slow variables only is straightforward. 
\end{remark}

In the following, we briefly present some well known results of Lagrangian dynamics that also hold for the Lagrangian multirate dynamics defined above. In the multirate setting, the standard expressions can be decomposed and sorted into terms taking the splitting of the variables and potentials into slow and fast parts into account.

Following \citep{MaWe01}, the mapping $D_\text{EL} L$ given in coordinates as $(D_\text{EL} L)_i=\frac{\partial L}{\partial {q^i}} -  \frac{d}{dt}\frac{\partial L}{\partial \dot{q}^i} $ is called Euler-Lagrange map and the Lagrangian vector field $X_L \colon TQ \rightarrow T(TQ)$ with its Lagrangian flow $F_L \colon TQ \times \mathbb{R}\rightarrow TQ$ satisfies $D_\text{EL} L \circ X_L =0$.

As shown in \citep{MaWe01}, the symplectic form $\Omega_L \colon TQ \rightarrow \mathbb{R} $ is preserved under the Lagrangian flow $F_L^t \colon TQ \rightarrow TQ$ (for a certain $t$)
\[
(F_L^t)^*(\Omega_L) = \Omega_L\quad \text{for all}\, t\in [t_0,t_N],
\]
making the Lagrangian flow $F_L$ symplectic.

Furthermore, the Lagrangian flow $F_L$ conserves momentum maps, what can be shown via Noether's theorem.
Suppose that a Lie group $G$ with Lie algebra $\mathfrak{g}$ acts on $Q $ by the action $\Phi \colon G\times Q \rightarrow Q$. 
If the Lagrangian $L \colon TQ \rightarrow \mathbb{R}$ is invariant under the lift of the action, the corresponding momentum map $J_L \colon TQ \rightarrow \mathfrak{g}^*$, where $\mathfrak{g}^*$ is the Lie algebra of the Lie group $g$, is a conserved quantity of the flow so that $J_L \circ F_L^t = J_L$ for all times $t$. A proof can be found in e.g.~\citep{MaWe01}.

\begin{remark}
In principle, the Lie group actions on slow and fast variables can be different and for that special case be divided into a slow and a fast action $\Phi^s \colon G^s \times Q^s \rightarrow Q^s$ and $\Phi^f \colon G^f \times Q^f \rightarrow Q^f$, with $ \Phi = (\Phi^s, \Phi^f)^T$.
For example, consider a three-dimensional pendulum (fast rotational degrees of freedom) being attached to a cart (slow translational degrees of freedom) that can move in a plane being perpendicular to the direction of gravity. The corresponding Lagrangian is invariant with respect to translation of the slow degrees of freedom, i.e.~$G^s=\mathbb{R}^2$, and rotation of the fast ones around the gravity axis, i.e.~$G^f\subset SO(3)$. 
\end{remark}

\section{Construction of multirate variational integrators}
To derive the multirate variational integrators, we transfer the continuous multirate Lagrangian dynamics to the discrete setting. Symplecticity of the discrete Lagrangian multirate flow is proven and a discrete version of Noether's Theorem is presented for the multirate case. Depending on the choice of the approximation rules, different integrating schemes are derived.

\subsection{Discrete variational principle}

\begin{figure}[htbp]
\centerline{\includegraphics[width=90mm]{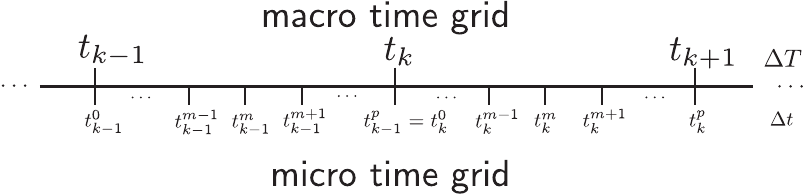}}
\caption{Macro and micro time grid}
\label{fig:time_grids}
\end{figure}

Rather than choosing one time grid for the approximation as for standard variational integrators, for the multirate integrator, two different time grids are introduced, see Fig.~\ref{fig:time_grids}.
With the time steps $\Delta T$ and $\Delta t$ (where $\Delta T=p \Delta t$, $p \in \mathbb{N}$), a macro time grid $\{t_k=k\Delta T \, | \, k=0,\ldots,N\}$ and a micro time grid $\{t^m_k=k\Delta T + m \Delta t \, | \, k=0,\ldots,N-1, m=0,\ldots,p\}$ are defined. Note that except for the boundary nodes $t_0, t_N$, two micro time nodes coincide with a macro time node, i.e.~$t_{k-1}^p=t_k^0=t_k$ for $k=1,\ldots,N-1$, see Fig.~\ref{fig:time_grids}.
The macro time grid provides the domain for the discrete macro trajectory of the slow variables 
\begin{equation*}
  \mathsf{q}^s_d   =\{q_k^s\}_{k=0}^N \quad \textrm{ with } \quad q^s_k \approx q^s(t_k)  
\end{equation*}
and the discrete micro trajectory of the fast variables lives on the micro grid
\begin{equation*}
 \mathsf{q}^f_d=\{q^f_k\}_{k=0}^{N-1}= \{\{q^{f,m}_k\}_{m=0}^{p}\}_{k=0}^{N-1} \quad \textrm{ with } \quad q^{f,m}_k \approx q^f(t^m_k)
\end{equation*}
Note that due to $t_{k-1}^{p}=t_{k}^{0}$ also $q_{k-1}^{f,p}=q_{k}^{f,0}$ holds.
The discrete macro trajectory of the variable $q$ defined on the macro grid is denoted by 
\begin{equation*}
 \mathsf{q}_d=  \lbrace q_k \rbrace_{k=0}^N \quad \textrm{ with } \quad q_k = (q^s_k, q^{f,0}_k )^T. 
\end{equation*} 

As an approximation to $\mathfrak{S}$ in \eqref{eq:ac}, the discrete action is defined as 
\begin{equation} \label{eq:discrete_action}
\mathfrak{S}_d(q_d^s, q^{f}_d)=  \sum_{k=0}^{N-1} L_d(q_k^s, q_{k+1}^s,q_{k}^{f}) 
\end{equation}
The discrete Lagrangian $L_d= T_d - V_d - W_d$ approximates $\int_{t_k}^{t_{k+1}}L(q, \dot{q}) \, dt$ and reads
\begin{equation} \label{eq:discrete_Lagrangian}
L_d(q_k^s, q_{k+1}^s,q_{k}^{f})=T_d(q_k^s, q_{k+1}^s,q_{k}^{f}) -V_d(q_k^s, q_{k+1}^s,q_{k}^{f}) - W_d(q_{k}^{f}).
\end{equation}
As we will see in Section~\ref{sec:approx}, a useful representation of $L_d(q_k^s, q_{k+1}^s,q_{k}^{f})$ is of the form
\begin{equation}\label{eq:Ldm}
L_d(q_k^s, q_{k+1}^s,q_{k}^{f}) = \sum_{m=0}^{p-1} L_d^m(q_k^s,q_{k+1}^s, q_k^{f,m}, q_{k}^{f,m+1}),
\end{equation}
where $L_d^m \approx \int_{t_k^m}^{t_{k}^{m+1}}L(q, \dot{q}) \, dt$ approximates the action on the micro interval $[t_k^m,t_k^{m+1}]$.

Omitting the arguments of $L_d$, the discrete Hamilton's principle leads to stationarity of the discrete action given by
\begin{align}\label{eq:discHam}
0 =& \delta \mathfrak{S}_d =  d \mathfrak{S}_d  \cdot (\delta q_d^s, \delta q_d^f) = \sum_{k=0}^{N-1}  \left( D_{q_k^s}  L_d   \cdot \delta q_k^s +  D_{q_{k+1}^s} L_d\cdot \delta q_{k+1}^s +  \sum_{m=0}^{p}  
 D_{q_{k}^{f,m}} L_d \cdot\delta q_{k}^{f,m} \right)\nonumber\\
  =& \sum_{k=1}^{N-1}  D_{q_k^s} \left(  L_d(q_k^s,{q_{k+1}^s}, q_k^{f}) +  L_d(q_{k-1}^s,{q_{k}^s}, q_{k-1}^{f}) \right) \cdot \delta q_k^s \nonumber\\
  & + \sum_{k=1}^{N-1} D_{q_k^{f,0}} \left(  L_d(q_k^s,{q_{k+1}^s}, q_k^{f})+ L_d(q_{k-1}^s,{q_{k}^s}, q_{k-1}^{f}) \right) \cdot \delta q_k^{f,0}\\
 & + \sum_{k=0}^{N-1} \sum_{m=1}^{p-1} D_{q_k^{f,m}} L_d(q_k^{s},{q_{k+1}^{s}},{q_k^{f}}) \cdot \delta q_k^{f,m} \nonumber\\
 & + D_{q_0^s}   L_d(q_0^s,{q_{1}^s}, q_0^{f}) \cdot \delta q_0^s + D_{q_N^s}   L_d(q_{N-1}^s,{q_{N}^s}, q_{N-1}^{f}) \cdot \delta q_N^s \nonumber\\
 & + D_{q_0^{f,0}} L_d(q_0^{s},{q_{1}^{s}},{q_0^{f}}) \cdot \delta q_0^{f,0} + D_{q_{N-1}^{f,p}} L_d(q_{N-1}^{s},{q_{N}^{s}},{q_{N-1}^{f}}) \cdot \delta q_{N-1}^{f,p} \nonumber
\end{align}
with independent variations $\delta q_k^s$ for $k=0,\ldots N$ and $\delta q_{k}^{f,m}$ for $k=0,\ldots,N-1$, $m=0,\ldots,p$ and fixed boundaries $\delta q_0^s=\delta q_N^s =0$ and $\delta q_{0}^{f,0} = \delta q_{N-1}^{f,p}=0$ and where $D_*$ denotes the partial derivative with respect to the argument $*$.
This yields the discrete multirate Euler-Lagrange equations (DEL) given by
\begin{equation}\label{eq:DEL}
\begin{array}{llccl}
\text{(DEL)}_k^s & D_{q_k^s} \left(  L_d(q_k^s,{q_{k+1}^s}, q_k^{f}) +  L_d(q_{k-1}^s,{q_{k}^s}, q_{k-1}^{f}) \right)& = &0, & k=1,\ldots,N-1,\\[3pt]
\text{(DEL)}_k^{f,0}& D_{q_k^{f,0}} \left(  L_d(q_k^s,{q_{k+1}^s}, q_k^{f})+ L_d(q_{k-1}^s,{q_{k}^s}, q_{k-1}^{f}) \right) &= &0,& k=1,\ldots,N-1,\\[3pt]
\text{(DEL)}_k^{f,m}& D_{q_k^{f,m}} L_d(q_k^{s},{q_{k+1}^{s}},{q_k^{f}})  &=&0, & k=0,\ldots,N-1 \\
& & & & m = 1,\ldots,p-1.
\end{array}
\end{equation}

The discrete slow and fast momenta on the macro nodes are defined by
\begin{align}
p_k^{s,-} &= -D_{q_k^s}   L_d(q_k^s,{q_{k+1}^s}, q_k^{f}), & p_k^{s,+} &= D_{q_k^s}   L_d(q_{k-1}^s,{q_{k}^s}, q_{k-1}^{f}),\label{eq:dis_ps}\\
p_k^{f,0,-} &= -D_{q_k^{f,0}}   L_d(q_k^s,{q_{k+1}^s}, q_k^{f}), & p_{k-1}^{f,p,+} &= D_{q_{k-1}^{f,p}}   L_d(q_{k-1}^s,{q_{k}^s}, q_{k-1}^{f})\label{eq:dis_pf}
\end{align}
and we denote the left and right momenta at the macro nodes consisting of slow and fast components by $p_k^- = \left(\begin{matrix} p_k^{s,-} \\ p_k^{f,0,-} \end{matrix}\right)$ and $p_k^+ = \left(\begin{matrix} p_k^{s,+} \\ p_{k-1}^{f,p,+}  \end{matrix}\right)$. 
If the discrete Lagrangian is of the form \eqref{eq:Ldm}, we define the discrete fast momenta on the micro nodes as
\begin{align}
p_k^{f,m,-} &= -D_{q_k^{f,m}} L_d^m(q_k^{s},{q_{k+1}^{s}},{q_k^{f,m}},{q_k^{f,m+1}}), \quad m=0,\ldots,p-1.\label{eq:pkfm-}\\
p_k^{f,m,+} &= D_{q_k^{f,m}} L_d^m(q_k^{s},{q_{k+1}^{s}},{q_k^{f,m-1}},{q_k^{f,m}}), \quad m=1,\ldots,p.\label{eq:pkfm+}
\end{align}
Note that there are no discrete slow momenta defined on micro nodes since we approximate the slow configuration only at macro nodes.

Now the discrete Euler-Lagrange equations \eqref{eq:DEL} can be rewritten in such a way that we have momenta matching at macro and micro nodes as
\begin{equation}\label{eq:DELp}
\begin{array}{llccl}
\text{(DEL)}_k^s & -p_k^{s,-} +  p_k^{s,+} & = &0, & k=1,\ldots,N-1,\\[3pt]
\text{(DEL)}_k^{f,0}& -p_k^{f,0,-} + p_{k-1}^{f,p,+}  &= &0,& k=1,\ldots,N-1,\\[3pt]
\text{(DEL)}_k^{f,m}& -p_k^{f,m,-} +  p_k^{f,m,+}  &=&0, & k=0,\ldots,N-1,\, m = 1,\ldots,p-1, 
\end{array}
\end{equation}
and we can define the unique momenta per macro and micro node as $p_k^s := p_k^{s,-} = p_k^{s,+}$, $p_k^{f,0} := p_k^{f,0,-} = p_{k-1}^{f,p,+}$, $p_{k-1}^{f,p} := p_k^{f,0,-} = p_{k-1}^{f,p,+}$, $p_k := p_k^- = p_k^+$ and $p_k^{f,m}:= p_k^{f,m,-}=p_k^{f,m,+}$.

As usual for variational integrators, for setting up the resulting integration scheme, we have to consider the first step separately.
Let $k=0$ and assume that an initial configuration $(q_0^s,q_0^{f,0})$ and an initial conjugate momentum $(p_0^s,p_0^{f,0})$ are given. Then for $k=0$, the unknowns ${q_1^s, q_0^{f,1},\ldots,q_0^{f,p}}$ are determined by solving the following set of equations for $m=1,\ldots,p-1$.
\begin{equation}\label{eq:DEL_ini}
\begin{array}{lrcl}
\text{(IC)}^s & D_{q_0^s}  L_d(q_0^s,{q_1^s}, q_0^{f}) &=& -p_0^s \\[3pt]
\text{(IC)}^f & D_{q_0^{f,0}} L_d(q_0^s,{q_1^s}, q_0^{f})&=& -p_0^{f,0}\\[3pt]
\text{(DEL)}_0^{f,m} & D_{q_0^{f,m}} L_d(q_0^{s},{q_1^{s}},{q_0^{f}}) &=&0.
 \end{array}
\end{equation}
These equations can be considered as initial conditions, since they determine the unknowns in the first macro time interval from given initial data. 
To proceed further in time for $k=1,\ldots,N-1$ (assuming that $q_{k-1}^s,q_k^s,q_{k-1}^{f,0},\ldots,q_{k-1}^{f,p}$ are given), solving the discrete Euler-Lagrange equations \eqref{eq:DEL} determines $q_{k+1}^s$, $q_k^{f,1},\ldots,q_k^{f,p}$.

The multirate discrete Lagrangian in \eqref{eq:discrete_Lagrangian} can be formed to a standard discrete Lagrangian, see \citep{MaWe01}. 

\begin{paragraph}{Equivalence of the multirate discrete Lagrangian to a standard discrete Lagrangian}
Define
\begin{align*}
L_d(q_k ,q_{k+1} )   
= \ext_{\begin{subarray}{c}  ( q^{f,1}_{{k}}, \ldots, q^{f,p-1}_{{k}} )    \end{subarray} } L_d ( q^{{s}}_k, q^{{s}}_{k+1}, q^f_k ) 
\end{align*}
where $q_k = ( q^s_k, q^{f,0}_k )$, $q_{k+1} = ( q^s_{k+1}, q^{f,p}_k )$. 
The standard {(two-point)} discrete Lagrangian $L_d(q_k ,q_{k+1})   $ 
is the multirate discrete Lagrangian defined in \eqref{eq:discrete_Lagrangian}, 
evaluated on the trajectory within the step which solves $\text{(DEL)}^{f,m}_k$ in \eqref{eq:DEL}. {Note that with a slight abuse of notation we use the same symbol $L_d$ to define the multirate Lagrangian and the Lagrangian depending only on the macro nodes.} 
\end{paragraph}

\begin{paragraph}{Discrete multirate Lagrangian flow}
The discrete multirate Lagrangian flow $F^{\Delta T}_{L_d} : Q \times Q \rightarrow Q \times Q$ is defined by scheme \eqref{eq:DEL} as
\begin{align} \label{eq:discLagFlow}
F^{\Delta T}_{L_d} (q_{k-1}, q_{k} ) =  (q_{k}, q_{k+1} )
\end{align}
and depends only on the macro nodes assumed that $q^{f,m}_{k}$, $m=1, \ldots,p-1$ satisfy the discrete Euler-Lagrange equations $\text{(DEL)}_k^{f,m}$ in \eqref{eq:DEL}. 
\end{paragraph}

\subsection{Structure preserving properties of the discrete Lagrangian flow}
\label{sec:conservation_properties}

The last two lines of the discrete Hamilton's principle \eqref{eq:discHam} provide the discrete slow and fast one-forms on $Q^s\times Q^s$ and $Q^f \times Q^f$, respectively, defined by
\begin{align}\label{eq:Theta_Ld}
\Theta_{L_d}^{s,-}(q_0^s,q_1^s,q_0^{f,0},q_0^{f,p}) &= -D_{q_0^s}   L_d(q_0^s,{q_{1}^s}, q_0^{f}) dq_0^s,  \\
 \Theta_{L_d}^{f,-}(q_0^s,q_1^s,q_0^{f,0},q_0^{f,p}) &= -D_{q_0^{f,0}}   L_d(q_0^s,{q_{1}^s}, q_0^{f}) dq_0^{f,0},\\
 \Theta_{L_d}^{s,+}(q_{N-1}^s,q_N^s,q_{N-1}^{f,0},q_{N-1}^{f,p}) &= D_{q_N^s}   L_d(q_{N-1}^s,{q_{N}^s}, q_{N-1}^{f}) dq_N^s,\\ 
 \Theta_{L_d}^{f,+}(q_{N-1}^s,q_N^s,q_{N-1}^{f,0},q_{N-1}^{f,p}) &= D_{q_{N-1}^{f,p}}   L_d(q_{N-1}^s,{q_{N}^s}, q_{N-1}^{f}) dq_{N-1}^{f,p},
\end{align}
depending on the macro nodes only assumed that $q^{f,m}_{k}$, $m=1, \ldots,p-1$, satisfy the discrete Euler-Lagrange equations $\text{(DEL)}_k^{f,m}$ in \eqref{eq:DEL}. The discrete one-forms can be summarized to $\Theta_{L_d}^{-}  = \Theta_{L_d}^{s,-}  +\Theta_{L_d}^{f,-}  $ and $\Theta_{L_d}^{+}  = \Theta_{L_d}^{s,+}  + \Theta_{L_d}^{f,+} $  
and are denoted by discrete Lagrangian one-forms on $Q\times Q$.

The discrete symplectic form is given by $\Omega_{L_d}(q_k, q_{k+1} ) = d \Theta^-_{L_d} = d \Theta^{+}_{L_d}$.  
In the following theorem it is shown that this symplectic form is preserved under the discrete multirate Lagrangian flow by means of a generating function argument.

\begin{theorem}
The discrete multirate Lagrangian flow \eqref{eq:discLagFlow} is symplectic.
\end{theorem}
\begin{proof} 
The proof takes the standard steps, see e.g.~\citep{HW}. Define the discrete action $S_d \colon Q \times Q \rightarrow \mathbb{R}$ as a function of the boundary values $q_0, q_N$, where $(q_k, q_{k+1}, \lbrace q^{f,m}_k \rbrace^p_{m=0} )^{N-1}_{k=0}$ is the solution of the discrete Euler-Lagrange equations \eqref{eq:DEL} with these boundary values. Determining the differential of $S_d$ shows that $S_d$ is a symplecticity generating function and thus the discrete multirate Lagrangian flow and the integration scheme are symplectic.  
\end{proof}
Symplectic integrators show a very good energy behaviour, meaning that the energy is not decreasing or increasing artifically over simulation time. The energy error is oscillating, but stays bounded even for very long simulation times. This favourable property of symplectic integrators can be shown by using techniques from backward error analysis, cf.~e.g.~\citep{hairer2006geometric}, and is verified numerically in Section \ref{sub:energies_and_preserved_quantities}.

A discrete version of Noether's theorem is proven in \citep{MaWe01} for the single rate case, i.e.~for $p=1$. Considering the multirate case, i.e.~$p>1$, the proof works analogously.
\begin{theorem}[Discrete Noether's theorem]\label{th_dNoether}
Consider a discrete multirate Lagrangian $L^{\Delta T}_d \colon Q \times Q \rightarrow \mathbb{R}$ which is invariant under the lift of the action $\Phi \colon G \times Q \rightarrow Q$. Then there exists a unique discrete momentum map $J_{L_d} \colon Q \times Q \rightarrow \mathfrak{g}^*$ which is a conserved quantity of the discrete multirate Lagrangian flow $F^{\Delta T}_{L_d} \colon Q \times Q \rightarrow Q \times Q $, so that $J_{L_d} \circ F^{\Delta T}_{L_d} = J_{L_d}$.  
\end{theorem}
\begin{proof}
The proof works analogously as in \citep{MaWe01}.  
\end{proof}

\subsection{Choice of approximation rules}\label{sec:approx}

Following the construction of Galerkin variational integrators (see e.g.~\citep{ober2014construction}), for the discrete formulation, we first approximate the space of trajectories and second, use quadrature rules to approximate the integral of the Lagrangian. In general, the space of piecewise smooth polynomials of order $l$, $q_d(t)\in \Pi^l$, for the trajectory space approximation together with a quadrature rule $(b_j,c_j)_{j=1}^r$ with $r$ quadrature points can be chosen. Having a clear separation of slow and fast variables, we are able to approximate them by different polynomials defined on different subintervals, i.e.~$q^s(t) \approx q_d^s(t)\in \Pi^{l^s}$ on the macro interval $[t_k,t_{k+1}]$ and $q^f(t) \approx q_d^{f,m}(t)\in \Pi^{l^f}$ on the micro interval $[t_k^m,t_{k}^{m+1}]$. Furthermore, it is also possible to use different quadrature rules, $(c_j,b_j)_{j=1}^r$, $(\tilde{c}_j,\tilde{b}_j)_{j=1}^{\tilde{r}}$ for the slow and fast contributions in the Lagrangian, such that a general approximation formula for the integral of the Lagrangian in \eqref{eq:multirateL} reads as 
\begin{align}
& L_d(q_k^s,q_{k+1}^s, q_k^f)  =  T_d(q_k^s,q_{k+1}^s, q_k^f) - V_d(q_k^s,q_{k+1}^s, q_k^f) - W_d(q_k^f) :=  \label{L_d} \\
&\approx   \sum_{m=0}^{p-1} \int_{t_k^m}^{t_{k}^{m+1}} T(\dot{q}^s(t),\dot{q}^f(t)) - V(q^s(t),q^f(t)) - W(q^f(t))\, dt  \nonumber \\
&= \int_{t_k}^{t_{k+1}} L(q^s(t),q^f(t),\dot{q}^s(t),\dot{q}^{f}(t))\, dt  \nonumber
\end{align}
with 
\begin{align*}
T_d(q_k^s,q_{k+1}^s, q_k^f)&= \Delta t  \sum_{m=0}^{p-1} \sum_{j=1}^r b_j T(\dot{q}_d^s(m\Delta t +c_j \Delta t), \dot{q}_d^{f,m}({m\Delta t} +c_j \Delta t)), \\
V_d (q_k^s,q_{k+1}^s, q_k^f)&= \Delta t  \sum_{m=0}^{p-1} \sum_{j=1}^r b_j V(q_d^s(m\Delta t + c_j \Delta t),q_d^{f,m}({m\Delta t} +c_j \Delta t)), \\
W_d (q_k^f)& =\Delta t  \sum_{m=0}^{p-1} \sum_{j=1}^{\tilde{r}} \tilde{b}_j W(q_d^{f,m}({m\Delta t} +\tilde{c}_j \Delta t))  
\end{align*} 
being the discrete kinetic energy, the discrete slow and the discrete fast potential {energy}.

The quadrature rules in use for the discrete kinetic energy, the discrete slow and the discrete fast potential {energy} determine the degree of coupling between the discrete Euler-Lagrange equations \eqref{eq:DEL}. This can range from a fully implicit scheme over variants being explicit in the slow and implicit in the fast quantities to fully explicit schemes as we will see in the following.
{However, when approximating the integral of the slow potential on the micro grid as in \eqref{L_d} the resulting scheme is always implicit and a nonlinear system of equations has to be solved per macro time step.}
To simplify matters in this work (which is in contrast to e.g.~\citep{wenger_variational_2016},
where polynomials of various higher degrees are used for the approximation of the slow and the fast degrees of freedom), in the following 
we restrict to the case of degree one polynomials, i.e.~$l^s=l^f=1$ together with either midpoint rule quadrature (being a Gauss quadrature rule with one quadrature point, i.e.~$b_1=1, c_1 =\frac{1}{2}$) {or the Lobatto quadrature rule with two quadrature points, i.e.~$b_1= b_2 = \frac{1}{2}, c_1=0, c_2=1$ corresponding to trapezoidal rule (also denoted as trapezoidal rule with $\alpha_{\{V,W\}}=\frac{1}{2}$, compare equation~\eqref{eq:slow_potential_midpoint_2}), or the rectangle rule based on a left respectively right point evaluation, i.e.~one quadrature point with $b_1=1, c_1=0$ respectively $b_1=1, c_1 =1$ (also denoted as trapezoidal rule with $\alpha \in \{ 0,1 \}$, compare equation~\eqref{eq:slow_potential_midpoint_2})}.
Using first order polynomials, the approximation polynomials for the slow and fast configurations and velocities are defined as

\begin{align} 
q_d^s(t)&= q_{k}^s + \frac{q_{k+1}^s-q_k^s}{\Delta T} (t - t_k), &\dot{q}_d^s(t) &=  \frac{q_{k+1}^s-q_k^s}{\Delta T} ,& t&\in [t_k,t_{k+1}]\quad \text{and}   \label{eq:qds_linear}  \\ 
q_d^{f,m}(t)&= q_{k}^{f,m} + \frac{q_{k}^{f,m+1}-q_k^{f,m}}{\Delta t} (t - t^m_k), &\dot{q}_d^{f,m}(t)& = \frac{q_{k}^{f,m+1}-q_k^{f,m}}{\Delta t}, & t&\in [t_k^m,t_{k}^{m+1}]  \label{eq:qdf_linear}
\end{align}
such that $q_d^s(t_k)= q_k^s$, $q_d^s(t_{k+1})= q_{k+1}^s$, $q_d^{f,m}(t_k^m)= q_k^{f,m}$, $m=0,\ldots,p$, and we further define $q_k^{s,m} := q_d^s(t_k^m)$ for $m=1,\ldots,p-1$.
Note that since the approximated slow and fast velocities are constant on each macro and micro interval, respectively, the discrete kinetic energy $T_d$ is always the same independent of the choice of quadrature rule since $L$ is assumed to be separable. In the following we assume that the kinetic energy can be decomposed in a contribution from the fast and the slow variables, i.e.
\begin{equation}\label{eq:T}
T(\dot{q}) =  \frac{1}{2} \dot{q}^T \cdot M \cdot  \dot{q}  = \frac{1}{2} \left(\dot{q}^s\right)^T \cdot  M^s  \cdot \dot{q}^s + \frac{1}{2} \left(\dot{q}^f\right)^T  \cdot M^f  \cdot \dot{q}^f,
\end{equation}
where $M^s$ and $M^f$ are the mass matrices for the slow and fast variables, respectively, yielding the total mass matrix as $M=\textrm{diag}(M^s, M^f)$.
Then the discrete kinetic energy is defined on the time interval $[t_k,t_{k+1}]$ as 
\begin{equation*}
 T_d(q_k^s,q_{k+1}^s, q_k^f)=  \frac{\Delta T}{2} \left( \frac{q^s_{k+1}-q^s_k}{\Delta T} \right)^T \cdot M^s  \cdot  \left( \frac{q^s_{k+1}-q^s_k}{\Delta T} \right)  +    \sum_{m=0}^{p-1} \frac{\Delta t}{2} \left( \frac{q^{f,m+1}_{k}-q^{f,m}_k}{\Delta t} \right)^T  \cdot M^f  \cdot  \left( \frac{q^{f,m+1}_{k}-q^{f,m}_k}{\Delta t} \right) 
\end{equation*}

Depending on the choice of quadrature rules for the slow and fast potential we obtain the following {implicit} integration schemes for configuration and momentum.

\begin{theorem}[midpoint fast, midpoint slow]\label{th:midfmids}
If we choose the midpoint rule for fast and slow potential as
\begin{eqnarray}
V_d(q_k^s,{q_{k+1}^s}, q_k^{f}) &=& \sum_{m=0}^{p-1} \Delta t \, V\left(\frac{q_{k}^{s,m}+q_{k}^{s,m+1}}{2},\frac{q_{k}^{f,m}+q_{k}^{f,m+1}}{2}\right)\label{eq:slow_potential_midpoint_1}\\
W_d(q_k^{f}) &=&  \sum_{m=0}^{p-1} \Delta t \, W\left(\frac{q_{k}^{f,m}+q_{k}^{f,m+1}}{2}\right),\label{eq:fast_potential_midpoint_1}
\end{eqnarray}
the resulting $(p,q)$-scheme is

\begin{align}
q^s_{k+1} & = q_k^s + \frac{\Delta T}{2} M^{s^{-1}} \left( p_k^{s} + p_{k+1}^{s} - \Delta t \sum_{m=0}^{p-1} \left(1-\frac{2m+1}{p}\right) \frac{\partial }{\partial q^s} V(\bar{q}_k^{s,m} , \bar{q}_k^{f,m} )  \right)   \label{midSlow1} \\
 p^s_{k+1} & = p_k^s - \Delta t \sum_{m=0}^{p-1} \frac{\partial }{\partial q^s} V(\bar{q}_k^{s,m} , \bar{q}_k^{f,m} ) \label{midSlow2}
\end{align}
\begin{align}
q^{f,m+1}_{k} & = q_k^{f,m} + \Delta t  M^{f^{-1}} \frac{p_k^{f,m} + p_{k}^{f,m+1}}{2}\\
 p^{f,m+1}_{k} & = p_k^{f,m}  -\Delta t \frac{\partial }{\partial q^f} V(\bar{q}_k^{s,m} , \bar{q}_k^{f,m} )- \Delta t \nabla W(\bar{q}_k^{f,m} ), 
\end{align}
with $\bar{q}_k^{f,m} := \frac{{q}_k^{f,m}+{q}_k^{f,m+1}}{2}$ and $\bar{q}_k^{s,m} := \frac{{q}_k^{s,m}+{q}_k^{s,m+1}}{2}$,
or equivalently, written as transformed $(p,q)$-scheme

\begin{align}
\tilde{p}_k^s & = p_k^s -  \Delta t \sum_{m=0}^{p-1} \left(1-\frac{2m+1}{p}\right) \frac{\partial }{\partial q^s} V(\bar{q}_k^{s,m} , \bar{q}_k^{f,m} ) \label{eq:tildeps}\\
q^s_{k+1} & = q_k^s + \Delta T  M^{s^{-1}} \frac{ \tilde{p}_k^{s} + p_{k+1}^{s}}{2} \label{eq:qpt_q}\\
 p^s_{k+1} & = \tilde{p}_k^s - \Delta t \sum_{m=0}^{p-1} \frac{2m+1}{p}\frac{\partial }{\partial q^s} V(\bar{q}_k^{s,m} , \bar{q}_k^{f,m} )\label{eq:qpt_p1}
\end{align}

\begin{align}
q^{f,m+1}_{k} & = q_k^{f,m} + \Delta t  M^{f^{-1}} \frac{p_k^{f,m} + p_{k}^{f,m+1}}{2}\\
 p^{f,m+1}_{k} & = p_k^{f,m}  -\Delta t \frac{\partial }{\partial q^f} V(\bar{q}_k^{s,m} , \bar{q}_k^{f,m} )- \Delta t \nabla W(\bar{q}_k^{f,m} )
\end{align}

\end{theorem}

\begin{proof}
First we derive the discrete Euler-Lagrange equations which are according to \eqref{eq:DEL} 
\begin{align*}
0  = &- M^s \frac{q_{k+1}^s-q_k^s}{\Delta T} +  M^s \frac{q_{k}^s-q_{k-1}^s}{\Delta T} - \frac{\Delta t}{2} \left[ \sum_{m=0}^{p-1} \frac{\partial}{\partial q^s} V(\bar{q}_k^{s,m},\bar{q}_k^{f,m}) \frac{2p-(2m+1)}{p} + \frac{\partial}{\partial q^s} V(\bar{q}_{k-1}^{s,m},\bar{q}_{k-1}^{f,m}) \frac{2m+1}{p} \right]\\
0  =& - M^f \frac{q_{k}^{f,1}-q_k^{f,0}}{\Delta t} +  M^f \frac{q_{k-1}^{f,p}-q_{k-1}^{f,p-1}}{\Delta t} - \frac{\Delta t}{2} \left[  \frac{\partial}{\partial q^f} V(\bar{q}_k^{s,0},\bar{q}_k^{f,0})  + \frac{\partial}{\partial q^f} V(\bar{q}_{k-1}^{s,p-1},\bar{q}_{k-1}^{f,p-1})\right]\\
&- \frac{\Delta t}{2} \left[  \nabla W(\bar{q}_k^{f,0})  + \nabla W(\bar{q}_{k-1}^{f,p-1}) \right]\\
0  =& - M^f \frac{q_{k}^{f,m+1}-q_k^{f,m}}{\Delta t} +  M^f \frac{q_{k}^{f,m}-q_{k}^{f,m-1}}{\Delta t} - \frac{\Delta t}{2} \left[  \frac{\partial}{\partial q^f} V(\bar{q}_k^{s,m},\bar{q}_k^{f,m})  + \frac{\partial}{\partial q^f} V(\bar{q}_k^{s,m-1},\bar{q}_k^{f,m-1})\right]\\
&- \frac{\Delta t}{2} \left[  \nabla W(\bar{q}_k^{f,m})  + \nabla W(\bar{q}_k^{f,m-1}) \right],\, m=1,\ldots,p-1
\end{align*}
Now we define the slow and fast momenta according to \eqref{eq:dis_ps}--\eqref{eq:pkfm+} as
\begin{align*}
p_k^{s,-} = &  M^s \frac{q_{k+1}^s-q_k^s}{\Delta T} +   \frac{\Delta t}{2} \sum_{m=0}^{p-1} \frac{\partial}{\partial q^s} V(\bar{q}_k^{s,m},\bar{q}_k^{f,m}) \frac{2p-(2m+1)}{p} \\
p_k^{s,+} = &  M^s \frac{q_{k}^s-q_{k-1}^s}{\Delta T} - \frac{\Delta t}{2}  \sum_{m=0}^{p-1} \frac{\partial}{\partial q^s} V(\bar{q}_{k-1}^{s,m},\bar{q}_{k-1}^{f,m}) \frac{2m+1}{p}  \\
p_k^{f,0,-} = &   M^f \frac{q_{k}^{f,1}-q_k^{f,0}}{\Delta t}  + \frac{\Delta t}{2}   \frac{\partial}{\partial q^f} V(\bar{q}_k^{s,0},\bar{q}_k^{f,0}) + \frac{\Delta t}{2}  \nabla W(\bar{q}_k^{f,0})   \\
p_{k-1}^{f,p,+} = &    M^f \frac{q_{k-1}^{f,p}-q_{k-1}^{f,p-1}}{\Delta t} - \frac{\Delta t}{2}  \frac{\partial}{\partial q^f} V(\bar{q}_{k-1}^{s,p-1},\bar{q}_{k-1}^{f,p-1})- \frac{\Delta t}{2} \nabla W(\bar{q}_{k-1}^{f,p-1}) \\
p_k^{f,m,-} = &  M^f \frac{q_{k}^{f,m+1}-q_k^{f,m}}{\Delta t} + \frac{\Delta t}{2}   \frac{\partial}{\partial q^f} V(\bar{q}_k^{s,m},\bar{q}_k^{f,m})  
+ \frac{\Delta t}{2}   \nabla W(\bar{q}_k^{f,m})  \\
p_k^{f,m,+} = &  M^f \frac{q_{k}^{f,m}-q_{k}^{f,m-1}}{\Delta t} - \frac{\Delta t}{2} \frac{\partial}{\partial q^f} V(\bar{q}_k^{s,m-1},\bar{q}_k^{f,m-1})- \frac{\Delta t}{2} \nabla W(\bar{q}_k^{f,m-1}) 
\end{align*}
Building the sums $p_k^{s,-}+p_{k+1}^{s,+}$, $p_k^{f,m,-}+p_{k}^{f,m+1,+}$ and differences $p_k^{s,-}-p_{k+1}^{s,+}$, $p_k^{f,m,-}-p_{k}^{f,m+1,+}$ of the discrete slow and fast momenta together with momenta matching $p^{s,-}_k = p^{s,+}_k$ and $p^{f,m,-}_k = p^{f,m,+}_k$, see \eqref{eq:DELp}, yields the desired scheme. The transformed $(p,q)$-scheme immediately follows by introducing the auxiliary slow momentum variable $\tilde{p}_k^s$ as defined in \eqref{eq:tildeps}. \qed
\end{proof}
By comparing the $(p,q)$-scheme with the continuous Hamiltonian equations~\eqref{eq:Ham}, we observe that the resulting scheme corresponds to a midpoint rule scheme for the fast dynamics as well as for the slow dynamics if $p=1$ as expected for Gauss quadrature rules. The slow transformed $(p,q)$-scheme \eqref{eq:tildeps}--\eqref{eq:qpt_p1} resembles a St\"ormer-Verlet method in onestep formulation (see e.g.~\citep{hairer2006geometric}), where $q$ is updated by the midpoint rule and averaged potentials on the the macro interval $[t_k,t_{k+1}]$ are used to update the intermediate momentum $\tilde{p}_k^s$ and final momentum $p_{k+1}^s$

\begin{theorem}[midpoint fast, trapezoidal slow]\label{th:midftraps} 
If we choose the midpoint rule for the fast and the trapezoidal rule for the slow potential as

\begin{equation}\label{eq:slow_potential_midpoint_2}
V_d(q_k^s,{q_{k+1}^s}, q_k^{f}) = \sum_{m=0}^{p-1} \Delta t \, \left( {\alpha_V} V(q_{k}^{s,m},q_{k}^{f,m}) + (1- {\alpha_V} ) V(q_{k}^{s,m+1},q_{k}^{f,m+1})\right)
\end{equation}
\begin{equation} \label{eq:fast_potential_midpoint_2}
W_d(q_k^{f}) =  \sum_{m=0}^{p-1} \Delta t \, W\left(\frac{q_{k}^{f,m}+q_{k}^{f,m+1}}{2}\right)
\end{equation}
with ${\alpha_V } \in[0,1]$, the resulting $(p,q)$-scheme is

\begin{align}
q^s_{k+1} & = q_k^s + \Delta T  M^{s^{-1}} \left[ p_{k}^{s} - \Delta t \left( \sum_{m=1}^{p-1} \frac{p-m}{p} \frac{\partial }{\partial q^s} V({q}_k^{s,m} , {q}_k^{f,m} ) + {\alpha_V} \frac{\partial }{\partial q^s} V({q}_k^{s,0} , {q}_k^{f,0} ) \right) \right]  \label{trapSlow1} \\
 p^s_{k+1} & = p_k^s - \Delta t \left(\sum_{m=1}^{p-1} \frac{\partial }{\partial q^s} V({q}_k^{s,m} , {q}_k^{f,m} ) + {\alpha_V} \frac{\partial }{\partial q^s} V({q}_k^{s,0} , {q}_k^{f,0} )  + (1- {\alpha_V}) \frac{\partial }{\partial q^s} V({q}_k^{s,p} , {q}_k^{f,p} )\right)   \label{trapSlow2}
\end{align}

\begin{eqnarray}
\text{kick}&\tilde{p}^{f,m}_k & = {p}^{f,m}_k - {\alpha_V} \Delta t   \frac{\partial }{\partial q^f} V({q}_k^{s,m} , {q}_k^{f,m} )\\
\text{oscillate}& q^{f,m+1}_{k} & = q_k^{f,m} + \Delta t  M^{f^{-1}} \frac{\tilde{p}_k^{f,m} + \hat{p}_{k}^{f,m+1}}{2}\\
 & \hat{p}^{f,m+1}_{k} & = \tilde{p}_k^{f,m}  - \Delta t \nabla W(\bar{q}_k^{f,m} )\\
 \text{kick}& p_{k}^{f,m+1}& =  \hat{p}^{f,m+1}_{k}  - (1- {\alpha_V}) \Delta t   \frac{\partial }{\partial q^f} V({q}_k^{s,m+1} , {q}_k^{f,m+1} ),
\end{eqnarray}

or equivalently, written as transformed $(p,q)$-scheme

\begin{align}
\tilde{p}_k^{s} & =  p_{k}^{s} - \Delta t \left( \sum_{m=1}^{p-1} \frac{p-m}{p} \frac{\partial }{\partial q^s} V({q}_k^{s,m} , {q}_k^{f,m} ) + {\alpha_V}  \frac{\partial }{\partial q^s} V({q}_k^{s,0} , {q}_k^{f,0} ) \right) \label{eq:pktilde} \\
q^s_{k+1} & = q_k^s + \Delta T  M^{s^{-1}} \tilde{p}_k^{s}\\
 p^s_{k+1} & = \tilde{p}_k^s - \Delta t \left( \sum_{m=1}^{p-1} \frac{m}{p}\frac{\partial }{\partial q^s} V({q}_k^{s,m} , {q}_k^{f,m} ) + (1- {\alpha_V})  \frac{\partial }{\partial q^s} V({q}_k^{s,p} , {q}_k^{f,p} ) \right)\label{eq:pks+1}
\end{align}
\begin{eqnarray}
\text{kick}&\tilde{p}^{f,m}_k & = {p}^{f,m}_k - {\alpha_V} \Delta t   \frac{\partial }{\partial q^f} V({q}_k^{s,m} , {q}_k^{f,m} )\\
\text{oscillate}& q^{f,m+1}_{k} & = q_k^{f,m} + \Delta t  M^{f^{-1}} \frac{\tilde{p}_k^{f,m} + \hat{p}_{k}^{f,m+1}}{2}\\
 & \hat{p}^{f,m+1}_{k} & = \tilde{p}_k^{f,m}  - \Delta t \nabla W(\bar{q}_k^{f,m} )\\
 \text{kick}& p_{k}^{f,m+1}& =  \hat{p}^{f,m+1}_{k}  - (1- {\alpha_V}) \Delta t   \frac{\partial }{\partial q^f} V({q}_k^{s,m+1} , {q}_k^{f,m+1} )
 \label{kickoscckick_fastPart}
\end{eqnarray}

\end{theorem}

The resulting fast scheme in {\eqref{kickoscckick_fastPart}} corresponds to a kick-oscillate-kick scheme, where the kick is only influenced by the slow potential and the oscillating part is only influenced by the fast potential, see \citep{hairer2006geometric}. The oscillating part is updated here with the midpoint rule.
Furthermore, for $p=1$, ${\alpha_V} = 0$ we obtain a symplectic Euler scheme for the slow dynamics and for ${\alpha_V} = 1$ the adjoint symplectic Euler scheme.
For $p=1$, ${\alpha_V} =\frac{1}{2} $ and recombining the slow and fast variables, the scheme results in the IMEX method \citep{SternGrin}. Similar as for the midpoint fast and midpoint slow scheme, the slow transformed $(p,q)$-scheme \eqref{eq:pktilde}--\eqref{eq:pks+1} resembles a St\"ormer-Verlet method in onestep formulation, where again averaged potentials on the macro interval $[t_k,t_{k+1}]$ are used to update the intermediate momentum $\tilde{p}_k^s$ and final momentum $p_{k+1}^s$

\begin{theorem}[trapezoidal fast, trapezoidal slow]\label{th:trapftraps}
If we choose the trapezoidal rule for fast and slow potential as

\begin{equation}
V_d(q_k^s,{q_{k+1}^s}, q_k^{f}) = \sum_{m=0}^{p-1} \Delta t \, \left( \alpha_{{V}} V(q_{k}^{s,m},q_{k}^{f,m}) + (1-\alpha_{{V}}) V(q_{k}^{s,m+1},q_{k}^{f,m+1})\right)  \label{eq:slow_potential_midpoint_3}
\end{equation}

\begin{equation} \label{eq:fast_potential_midpoint_3}
W_d(q_k^{f}) =  \sum_{m=0}^{p-1} \Delta t \left(\alpha_{{W}} W(q_{k}^{f,m}) + (1-\alpha_{{W}})W(q_{k}^{f,m+1}) \right)
\end{equation}

with $\alpha_{{V,W}}  \in[0,1]$, the resulting $(p,q)$-scheme is

\begin{align}
q^s_{k+1} & = q_k^s + \Delta T  M^{s^{-1}} \left[ p_{k}^{s} - \Delta t \left( \sum_{m=1}^{p-1} \frac{p-m}{p} \frac{\partial }{\partial q^s} V({q}_k^{s,m} , {q}_k^{f,m} ) + \alpha_{{V}} \frac{\partial }{\partial q^s} V({q}_k^{s,0} , {q}_k^{f,0} ) \right) \right] \label{trapSlow3} \\
 p^s_{k+1} & = p_k^s - \Delta t \left(\sum_{m=1}^{p-1} \frac{\partial }{\partial q^s} V({q}_k^{s,m} , {q}_k^{f,m} ) + \alpha_{{V}} \frac{\partial }{\partial q^s} V({q}_k^{s,0} , {q}_k^{f,0} )  + (1-\alpha_{{V}}) \frac{\partial }{\partial q^s} V({q}_k^{s,p} , {q}_k^{f,p} )\right) \label{trapSlow4}
\end{align}

\begin{align}
{q^{f,m+1}_{k}} & = { q_k^{f,m} + \Delta t  M^{f^{-1}}  \left[  p_{k}^{f,m} -  \Delta t \left( \alpha_V \frac{\partial }{\partial q^f} V({q}_k^{s,m} , {q}_k^{f,m} ) + \alpha_W \nabla W({q}_k^{f,m} )\right)  \right] } \\
  { p^{f,m+1}_{k} } & = { p_k^{f,m}  -\Delta t \left[  \left( \alpha_V \frac{\partial }{\partial q^f} V({q}_k^{s,m} , {q}_k^{f,m} )+ \alpha_W \nabla W({q}_k^{f,m} )\right) \right.  }   \\
 & { \left.  +     \left(  (1-\alpha_V) \frac{\partial }{\partial q^f} V({q}_k^{s,m+1} , {q}_k^{f,m+1} )+ (1-\alpha_W)  \nabla W({q}_k^{f,m+1} )\right)  \right] }
\end{align}

or equivalently, written as transformed $(p,q)$-scheme

\begin{align}
\tilde{p}_k^{s} & =  p_{k}^{s} - \Delta t \left( \sum_{m=1}^{p-1} \frac{p-m}{p} \frac{\partial }{\partial q^s} V({q}_k^{s,m} , {q}_k^{f,m} ) + { \alpha_V}  \frac{\partial }{\partial q^s} V({q}_k^{s,0} , {q}_k^{f,0} ) \right) \\
q^s_{k+1} & = q_k^s + \Delta T  M^{s^{-1}} \tilde{p}_k^{s}\\
 p^s_{k+1} & = \tilde{p}_k^s - \Delta t \left( \sum_{m=1}^{p-1} \frac{m}{p}\frac{\partial }{\partial q^s} V({q}_k^{s,m} , {q}_k^{f,m} ) + (1- {  \alpha_V})  \frac{\partial }{\partial q^s} V({q}_k^{s,p} , {q}_k^{f,p} ) \right)
\end{align}

\begin{align}
{q^{f,m+1}_{k}} & = { q_k^{f,m} + \Delta t  M^{f^{-1}}  \left[  p_{k}^{f,m} -  \Delta t \left( \alpha_V \frac{\partial }{\partial q^f} V({q}_k^{s,m} , {q}_k^{f,m} ) + \alpha_W \nabla W({q}_k^{f,m} )\right)  \right] } \\
 { p^{f,m+1}_{k} } & = { p_k^{f,m}  -\Delta t \left[  \left( \alpha_V \frac{\partial }{\partial q^f} V({q}_k^{s,m} , {q}_k^{f,m} )+ \alpha_W \nabla W({q}_k^{f,m} )\right) \right.  }   \\
 & { \left.  +     \left(  (1-\alpha_V) \frac{\partial }{\partial q^f} V({q}_k^{s,m+1} , {q}_k^{f,m+1} )+ (1-\alpha_W)  \nabla W({q}_k^{f,m+1} )\right)  \right] }
\end{align}
\end{theorem}

The resulting scheme for the fast dynamics corresponds for ${\alpha_V = \alpha_W  = 0}$ to the symplectic Euler and for ${\alpha_V = \alpha_W  =1}$ to the adjoint symplectic Euler scheme as expected. For $p=1$ the same correspondence holds also for the slow dynamics scheme which is the same as the slow scheme for the midpoint fast trapezoidal slow scheme.

Proofs of Theorems~\ref{th:midftraps} and \ref{th:trapftraps} follow the same lines as the proof of Theorem~\ref{th:midfmids}, i.e.~by deriving discrete Euler-Lagrange equations, discrete momenta and rearranging the discrete momenta equations.

Note that all the schemes presented so far are implicit due to the interpolation of the slow configurations at the micro nodes which always depend on the unknown configuration $q_{k+1}^s$ for polynomial approximations $q_d^s\in \Pi^{l_s}$ for $l_s\ge 1$. 
Other integration schemes can be derived by approximating the integral of the slow potential using only macro nodes rather than using also (interpolated) micronodes for slow and fast configurations. 
When standard variational integrators are used for problems with, for example, very stiff potentials, their discrete counterparts are often based on midpoint evaluations of the continuous potentials such that the corresponding integration scheme is implicit. On the other hand, the integral of softer potentials can be approximated by evaluations of the continuous potential on the left or right node yielding explicit schemes (at least as long as there are no constraints present), which are of course much cheaper regarding the computational costs. 
Different combination{s} are given in Examples~\ref{ex:slow_ex__fast_im} and \ref{ex:fully_ex} (see also \citep{leyendecker13_1}). Depending on the complexity of the evaluation of the potential functions and their gradients, the computational costs of the overall simulation is heavily influenced by the choice of quadrature.

\begin{example}[Implicit fast and explicit slow forces]\label{ex:slow_ex__fast_im}
Choosing an affine combination as approximation in the slow potential that involves only macro nodes
\begin{equation}\label{eq:slow_potential_affine}
V_d(q_k^s,{q_{k+1}^s}, q_k^{f}) =  \Delta T\left( {\alpha_V} V((q_k^s, q_k^{f,0})) + (1- {\alpha_V}) V((q_{k+1}^s,q_k^{f,p}))\right)
\end{equation}
with ${\alpha_V} \in[0,1]$ and a micro node based midpoint rule in the fast potential 
\begin{equation} \label{eq:fast_potential_midpoint_4}
W_d(q_k^{f}) =  \sum_{m=0}^{p-1} \Delta t \, W\left(\frac{q_{k}^{f,m}+q_{k}^{f,m+1}}{2}\right)
\end{equation}
leads to discrete conservative forces in the discrete Euler-Lagrange equations which are explicit for the slow potential and implicit for the fast one. Thus, only few evaluations of the gradient of $V$ are necessary which is advantageous when the slow potential's evaluation is very costly compared to the fast one. The resulting scheme can be interpreted as a variational splitting method which is symmetric and symplectic, since it is a symmetric composition of symmetric and symplectic methods. When this method is formulated with ${\alpha_V} =\frac{1}{2}$ on only one time grid with a constant time step (i.e.~$\Delta t = \Delta T$ and $p=1$) and without splitting the configuration variable into fast and slow variables, one obtains the IMEX method in \citep{SternGrin} which is an example of an impulse method, see \citep{McLachStern}, \citep{hairer2006geometric} and references therein.
\end{example}

\begin{example}[Fully explicit scheme]\label{ex:fully_ex}
Using the affine combination of the slow potential evaluated at the macro nodes in \eqref{eq:slow_potential_affine} and the affine combination of micro node evaluations of the fast potential
\begin{equation} \label{eq:fast_potential_trapezoidal}
W_d(q_k^{f}) =  \sum_{m=0}^{p-1} \Delta t \left({\alpha_W} W(q_{k}^{f,m}) + (1-{\alpha_W})W(q_{k}^{f,m+1}) \right)
\end{equation}
with ${\alpha_{V,W}} \in[0,1]$ leads to discrete Euler-Lagrange equations \eqref{eq:DEL} that can subsequently be solved without iteration, i.e.~first $q_k^{f,1}$ is obtained from (DEL)$_k^{f,0}$, then (DEL)$_k^{f,m}$ yields $q_k^{f,m+1}$ for $m=1,\ldots,p-1$. At any time, $q_{k+1}^s$ can be computed from (DEL)$_k^s$. For ${\alpha_{V,W}} =\frac{1}{2}$, this choice of quadrature leads to the scheme in \citep{lew03_1} for the special case that a synchronized time grid is used there.
\end{example}

{
\section{Convergence analysis}}
\label{Con_analysis}

In this section, the order of accuracy of the variational multirate integrators is analysed. Rather than considering how closely the  discrete trajectory matches the exact trajectory, the variational error analysis considers how closely a discrete Lagrangian matches the exact discrete Lagrangian that is defined by
\begin{equation}
L^E_d(q_0, q_1, \Delta T) = \int^{\Delta T}_{0} L(q(t),\dot{q}(t)) dt
\label{eq:exactLd}
\end{equation}
for sufficiently small $\Delta T$ and close $q_0$ and $q_1$ and $q(t)$ being the unique solution of the Euler-Lagrange equations satisfying $q(0)=q_0$ and $q(\Delta T)=q_1$. The result of the variational error analysis presented in \citep{MaWe01} and refined in \citep{patrick2009error} is given in the next Theorem. 
\begin{theorem}[Variational error analysis]\label{th:varerror}
Let the discrete Lagrangian $L_d$ and the exact discrete Lagrangian $L^E_d$ be two discretizations of a regular Lagrangian $L$. If $L_d$ approximates $L_d^E$ to order $r$, i.e.~$L_d(q_0,q_1,\Delta T) = L^E_d(q_0,q_1, \Delta T) + \mathcal{O}(\Delta T^{r+1} )$, then the variational integrator induced by $L_d$ is of order $r$.
\end{theorem}

For the estimate of the discrete Lagrangian, we make use of the variational error analysis for Galerkin variational integrators {\cite{ober2021superconvergence}} where we have to adapt several steps for the multirate case.
In the following, we give a detailed proof along the lines of {\cite{ober2021superconvergence}}.
Let us introduce the space of polynomials
\begin{align*}
\mathcal{P}^l([0,\Delta T],Q) = \{ q \in C^\infty([0,\Delta T],Q) \,|\, q \text{ a polynomial of degree at most } l \} .
\end{align*}
We fix $l+1$ control points $\Delta T d_0 < \Delta T d_1 < \ldots < \Delta T d_l$, where $d_0 = 0$ and $d_l = 1$. If for each of these control points a value $q(\Delta T d_j) = q_j$ is prescribed, then the polynomial $q \in \mathcal{P}^l([0,\Delta T],Q)$ is uniquely determined. We denote by $\hat{q}(\,\cdot\, ; q_0,\ldots,q_l,\Delta T)$ the polynomial thus obtained. To analyse the multirate case, we also introduce the space of piecewise polynomials
\begin{align*}
\mathcal{P}_p^l([0,\Delta T],Q) =& \{ q \in PC^\infty([0,\Delta T],Q) \,|\, q|_{[m\Delta t,(m+1)\Delta t)} \in \mathcal{P}^l([m\Delta t,(m+1)\Delta t),Q), \,m=0,\ldots,p-2, \\
& q|_{[(p-1)\Delta t,p\Delta t]} \in \mathcal{P}^{\rm l}([(p-1)\Delta t,p\Delta t],Q)\\
&\text{ and } q \text{ is continuous at } m\Delta t,\, m=1,\ldots,p-1\} 
\end{align*} 
where $PC^\infty$ is the space of piecewise smooth functions.
Similar as before we denote by 
$\check{q}(\,\cdot\, ; \{q^m_0,\ldots,q^m_l\}_{m=0}^{p-1},\Delta T)$ the piecewise polynomial designed by $\check{q}(m\Delta t + \Delta t d_j) = q_j^m,\, m=0,\ldots,p-1$.
For the construction of the multirate integrator, we approximate the slow variable $q^{s}$ by a polynomial $\hat{q}^{s}$ defined on the macro interval $[0,\Delta T]$ and the fast variable $q^{f}$ by a polynomial on the micro interval $[m\Delta t, (m+1) \Delta t]$ which corresponds to piecewise polynomial $\check{q}^{f}$ on $[0,\Delta T]$.

For the error estimate we need the following lemma (for a proof see {\cite{ober2021superconvergence}}). 
\begin{lemma}\label{lemma-poly}
	Let $q$ be smooth curve and $\hat{q}$ a family of polynomials of degree $l$, parametrized by $\Delta T$, which equals $q$ at the control points $0 = \Delta T d_0 < \Delta T d_1 < \ldots < \Delta T d_l = \Delta T$.
	Then for any $k \leq l$ there holds
	\[ \big\| q^{(k)} - \hat{q}^{(k)} \big\|_\infty = \mathcal{O}(\Delta T^{l+1-k}) , \]
	where $\|\cdot\|_\infty$ denotes the maximum norm on $[0,\Delta T]$.
\end{lemma}

For the approximation of the action we introduce a quadrature rule. We fix quadrature points $c_i \in [0,1]$ and weights $b_i \in \mathbb{R}$, with $\sum_i b_i = 1$. We denote by $u$ the order of the corresponding quadrature rule. Then for any smooth function $f$ there holds 
\[ \int_0^{\Delta T} f(t) \,d t - \Delta T \sum_i b_i f(\Delta T c_i) = \mathcal{O}(\Delta T^{u+1}). \]
Based on the polynomial spaces and this quadrature rule we investigate the multirate discrete Lagrangian defined as 
\begin{align}\label{eq:hatLext}
&{L}_d(q_0,q_1) = \\  &\underset{\substack{\hat{q}^{s} \in \mathcal{P}^l([0,\Delta T],Q^{\rm s})\\ \check{q}^{f} \in \mathcal{P}_{p}^l([0,\Delta T],Q^{f})\\ \hat{q}^{s}(0)= q^{s}_0, \hat{q}^{s}(\Delta T)=q^{s}_1\\\check{q}^{f}(0)= q^{f,0}_0, \check{q}^{f}(\Delta T)=q^{f,p}_0}}{\textrm{ext}}\Delta t \sum_{m=0}^{ {p }-1}\sum_{i=1}^r b_i {L}(\hat{q}^{s}(m \Delta t + c_i \Delta t), \check{q}^{f}(m \Delta t + c_i \Delta t), \dot{\hat{q}}^{s}(m \Delta t + c_i \Delta t),\dot{\check{q}}^{f}(m \Delta t + c_i \Delta t) ) \notag .
\end{align}
which corresponds to the standard (two-point) discrete Lagrangian introduced in \eqref{L_d} with $l=l_s = l_f$, $\tilde{c}_i = \tilde{c}_i$ and $\tilde{b}_i=b_i$, $i=1,\ldots,r$. 
For the error analysis we assume that the continuous and discrete critical curves are not only extremisers but indeed minimise their respective actions and that minimisers are unique. This is the case in particular if the Lagrangian is of mechanical type and the quadrature rule is sufficiently accurate (see \cite{hall2015spectral}, Theorem 3.5).
%
Now we are ready to determine the approximation error of the discrete Lagrangian. 
\begin{theorem}
	\label{thm-superh}
	Let ${L}_d$ as given in \eqref{eq:hatLext} 
 be a multirate discretisation with of the Lagrangian ${L}$ 
 based on polynomials of degree $l$ and a quadrature rule of degree $u$. Assume that all discrete and continuous critical curves minimise their respective actions. 
 Then the difference between the discrete and the exact discrete Lagrangian is
 \begin{align}\label{eq:LELd}
{L}_{d}(q(0),q(\Delta T)) = L_{d}^E(q(0),q(\Delta T))+\mathcal{O}((\Delta T)^{2l+1}) +  p \mathcal{O}((\Delta t)^{u+1})
\end{align}
\end{theorem}

\begin{proof}
The proof strategy goes along the lines of {\cite{ober2021superconvergence}}.
	%
	Let $q_\mathrm{EL}$ denote the solution of the continuous forced Euler-Lagrange equations with $q_\mathrm{EL}(0) = q(0)$ and $q_\mathrm{EL}(\Delta T) = q(\Delta T)$.  Let $\hat{q}^{s} \in \mathcal{P}^l([0,\Delta T],Q^s)$ be the polynomial that agrees with $q_\mathrm{EL}^{s}$ at the control points $0 = \Delta T d_0 < \Delta T d_1 < \ldots < \Delta T d_s = \Delta T$ and $\check{q}^{f} \in \mathcal{P}_{p}^l([0,\Delta T],Q^f)$ be the piecewise polynomial that agrees with $q_\mathrm{EL}^{f}$ at the control points $0 = m\Delta t + \Delta t d_0 < m\Delta t + \Delta t d_1 < \ldots < m\Delta t + \Delta t d_s = \Delta T$ for $m=0,\ldots,p-1$. Furthermore, let $\tilde{q}^{s} \in \mathcal{P}^l([0,\Delta T],Q^{s})$ and $\tilde{q}^{f} \in \mathcal{P}_p^l([0,\Delta T],Q^{f})$ be the (piecewise) polynomials that minimise dicrete action
 \[
\Delta t \sum_{m=0}^{p-1}\sum_{i=1}^r b_i {L}(\tilde{q}^{s}(m \Delta t + c_i \Delta t), \tilde{q}^{f}(m \Delta t + c_i \Delta t),\dot{\tilde{q}}^{s}(m \Delta t + c_i \Delta t), \dot{\tilde{q}}^{f}(m \Delta t + c_i \Delta t) )
 \]
	in $\mathcal{P}^l([0,\Delta T],Q^{s}) \times \mathcal{P}_{p}^l([0,\Delta T],Q^{f})$. With \eqref{eq:exactLd}
	we have to estimate the relation
	\begin{equation}\label{toshow}
	\int_0^{\Delta T} {L}(q_\mathrm{EL},\dot{q}_\mathrm{EL}) \,d t - \Delta t \sum_{m=0}^{p-1}\sum_{i=1}^r b_i {L}(\tilde{q}^{s}(m \Delta t + c_i \Delta t), \tilde{q}^{f}(m \Delta t + c_i \Delta t), \dot{\tilde{q}}^{s}(m \Delta t + c_i \Delta t),\dot{\tilde{q}}^{f}(m \Delta t + c_i \Delta t) ) .
	\end{equation}
	With a slight abuse of notation using the same symbol $L$ for the functions $L(q,\dot{q})$ and  $L(q^s,q^f,\dot{q}^s,\dot{q}^f)$ we expand this difference as
 \begin{align}
	\label{difference-expansion}
	 &\left( \int_0^{\Delta T} {L}(q^{s}_\mathrm{EL},q^{f}_\mathrm{EL},\dot{q}^{s}_\mathrm{EL},\dot{q}^{f}_\mathrm{EL}) \,d t - \int_0^{\Delta T} {L}(\hat{q}^{s},\check{q}^{f},\dot{\hat{q}}^{s},\dot{\check{q}}^{f}) \,d t \right)  \\
	 &+ \left( \int_0^{\Delta T} {L}(\hat{q}^{s},\check{q}^{f},\dot{\hat{q}}^{s},\dot{\check{q}}^{f})\,d t - \Delta t \sum_{m=0}^{p-1}\sum_{i=1}^r b_i {L}(\tilde{q}^{s}(m \Delta t + c_i \Delta t), \tilde{q}^{f}(m \Delta t + c_i \Delta t), \dot{\tilde{q}}^{s}(m \Delta t + c_i \Delta t),\dot{\tilde{q}}^{f}(m \Delta t + c_i \Delta t)) \right).\nonumber
	\end{align}
	%
	We start with the first term of \eqref{difference-expansion}. From Lemma \ref{lemma-poly} we know that $q^{s} - \hat{q}^{s} = \mathcal{O}(\Delta T^{l+1})$ and $\dot{q}^{s} - \dot{\hat{q}}^{s} = \mathcal{O}(\Delta T^l)$ on $[0,\Delta T]$ and $q^{f} - \check{q}^{f} = \mathcal{O}(\Delta t^{l+1})$ and $\dot{q}^{f} - \dot{\check{q}}^{f} = \mathcal{O}(\Delta t^l)$ on $[m\Delta t ,(m+1)\Delta t], m=0,\ldots,p-1$, hence
	\begin{align*}
	\int_0^{\Delta T}& L(q^{s}_\mathrm{EL},q^{f}_\mathrm{EL},\dot{q}^{s}_\mathrm{EL},\dot{q}^{f}_\mathrm{EL}) \,d t - \int_0^{\Delta T} L(\hat{q}^{s},\check{q}^{f},\dot{\hat{q}}^{s},\dot{\check{q}}^{f}) \,d t  \\ 
	=& \int_0^{\Delta T} \bigg( \frac{\partial L(q^{s}_\mathrm{EL},q^{f}_\mathrm{EL},\dot{q}^{s}_\mathrm{EL},\dot{q}^{f}_\mathrm{EL})}{\partial q^{s}}(q^{s}_\mathrm{EL} - \hat q^{s}) +  \frac{\partial L(q^{s}_\mathrm{EL},q^{f}_\mathrm{EL},\dot{q}^{s}_\mathrm{EL},\dot{q}^{f}_\mathrm{EL})}{\partial \dot{q}^{s} } (\dot{q}^{s}_\mathrm{EL} - \dot{\hat q}^{s}) + \mathcal{O}\left(\Delta T^{2l}\right) \bigg) d t \\
 &+\sum_{m=0}^{p-1}\int_{m\Delta t}^{(m+1)\Delta t} \bigg( \frac{\partial L(q^{s}_\mathrm{EL},q^{f}_\mathrm{EL},\dot{q}^{s}_\mathrm{EL},\dot{q}^{f}_\mathrm{EL})}{\partial q^{f}}(q^{f}_\mathrm{EL} - \check q^{f}) +  \frac{\partial L(q^{s}_\mathrm{EL},q^{f}_\mathrm{EL},\dot{q}^{s}_\mathrm{EL},\dot{q}^{f}_\mathrm{EL})}{\partial \dot{q}^{f} } (\dot{q}^{f}_\mathrm{EL} - \dot{\check q}^{f}) + \mathcal{O}\left(\Delta t^{2l}\right) \bigg) d t \\
	=& \int_0^{\Delta T} \bigg( \left(\frac{\partial L(q^{s}_\mathrm{EL},q^{f}_\mathrm{EL},\dot{q}^{s}_\mathrm{EL},\dot{q}^{f}_\mathrm{EL})}{\partial q^{s}} - \frac{d}{d t} \frac{L(q^{s}_\mathrm{EL},q^{f}_\mathrm{EL},\dot{q}^{s}_\mathrm{EL},\dot{q}^{f}_\mathrm{EL})}{\partial \dot{q}^{s} } \right)(q^{s}_\mathrm{EL} - \hat q^{s})  + \mathcal{O}\left(\Delta T^{2l}\right) \bigg) d t\\ &+ \left( \frac{\partial L(q^{s}_\mathrm{EL},q^{f}_\mathrm{EL},\dot{q}^{s}_\mathrm{EL},\dot{q}^{f}_\mathrm{EL})}{\partial \dot{q}^{s} } (q^{s}_\mathrm{EL} - \hat q^{s}) \middle) \right|_0^{\Delta T}\\
 &+ \sum_{m=0}^{p-1}\int_{m\Delta t}^{(m+1)\Delta t} \bigg( \left(\frac{\partial L(q^{s}_\mathrm{EL},q^{f}_\mathrm{EL},\dot{q}^{s}_\mathrm{EL},\dot{q}^{f}_\mathrm{EL})}{\partial q^{f}} - \frac{d}{d t} \frac{L(q^{s}_\mathrm{EL},q^{f}_\mathrm{EL},\dot{q}^{s}_\mathrm{EL},\dot{q}^{f}_\mathrm{EL})}{\partial \dot{q}^{f} } \right)(q^{f}_\mathrm{EL} - \check q^{f})  + \mathcal{O}\left(\Delta t^{2l}\right) \bigg) d t\\ &+ \sum_{m=0}^{p-1}\left( \frac{\partial L(q^{s}_\mathrm{EL},q^{f}_\mathrm{EL},\dot{q}^{s}_\mathrm{EL},\dot{q}^{f}_\mathrm{EL})}{\partial \dot{q}^{f} } (q^{f}_\mathrm{EL} - \check q^{f}) \middle) \right|_{m\Delta t}^{(m+1)\Delta t}\\
  =&  \int_0^{\Delta T}  \mathcal{O}\left(\Delta T^{2l}\right)\, dt + \sum_{m=0}^{p-1}\int_{m \Delta t}^{(m+1)\Delta t}  \mathcal{O}\left(\Delta t^{2l}\right)\, dt.
	\end{align*}
 For the first equality we Taylor expanded the second integrand around $(q^{s}_\mathrm{EL},q^{f}_\mathrm{EL})$ and used the estimates for the polynomial interpolation (Lemma~\ref{lemma-poly}) where we have rewritten the integral involving fast derivatives as sum of integrals defined on the micro intervals. The second equality follows from applying integration by parts.
	The last equality follows since the boundary terms vanish because $\hat q^{s}(0) = q^{s}_\mathrm{EL}(0)$ and $\hat q^{s}(\Delta T) = q^{s}_\mathrm{EL}(\Delta T)$ as well as $\check q^{f}(m\Delta t) = q^{f}_\mathrm{EL}(m\Delta t)$ for $m=0,\ldots,p$ and since $(q^{s}_\mathrm{EL},q^{f}_\mathrm{EL})$ solves the Euler-Lagrange equation. Thus, we find 
 \begin{align}\label{eq:firstterm}
 \int_0^{\Delta T}& L(q^{s}_\mathrm{EL},q^{f}_\mathrm{EL},\dot{q}^{s}_\mathrm{EL},\dot{q}^{f}_\mathrm{EL}) \,d t - \int_0^{\Delta T} L(\hat{q}^{s},\check{q}^{f},\dot{\hat{q}}^{s},\dot{\check{q}}^{f}) \,d t \\
 =&  \mathcal{O}\left( \Delta T^{2l+1} \right) + \sum_{m=0}^{p-1} \mathcal{O}\left( (m+1)^{2l+1} - m^{2l+1}) \Delta t^{2l+1} \right) = \mathcal{O}\left( \Delta T^{2l+1} \right)\nonumber
 \end{align}
 since $\sum_{m=0}^{p-1} \left( (m+1)^{2l+1} - m^{2l+1})\right) = {p}^{2l+1}$ and ${p}\Delta t =\Delta T$.
	To bound the second term of \eqref{difference-expansion} we follow the arguments of \cite[Theorem 3.3]{hall2015spectral}. Since $(\tilde{q}^{s},\tilde{q}^{f})$ is the minimising element of $\mathcal{P}^l([0,\Delta T],Q^{s}) \times \mathcal{P}_{p}^l([0,\Delta T],Q^{f})$, we have
	\begin{align*}
 &\Delta t \sum_{m=0}^{p-1}\sum_{i=1}^r b_i L(\tilde{q}^{s}(m \Delta t + c_i \Delta t),\tilde{q}^{f}(m \Delta t + c_i \Delta t), \dot{\tilde{q}}^{s}(m \Delta t + c_i \Delta t),\dot{\tilde{q}}^{f}(m \Delta t + c_i \Delta t))\\
 &\leq \Delta t \sum_{m=0}^{p-1}\sum_{i=1}^{s} b_i L(\hat{q}^{s}(m \Delta t + c_i \Delta t),\check{q}^{f}(m \Delta t + c_i \Delta t), \dot{\hat{q}}^{s}(m \Delta t + c_i \Delta t),\dot{\check{q}}^{f}(m \Delta t + c_i \Delta t) ) \\
  &= \sum_{m=0}^{p-1}\int_{m\Delta t}^{(m+1)\Delta t} L(\hat q^{s}, \check q^{f} ,\dot{\hat q}^{s},\dot{\check q}^{f}) \, d t + \mathcal{O}(\Delta t^{u+1})\\
   &= \int_0^{\Delta T} L(\hat q^{s}, \check q^{f}, \dot{\hat q}^{s},\dot{\check q}^{f}) \, d t + {p} \,\mathcal{O}( \Delta t^{u+1})
	\end{align*}
	On the other hand, since $q_\mathrm{EL}$ minimises the continuous action, there holds
	\begin{align*}
 &\Delta t \sum_{m=0}^{p-1}\sum_{i=1}^r b_i L(\tilde{q}^{s}(m \Delta t + c_i \Delta t),\tilde{q}^{f}(m \Delta t + c_i \Delta t),  \dot{\tilde{q}}^{s}(m \Delta t + c_i \Delta t),\dot{\tilde{q}}^{f}(m \Delta t + c_i \Delta t))\\
 &=\int_0^{\Delta T} L(\tilde q^{s}, \tilde q^{f}, \dot{\tilde q}^{s},\dot{\tilde q}^{f}) \, d t + p \,\mathcal{O}( \Delta t^{u+1}) \geq \int_0^{\Delta T} L(q^{s}_\mathrm{EL},q^{f}_\mathrm{EL}, \dot{q}^{s}_\mathrm{EL},\dot{q}^{f}_\mathrm{EL}) \,d t + p \,\mathcal{O}( \Delta t^{u+1})  \\
 &=\int_0^{\Delta T} L(\hat q^{s},\check q^{f}, \dot{\hat q}^{s}, \dot{\check q}^{f}) \, d t + \mathcal{O}\left( \Delta T^{2l+1} \right) + p\,\mathcal{O}( \Delta t^{u+1}),
	\end{align*} 
	where the last line follows from \eqref{eq:firstterm}. Combining both inequalities we find
 \begin{align}\label{difference2}
& \int_0^{\Delta T} L(\hat{q}^{s},\check{q}^{f},\dot{\hat{q}}^{s},\dot{\check{q}}^{f})\,d t - \Delta t \sum_{m=0}^{p-1}\sum_{i=1}^r b_i L(\tilde{q}^{s}(m \Delta t + c_i \Delta t),\tilde{q}^{f}(m \Delta t + c_i \Delta t),  \dot{\tilde{q}}^{s}(m \Delta t + c_i \Delta t),\dot{\tilde{q}}^{f}(m \Delta t + c_i \Delta t))\nonumber\\
 &=\mathcal{O}\left( \Delta T^{2l+1} \right) + p \,\mathcal{O}( \Delta t^{u+1}).
 \end{align}
	Equations \eqref{eq:firstterm} and \eqref{difference2} together imply the desired result \eqref{eq:LELd}.
 \end{proof}
  By rewriting $p \,\mathcal{O}( \Delta t^{u+1}) = \frac{1}{p^u} \mathcal{O}( \Delta T^{u+1})$ we directly see that for $p=1$ the error estimate is identical to the error estimate for classical Galerkin variational integrators \cite{ober2021superconvergence} leading to superconvergence of the variational integrator, namely a convergence order of  $\min{(2l,u)}$. For $p>1$ the second error term in \eqref{eq:LELd} is decreased, nevertheless the order does not change for constant $p$. Thus, with Theorem~\ref{th:varerror} we obtain

\begin{theorem}[Superconvergence of multirate variational integrator]\label{th:super}
Let ${L}_d$ as given in \eqref{eq:hatLext} 
 be a multirate discretisation of the Lagrangian ${L}$ 
 based on polynomials of degree $l$ and a quadrature rule of degree $u$. Assume that all discrete and continuous critical curves minimise their respective actions. Then the corresponding multirate variational integrator \eqref{eq:DEL} is of order $\min{(2l,u)}$.
\end{theorem}

The different potential approximation in Theorems~\ref{th:midfmids}--\ref{th:trapftraps} can be generalised to
\begin{align}
Z_d(q_0^0,\ldots,q_0^p) &= \Delta t \sum_{m=0}^{p-1}\alpha Z(\gamma q_0^m + (1-\gamma) q_0^{m+1}) + (1-\alpha) Z((1-\gamma) q_0^m + \gamma q_0^{m+1} )\label{eq:Zapprox}
\end{align}
with $\alpha,\gamma \in [0,1]$ and $Z$ being fast or slow potential. Table~\ref{tab:coeff} gives an overview of the coefficients $\alpha$, $\gamma$ and the resulting quadrature rule with the corresponding coefficients $(b_i,c_i)^r_{i=1}$ that are used in \eqref{L_d}. All the underlying polynomials for the approximation of curves are of degree $1$. 
\begin{table}
		\caption{Coefficients of the quadrature rules}
		\label{tab:coeff}
\begin{tabular}{ll}
\hline\noalign{\smallskip}
$\alpha \in [0,1]$, $ \gamma_= \frac{1}{2}$ &  Gauss quadrature, $b_1=1, c_1=\frac{1}{2}$ (midpoint rule) \\
\noalign{\smallskip}\hline\noalign{\smallskip}
 $\alpha \in [0,1]$, $\gamma \in\{0, 1\}$ & trapezoidal rule \\[6pt]
$\quad \: \alpha=\frac{1}{2}$, $\gamma \in\{0, 1\}$ & Lobatto quadr., $b_1= b_2 = \frac{1}{2}, c_1=0, c_2=1$ \\[6pt]
$\quad \:  \alpha=1$, $\gamma =1$   & \multirow{2}{*}{{left rectangle rule, $b_1= 1, c_1=0 $}}  \\
$\quad \:   \alpha=0$, $\gamma =0$  & \\[6pt]
$\quad \:  \alpha=1$, $\gamma =0$  & \multirow{2}{*}{{right rectangle rule, $b_1= 0, c_1=1$}}   \\
$\quad \:  \alpha=0$, $\gamma =1$   &  \\
\noalign{\smallskip}\hline\noalign{\smallskip}
\end{tabular}
	\end{table}	 
The three approximations given in Theorems~\ref{th:midfmids}--\ref{th:trapftraps} can be constructed by the following choices of $\alpha_{\{W,V\}}$ and $\gamma_{\{W,V\}}$, where $(\alpha_W,\gamma_W)$ determine the approximation of the fast potential $W$ and $(\alpha_V,\gamma_V)$ determine the approximation of the slow potential $V$: 
\begin{itemize}
\item midpoint fast, midpoint slow: $\alpha_W, \alpha_V \in [0,1], \gamma_W = \gamma_V = \frac{1}{2}$ 
 \item midpoint fast, trapezoidal slow:  $\alpha_W, \alpha_V \in [0,1], \gamma_W = \frac{1}{2},\gamma_V \in\{0, 1\}$  
  \item trapezoidal fast, trapezoidal slow: $\alpha_W,\alpha_V\in [0,1],\gamma_W,\gamma_V \in \{0, 1\}$.
\end{itemize}

With this notation, we can summarise the local error resulting by applying the three presented variational integrators in Theorems~\ref{th:midfmids}--\ref{th:trapftraps} to simulate the multirate Lagrangian system.

\begin{lemma}[Approximation error]\label{th:approx}
Consider the multirate Lagrangian \eqref{eq:multirateL}. Using the approximation rules given in Theorems~\ref{th:midfmids}--\ref{th:trapftraps} leads to the following maximal local error of the corresponding discrete Lagrangian $L_d$ and therefore of the discrete evolution map $(q_k, p_k) \mapsto (q_{k+1}, p_{k+1} )$ of the corresponding $(p,q)$--scheme 
\begin{align} 
&  \|   L^E_d (q_0,q_1, \Delta T) - L_d(q_0, q_{1}, \Delta T)  \|  = \\[5pt]
 &\left\{\begin{array}{cll} \mathcal{O}(\Delta T^3) + p\mathcal{O}(\Delta t^3) & \text{for midpoint and trapezoidal with }{\alpha_{V,W}} =\frac{1}{2} \text{ combinations} \\ 
 &\text{in slow and fast}  \\[1ex] 
\frac{1}{p}\mathcal{O}(\Delta T^2) & \text{for all the other combinations mentioned in Table \ref{tab:coeff}} \end{array} \right.
\end{align} 
leading to a second order integrator for midpoint and trapezoidal approximation (${\alpha_{V,W}} =\frac{1}{2}$) and  to an first order integrator otherwise.
\end{lemma}

\begin{proof}
This directly follows from Theorem~\ref{thm-superh} and Theorem~\ref{th:super} since $l=1$ for all cases, and $u=2$ for midpoint and trapezoidal (${\alpha_{V,W}} =\frac{1}{2}$) and $u=1$ for all other combinations mentioned in Table \ref{tab:coeff}.
\end{proof}


\begin{remark}
In summary, the theorem states that as soon as the part of the action integral of either the fast or the slow potential is approximated by any quadrature other than midpoint or trapezoidal with ${\alpha_{V,W}}=\frac{1}{2}$, the order of convergence is reduced from 2 to 1.
\end{remark}

Note that, given the multirate variational scheme is stable (cf.~Section~\ref{sec:stability}), the local errors given in Theorem~\ref{th:approx} determine the convergence order for finer macro and micro discretizations, i.e.~$\Delta T\to 0$ and $\Delta t \to 0$. Considering the local error as function of $\Delta t$, $\Delta T$ and $p$, we have
\begin{equation}\label{eq:lerror}
e(\Delta t,\Delta T,p) = p \mathcal{O}(\Delta t^q) + \mathcal{O}(\Delta T ^q) = c p \Delta t^q + C \Delta T ^q
\end{equation}
for some constants $c,C >0$, which we assume w.l.o.g.~to be equal to $1$. 
Let us first consider the case in which both, micro and macro grids will get refined by the same {amount}, i.e.~$\Delta T, \Delta t \to 0$ and $p=\text{constant}$. Then, in \eqref{eq:lerror} the dominant term is $\mathcal{O}(\Delta T ^q)$ and we expect a convergence order of $q-1$. To see this more concretely, we compute
\begin{equation}
\frac{e\left(\frac{\Delta t}{2},\frac{\Delta T}{2},p\right)}{e(\Delta t,\Delta T,p)} =  \frac{p \left( \frac{\Delta t}{2} \right)^q + \left( \frac{\Delta T}{2} \right)^q}{p \Delta t^q + \Delta T^q} = \frac{1}{2^q} \frac{p \Delta t^q + \Delta T^q}{p \Delta t^q + \Delta T^q} = \frac{1}{2^q}.
\end{equation}
The local error rate can be conveniently described by the following logarithmic representation
\begin{equation}
\frac{\log_2{\left(e(\Delta t,\Delta T,p)\right)}-\log_2{\left(e\left(\frac{\Delta t}{2},\frac{\Delta T}{2},p\right)\right)}}{\log_2{\Delta t} - \log_2{\frac{\Delta t}{2}}}= \log_2{\frac{e(\Delta t,\Delta T,p)}{e\left(\frac{\Delta t}{2},\frac{\Delta T}{2},p\right)}}=\log_2{2^q} = q.
\end{equation}
Thus, for refining macro and micro grid at the same time, the rate is $q$ providing an error of order $q-1$ as expected.
Let us now consider the case, in which only the micro grid is refined to improve the approximation of the fast variables, whereas the macro grid stays the same, i.e.~$\Delta T=\text{constant}$, $\Delta t \to 0$ and $p\to\infty$. In this case, only the first term in \eqref{eq:lerror}, $p \mathcal{O}(\Delta t^q)$, will decrease while the second term $\mathcal{O}(\Delta T ^q)$ stays constant for constant $\Delta T$ and will dominate eventually for $\Delta t\to 0$ and $p\to \infty$. Thus, for increasing $p$ we expect the error rate to get smaller. To see this, we compute
\begin{equation}
\frac{e\left(\frac{\Delta t}{2},\Delta T,2p\right)}{e(\Delta t,\Delta T,p)} =  \frac{2p \left( \frac{\Delta t}{2} \right)^q + \left( {\Delta T} \right)^q}{p \Delta t^q + \Delta T^q} = \frac{1}{2^q} \frac{2p \Delta t^q + 2^q\Delta T^q}{p \Delta t^q + \Delta T^q} \stackrel{\Delta T = p \Delta t}{=} \frac{1}{2^q} \frac{2p  + 2^q p^q}{p  + p^q}
\end{equation}
and thus
\begin{equation}
\frac{\log_2{\left(e(\Delta t,\Delta T,p)\right)}-\log_2{\left(e\left(\frac{\Delta t}{2},\Delta T,2p\right)\right)}}{\log_2{\Delta t} - \log_2{\frac{\Delta t}{2}}} = q-\log_2{\frac{2p  + 2^q p^q}{p  + p^q}}= q-1-\log_2{\frac{1  + (2p)^{q-1}}{1  + p^{q-1}}}.
\end{equation}
As $p\to \infty$ we have $\lim_{p\to\infty} \log_2{\frac{1  + (2p)^{q-1}}{1  + p^{q-1}}} = \lim_{p\to\infty} \log_2{\frac{1/p^{q-1}  + 2^{q-1}}{1/p^{q-1}  + 1}}=$ \\ $\log_2{2^{q-1}} = q-1$ and the rate becomes $0$, i.e., there is no improvement of the error for decreasing $\Delta t$ since the error in the slow discretized dynamics is not getting reduced if $\Delta T$ stays constant.

\section{Linear stability analysis}\label{sec:stability}

A general linear stability analysis of the multirate variational integrator is not carried out, since the $[4 \times 4]$-propagation matrix depends on $p$ and thus, finding a compact representation such that the propagation matrix does not have to be recalculated for each specific value of $p$ is difficult. However, the influence of the linearly interpolated slow variable on the linear stability is investigated. Therefore, combinations where the variables are not split but linearly interpolated with function evaluations on the micro grid are considered. Linearly interpolating the variables between the macro nodes together with potential evaluation on the micro grid can be interpreted as a composition rule on the basis of the trapezoidal respectively the midpoint rule. 
As a test system, a simple harmonic oscillator with eigenfrequency $\omega_s$ is used. The corresponding Lagrangian reads
\begin{equation}
L(q, \dot{q} ) = \frac{1}{2}  \dot{q}^2 - \omega^2_s q^2.
\label{harOsci_Stabili_slow}
\end{equation}
The results of the propagation matrix analysis are given in Lemma \ref{lemma:stability_slow}. {When using the trapezoidal rule, the linear stability region increases by increasing the number $p$ of microsteps, until a supremum is reached. When using the midpoint rule, the unconditional stability of the single rate midpoint rule, $p=1$, is lost by increasing $p$. For $p = 2$ linear stability is guaranteed by the choice $\omega_s \Delta T  < 4$. By increasing $p$ further the stability region reduces until a infimum is reached.}

\begin{lemma} \label{lemma:stability_slow}
Linearly interpolating the non-split variable between the macro nodes together with
the trapezoidal rule to approximate the action of the harmonic oscillator (with eigenfrequency $\omega_s$) on the micro grid, cf.~\eqref{trapSlow1} and \eqref{trapSlow2} in Theorem \ref{th:midftraps} respectively \eqref{trapSlow3} and \eqref{trapSlow4} in Theorem \ref{th:trapftraps} with $q^s=q$ and $V=V(q^s)$, the algorithm is linearly stable, if
\begin{align*}
\omega^2_s \Delta T^2 < \frac{12 p^2}{p^2 +2} 
\end{align*}
is fulfilled. While for $p=1$ the upper bound is 4, {the supremum of the upper bound is $12$}. 

Linearly interpolating the non-split variable between the macro nodes together with
the midpoint rule to approximate the action of the harmonic oscillator (with eigenfrequency $\omega_s$) on the micro grid, cf.~\eqref{midSlow1} and \eqref{midSlow2} in Theorem \ref{th:midfmids} with $q^s=q$ and $V=V(q^s)$, the algorithm is linearly stable, if
\begin{align*}
{\omega^2_s \Delta T^2  < \frac{12 p^2}{p^2 - 1 } }
\end{align*}
is fulfilled. For $p=1$, the scheme is unconditionally stable as this yields the midpoint rule. {For $p \rightarrow \infty$, the infimum of the upper bound is $12$.}
\label{re:stabil_slow}
\end{lemma}
\begin{proof}
{For the proof, the equations in \eqref{trapSlow1} and \eqref{trapSlow2} resp.~in \eqref{midSlow1} and \eqref{midSlow2} are applied to the harmonic oscillator as given in \eqref{harOsci_Stabili_slow} and reformulated to
\begin{equation}
\begin{pmatrix}
q_{k+1} \\ p_{k+1} 
\end{pmatrix} = P^{\text{trap,s}/\text{mid,s}} \begin{pmatrix}
q_k \\ p_k
\end{pmatrix}
\end{equation}
where $P^{\text{trap,s}/\text{mid,s}}$ is the so called propagation matrix. The restriction of the eigenvalues of the propagation matrix to be distinct and of modulus one guarantees linear stability.} 
The corresponding propagation matrices $P^{\text{trap,s}}$ resp.~$P^{\text{mid,s}}$ read
\begin{align*}
P^{\text{trap,s}} &= \begin{pmatrix}
- \frac{\omega^2_s \Delta t^2 ( \frac{p^2}{3} + ({ \alpha_V} - \frac{1}{2}) p + \frac{1}{6})  - 1}{ \frac{\omega^2_s \Delta t^2 (p^2 - 1) }{6}    + 1} & 
\frac{\Delta t p}{\frac{\omega^2_s  \Delta t^2 (p^2 - 1)}{6}    + 1 }  \\
\frac{\Delta t \omega^2_s ( \omega^2_s ( \frac{p^2}{3}   + ({\alpha_V} - \frac{1}{2}) p + \frac{1}{6}) \Delta t^2 - 1) ( \frac{p}{2} - {\alpha_V} + \frac{1}{2})}{
\frac{(\omega^2_s (p^2 - 1) \Delta t^2)}{6} + 1}  - \Delta t \omega^2_s ({\alpha_V} + \frac{p}{2} - \frac{1}{2} ) &
1 - \frac{\Delta t^2 p \omega^2_s (\frac{p}{2} - {\alpha_V} + \frac{1}{2})}{\frac{ (\omega^2_s (p^2 - 1) \Delta t^2)}{6} + 1}
\end{pmatrix} \\[4pt]
{P^{\text{mid,s}}} &= \begin{pmatrix}
 \frac{ 12 -4 \Delta t^2 p^2  \omega^2_s +  \Delta t^2 \omega_s^2 }{12 + 2 \Delta t^2 p^2 \omega_s^2 + \Delta t^2 \omega_s^2} &  \frac{12 p \Delta t}{12 + 2 \Delta t^2 p^2 \omega_s^2 + \Delta t^2 \omega_s^2} \\[4pt]
\frac{-12 \Delta t \omega_s^2 p + \Delta t^3 p^3 \omega_s^4 - \Delta t^3 \omega_s^4 p)}{12 + 2 \Delta t^2 p^2 \omega_s^2 + \Delta t^2 \omega_s^2}  &  \frac{12-4 \Delta t^2 p^2 \omega_s^2 + \Delta t^2 \omega_s^2}{12 + 2 \Delta t^2 p^2 \omega_s^2 + \Delta t^2 \omega_s^2}
\end{pmatrix}.
\end{align*}
As $\det(P^{\text{trap,s}}) = \det(P^{\text{mid,s}}) = 1$, the schemes are stable, if $| \text{tr}(P^{\text{trap,s}}) | < 2$ resp.~$| \text{tr}(P^{\text{mid,s}}) | < 2$ with 
\begin{align*}
\text{tr}(P^{\text{trap,s}}) &= - \frac{2 (2 \omega^2_s \Delta t^2 p^2 + \omega^2_s \Delta t^2 - 6)}{\omega^2_s \Delta t^2 p^2 - \omega^2_s \Delta t^2 + 6} \\
{\text{tr}(P^{\text{mid,s}}) } &=   2 \frac{ 12 -4 \Delta t^2 p^2  \omega_s^2 +  \Delta t^2 \omega_s^2 }{12 + 2 \Delta t^2 p^2 \omega_s^2 + \Delta t^2 \omega_s^2}
\end{align*}
yielding the inequalities in Lemma \ref{lemma:stability_slow}. 
\end{proof}

\begin{figure}[h]
		\centering
		\includegraphics[width=.49\textwidth]{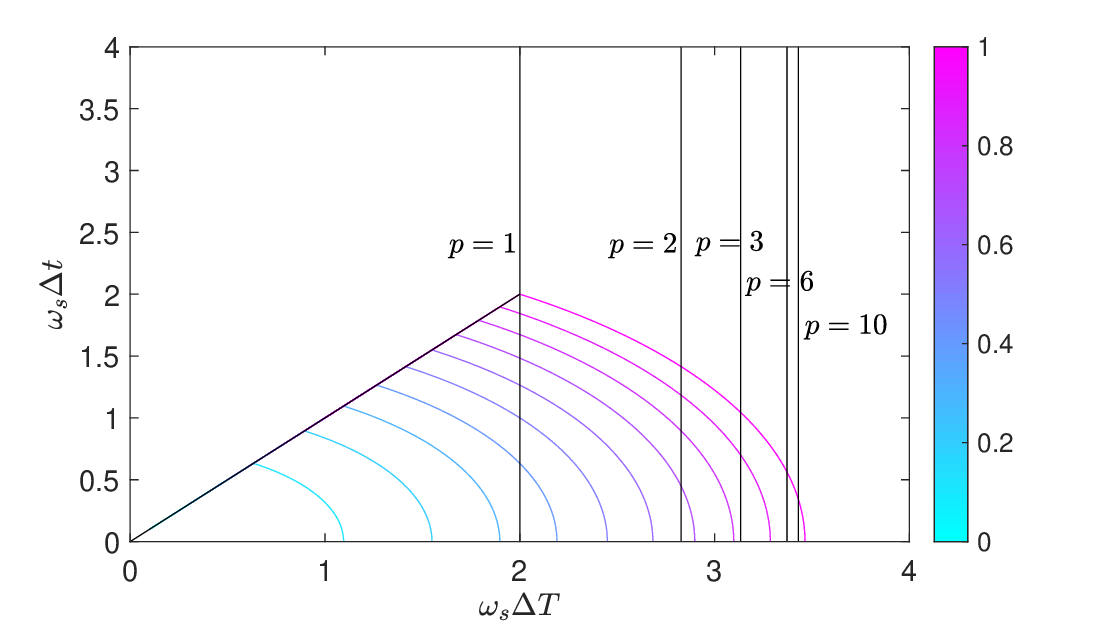}
		\includegraphics[width=.49\textwidth]{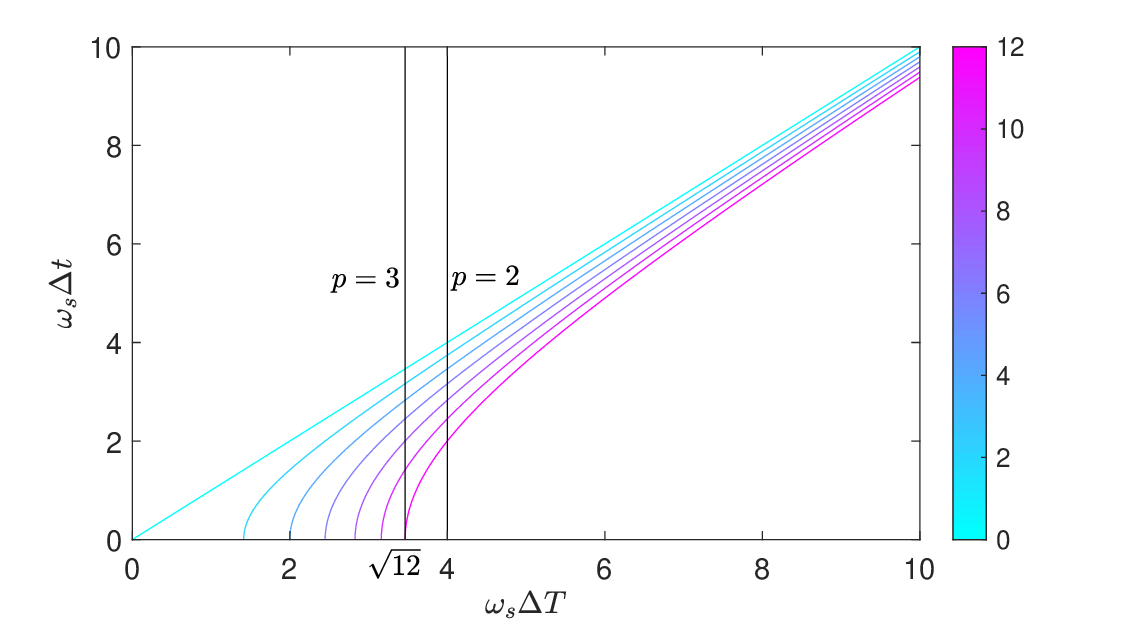}
		\caption{Linear stability regions depending on the ratio $p = \frac{\omega_s \Delta
		 T}{\omega_s \Delta t}$ for using the trapezoidal rule (left plot) resp.~the midpoint rule (right plot) to approximate the action on the micro grid.}
		\label{fig:linear_stability}
	\end{figure}

{Substituting $p = \frac{\Delta T}{\Delta t}$ in the inequalities given in Lemma \ref{re:stabil_slow} yields $\frac{1}{12} \omega_s^2 \Delta T^2 + \frac{1}{6} \omega^2_s \Delta t^2 < 1$ (maximum of coloured level sets in Figure \ref{fig:linear_stability}, left), when using the trapezoidal rule, and $\omega^2_s \Delta T^2 - \omega^2_s \Delta t < 12$ (maximum of coloured level sets in Figure \ref{fig:linear_stability}, right), when using the midpoint rule to approximate the integral of the potential on the micro grid. Figure \ref{fig:linear_stability} shows the stability regions when approximating via the trapezoidal rule (left plot) resp.~via the midpoint rule (right plot). In the left plot, the black inclined line represents the ratio $p =1$. The stability region for $p=1$ lies to the left of the vertical line through its intersection with the level set of value 1. Similarly, the other vertical lines bound the stability regions for higher values of $p$. It can be seen that the stability region increases by increasing $p$ until the supremum of $\sqrt{12}$ is reached. The right plot shows that the stability region decreases when increasing $p$. The turquoise inclined straight line represents $p=1$, i.e.~$\Delta T = \Delta t$, being of level $0$ and therefore fulfilling the linear stability condition for all choices of $\omega_s \Delta T$. By increasing $p$, the unconditional linear stability is lost and for $p=2$, $\omega_s \Delta T$ can be at most 4 and for the choice $p \geq 3$, $\omega_s \Delta T$ can be at most $\sqrt{12}$. Deductively, for higher values of $p$, $ p \gtrsim 10$, the choice whether to use the trapezoidal rule or the midpoint rule to approximate the action on the micro grid has negligible impact on the linear stability.}

The stability analysis results in Lemma \ref{re:stabil_slow} also provide an orientation for the choice of the number of intermediate steps $p$ for the general approach of multirate variational integrators (i.e.~with split variables and split potentials, see Theorem \ref{th:midfmids}, \ref{th:midftraps} resp.~\ref{th:trapftraps}).

\section{Numerical examples}
\label{sec:numerical_Examples}

Some of the theoretical results proven in the preceding Sections are demonstrated numerically. 
The good energy behaviour of the multirate variational integrators due to the conservation of the symplectic structure is shown. The conservation of momentum maps proven via Noether's theorem in Section \ref{sec:conservation_properties} is demonstrated. Further, the convergence order of the integrators proven in Theorem \ref{th:approx} is verified numerically. Computing time investigations show the performance of the integrators regarding runtime savings.

In the following we refer to the integrator defined in Theorem \ref{th:midfmids} as midpoint-midpoint scheme, to the integrator defined in Theorem \ref{th:midftraps} as trapezoidal-midpoint scheme and to the integrator defined in Theorem \ref{th:trapftraps} as trapezoidal-trapezoidal scheme. 
We use the left rectangle rule for the slow potential ($\{ \alpha_V, \gamma_V \} \in \{0,1\}, \gamma_V=\alpha_V)$ in the trapezoidal-midpoint scheme and for both the slow and the fast potential ($\{ \alpha_V , \gamma_V \} \in \{0,1\}, \gamma_V=\alpha_V$ and $\{ \alpha_W , \gamma_W \} \in \{0,1\}, \gamma_W=\alpha_W)$ in the trapezoidal-trapezoidal scheme.
Two dynamical systems are considered, the famous Fermi-Pasta-Ulam problem and the so called spring ring example.

\begin{figure}[h]
			\centering
			\includegraphics[width=1.0\textwidth]{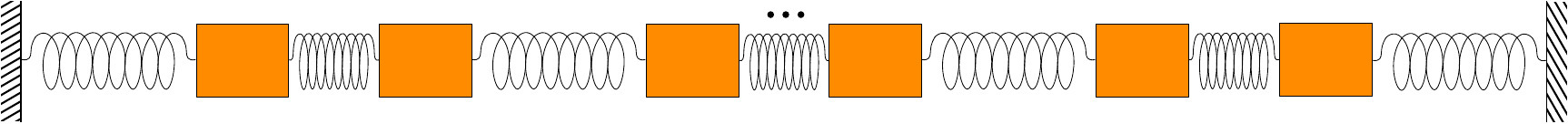}
			\caption{Fermi-Pasta-Ulam problem}
			\label{fig:FPU}
	\end{figure}

The Fermi-Pasta-Ulam problem (FPU) was originally introduced by \citep{fpu55}. We consider a modification here, introduced by \citep{galgani1992problem}, consisting of \(2\ell\) masses linked with alternating soft and stiff springs as shown in Fig.~\ref{fig:FPU}.
The slow variable $q^s_i$, $i=1,\ldots,l$, represents a scaled displacement of the centre of the $i$th stiff spring and the fast variables $q^f_i$, $i=1,\ldots,l$, a scaled expansion (or compression) of the $i$th stiff spring, cf.~\citep{hairer2006geometric}. 
	For \(2\ell\) masses, the kinetic energy and the slow and fast potential read
	\begin{align*}
	T(\dot{q}) &= \frac{1}{2} \dot{q}^T \cdot M \cdot \dot{q} \\
	V(q) &= \frac{1}{4}\left[(q^{s}_{1}-q^{f}_{1})^{4}+\sum^{\ell-1}_{i=1}(q^{s}_{i+1}-q^{f}_{i+1}-q^{s}_{i}-q^{f}_{i})^{4}+(q^{s}_{l}+q^{f}_{l})^{4}\right]\\
	W(q^f) &= \frac{\omega^2}{2}\sum^{\ell}_{i=1}(q^{f}_{i})^{2}
	\end{align*}
where \(M \in \mathbb{R}^{2\ell\times 2\ell}\) is the mass matrix. For the simulations six masses are included, thus $l=3$. All masses are set equal to one. The stiffness of the stiff springs is \(\omega^2 = 2500\). The initial positions are \(q^s_0 = \left[1,0,0\right]^T \text{, } q_0^{f,0} = \left[\frac{1}{|\omega|},0,0\right]^T\) and the initial 	
velocities are \(v^s_0 = \left[1,0,0\right]^T \text{, } v^{f,0}_0 = \left[1,0,0\right]^T\). These values are taken from \citep{hairer2006geometric}. In all simulations of the FPU we set the tolerance of the Newton-Raphson iteration to \( 10^{-9}\).

	\begin{figure}[h]
	
		\centering
		\includegraphics[width=0.5\textwidth]{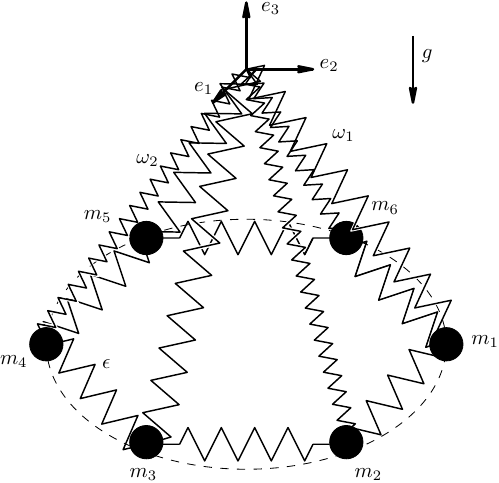}
		\caption{Spring ring example system with slow and fast configuration}
		\label{fig:spring_ring_sketch}
	\end{figure}
The spring ring (SR) is illustrated in Fig.~\ref{fig:spring_ring_sketch}.
	It consists of $2{ l}$ masses arranged in a circle, which are connected by soft springs of stiffness \(\epsilon\).
	Further, the masses are connected to the point of origin by alternating soft and stiff springs with stiffness \(\omega_1\) and \(\omega_2\) respectively.
The masses $m_1, m_3, m_5, \ldots , m_{ 2l-1}$ are connected to the origin by the soft springs with stiffness $\omega_1$ and their positions are assigned to the slow variables, while the masses $m_2 , m_4 , m_6 , \ldots , m_{{ 2l}}$ are connected to the origin by the stiff springs with stiffness $\omega_2$ and their positions are assigned to the fast variables.
	The kinetic and potential energies for $2{ l} $ masses are given by
	\begin{align*}
	T(\dot{q}^s,\dot{q}^f) 	&= \frac{1}{2} \dot{q}^T \cdot M \cdot \dot{q}\\
	V(q^s,q^f)		&= \frac{{\omega_1}}{2} ({q^s})^T   \cdot  q^s + \frac{\epsilon}{4} \left(\sum_{i=1}^{2 l - 1} \left\lbrace \|q_{i+1} - q_i\|^4 \right\rbrace + \|q_1 - q_{ 2l } \|^4 \right)+ {q}^T   \cdot  M   \cdot  g  \\
	W(q^f)		&= \frac{{\omega_2}}{2} ({q^f})^T  \cdot q^f
	\end{align*}
	with the gravitational acceleration vector $g$, pointing in the negative $e_3$ direction, $\vert g \vert$ $=9.81$, $q_i \in \mathbb{R}^3$ being the position vector of the $i$-th mass and $M\in \mathbb{R}^{6N \times 6N}$ being the mass matrix. 
For the simulations, we set $2{ l}=6$ and $m_i = 2$, $i= 1, \ldots, 6$. The stiffnesses of the springs are \(\epsilon = 5, \ \omega_1 = 2 \), and \(\omega_2 = 4000\). 
	The masses are arranged in a circle with radius \(r=2\). 
The midpoint of the circle is positioned under the origin in a distance of \(h = 2\).
The position of mass $i$ then reads $q_{i+1} = [r\sin(i  \pi/N ), -r\cos(i  \pi/N),-h]^T\) with $i=0,\ldots,2N-1$. \\
{The position of mass $i$ then reads $q_{i} = [r\sin((i-1)  \pi/N ), -r\cos((i-1)  \pi/N),-h]^T\) with $i=1,\ldots,2l$.}
	Further, to have an initial stretch in the springs, the masses \(m_3,m_5\) and \(m_2,m_4\) are repositioned. 
	\begin{alignat*}{2}
	{q_3}_0 &= {q_3}  + [0.3,0.3,0]^T     &\qquad {q_2}_0 &= {q_2}  + [0.2,-0.2,0]^T\\
	{q_5}_0 &= {q_5}  + [0.2,-0.3,-0.3]^T &       {q_4}_0 &= {q_4}  + [-0.3,0.4,0]^T\\
	\end{alignat*}
	The initial velocities are given in the following equations.
	\begin{alignat*}{2}
	{v_1}_0 & = \phantom{-}5\frac{q_1-q_2}{\| q_1 - q_2\|} &\qquad {v_2}_0 & = -30\frac{q_1-q_2}{\| q_1 - q_2\|}  \\
	{v_3}_0 & = -5\frac{q_3 - q_6}{\|q_3 - q_6\|}		   &       {v_4}_0 & = [50,40,-10]^T  \\
	{v_5}_0 & = [0,0,0]^T 								   &       {v_6}_0 & = [50, 40, 10]^T
	\end{alignat*}
	For all calculations we set the tolerance of the Newton-Raphson iteration to \(10^{-8}\).

\subsection{Conservation properties}
\label{sub:energies_and_preserved_quantities}
In this section, the structure preserving properties of the variational multirate integrators are presented. Further, it is shown that increasing $p$, with $\Delta T$ fixed, improves the simulations results.

\begin{figure}[h]
	
		\centering
		\includegraphics[width=.32\textwidth]{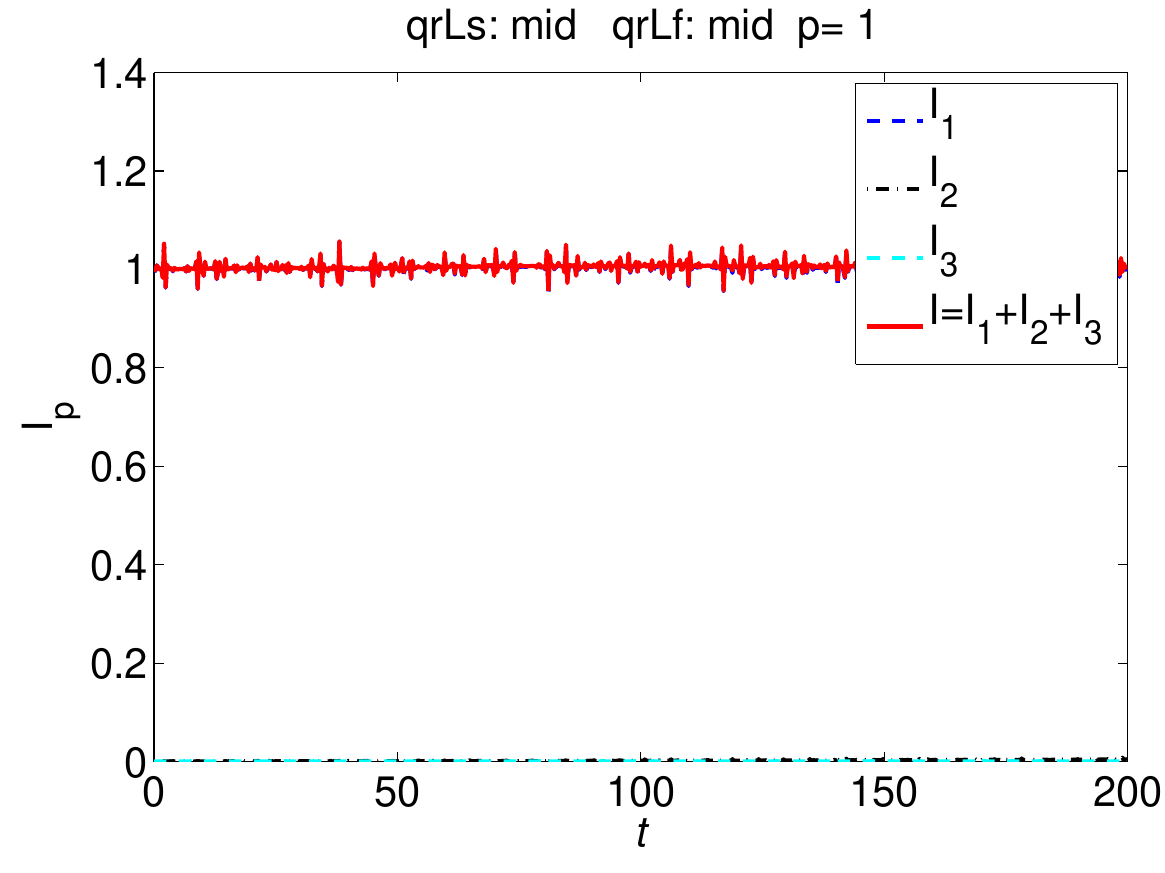}
		\includegraphics[width=.32\textwidth]{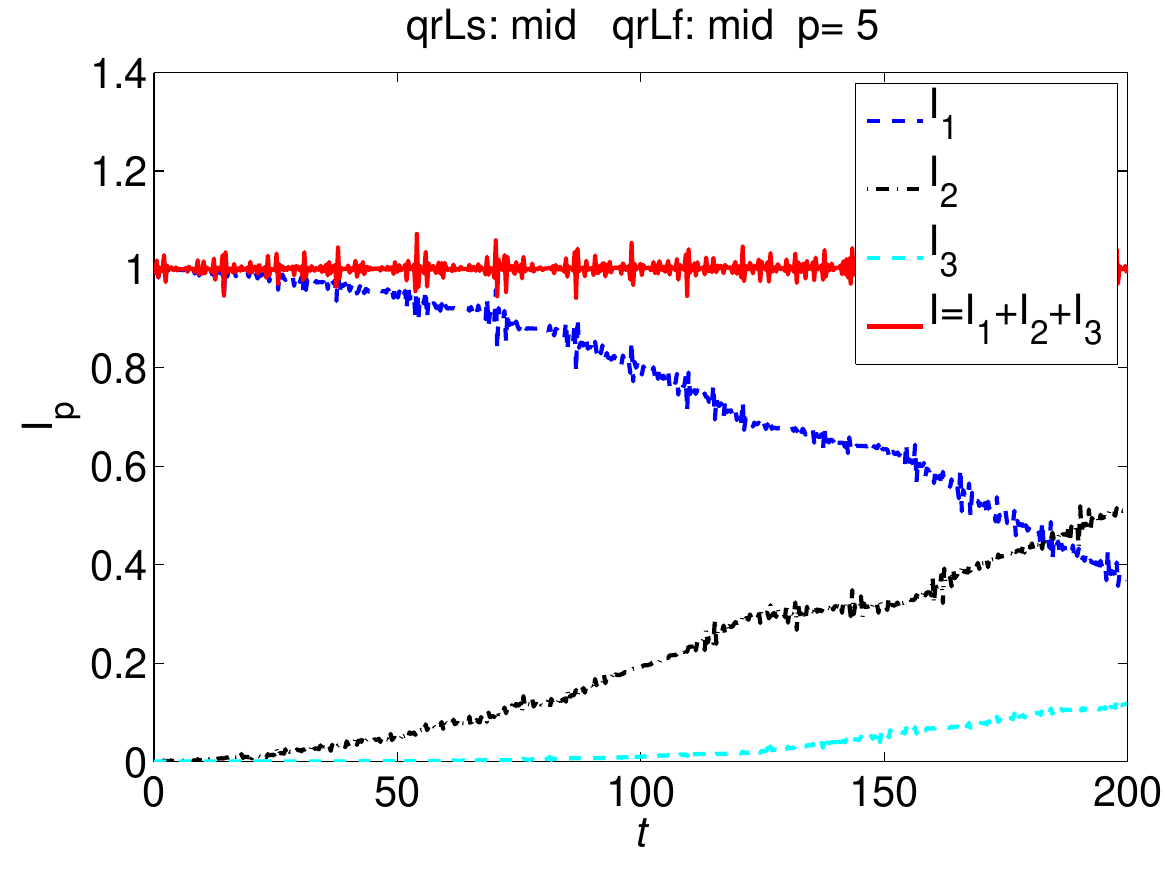}
		\includegraphics[width=.32\textwidth]{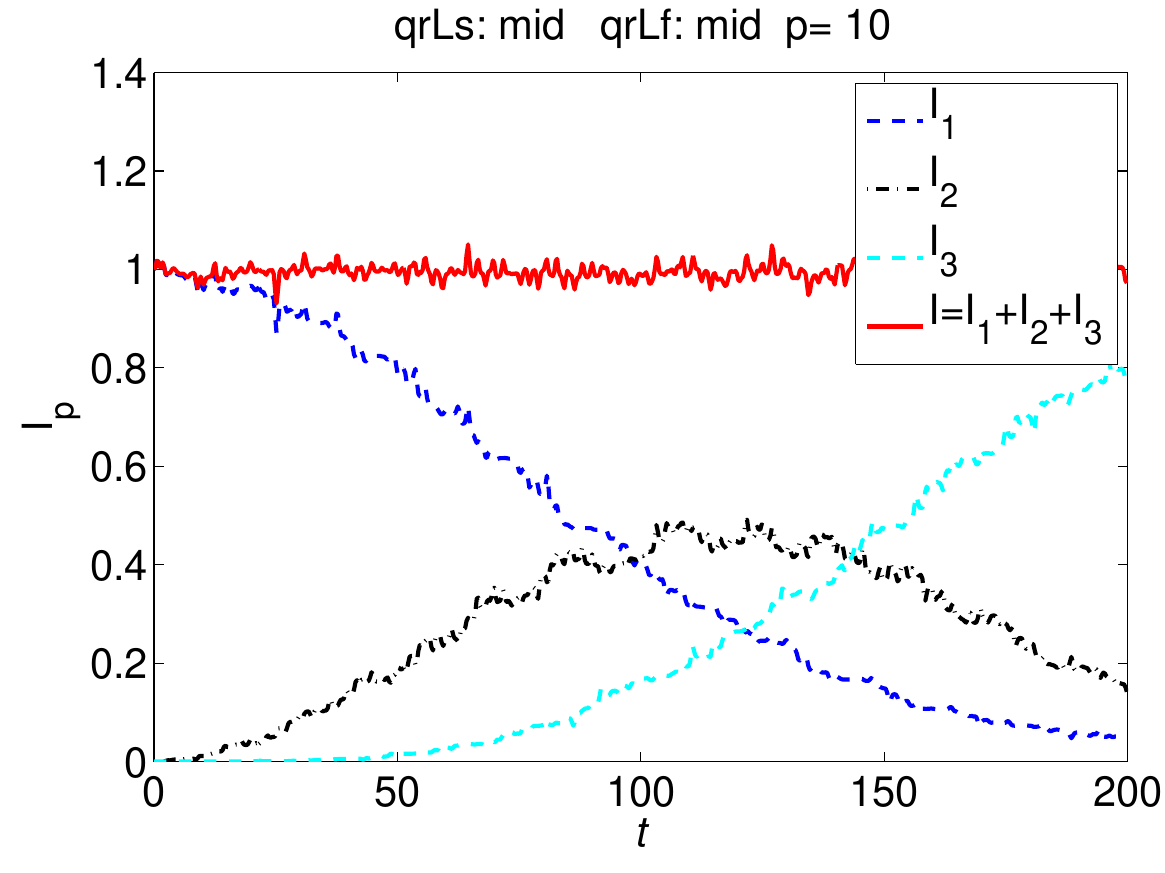}
		\caption{FPU: total energy in the stiff springs and sum of all stiff energies for \(p=1\) (left), \(p=5\) (centre), and \(p=10\) (right) with \(\Delta T = 0.3 \) for midpoint-midpoint scheme}
		\label{fig:fpu invariance imim}
	\end{figure}

\begin{figure}[h]

	\centering
	\includegraphics[width=.32\textwidth]{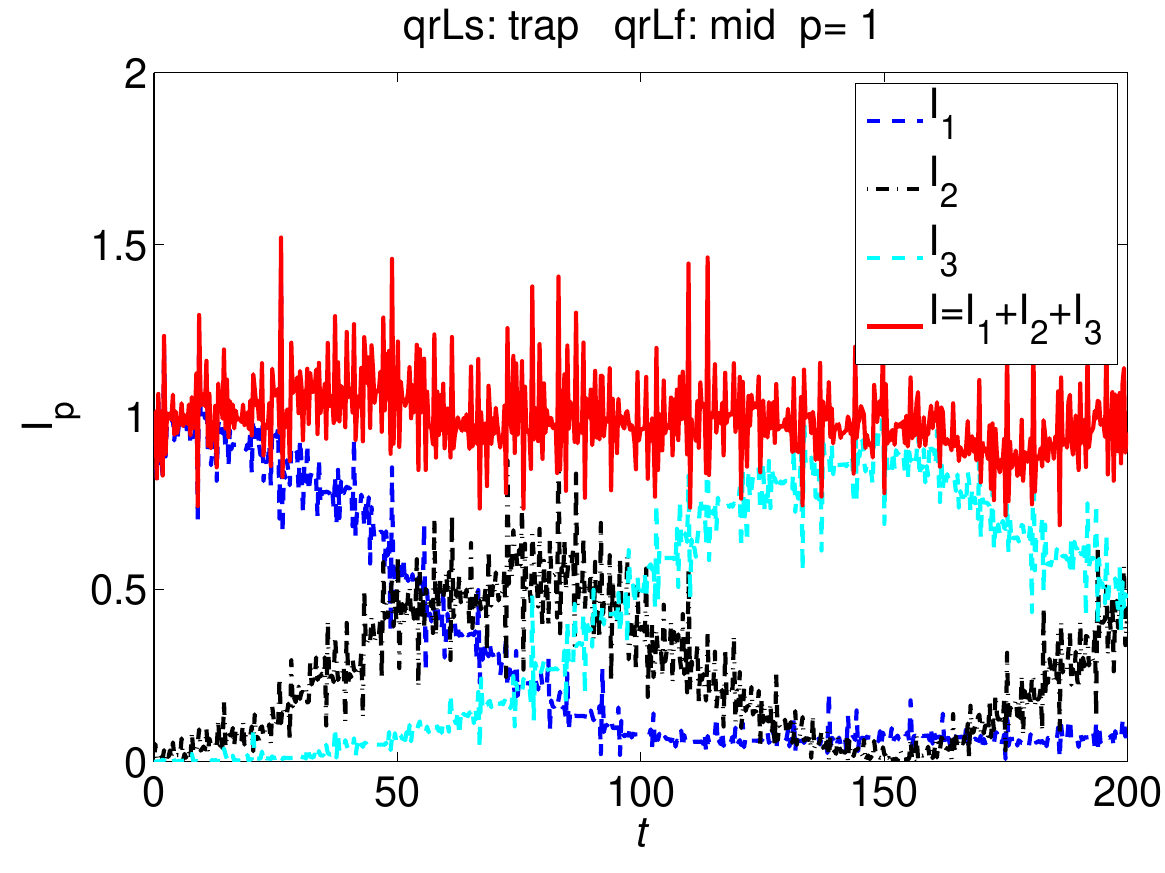}
	\includegraphics[width=.32\textwidth]{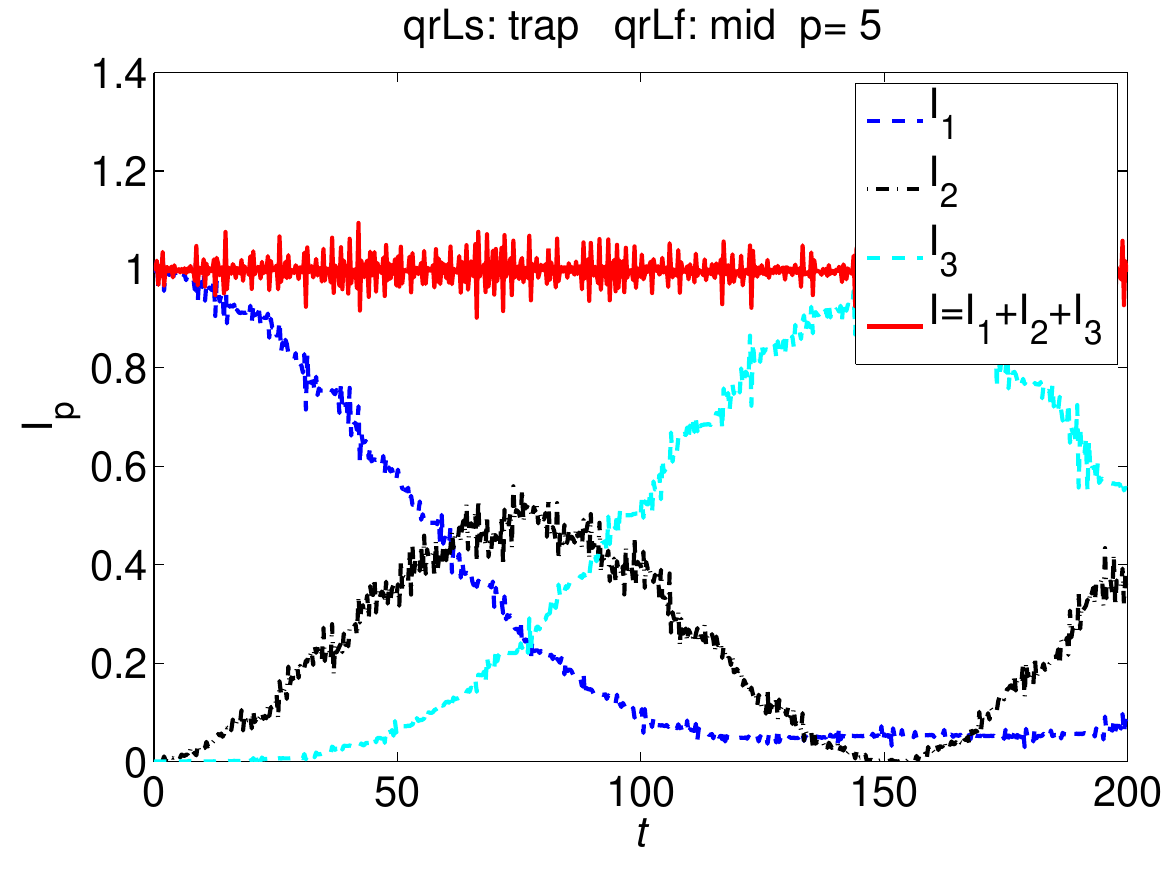}
	\includegraphics[width=.32\textwidth]{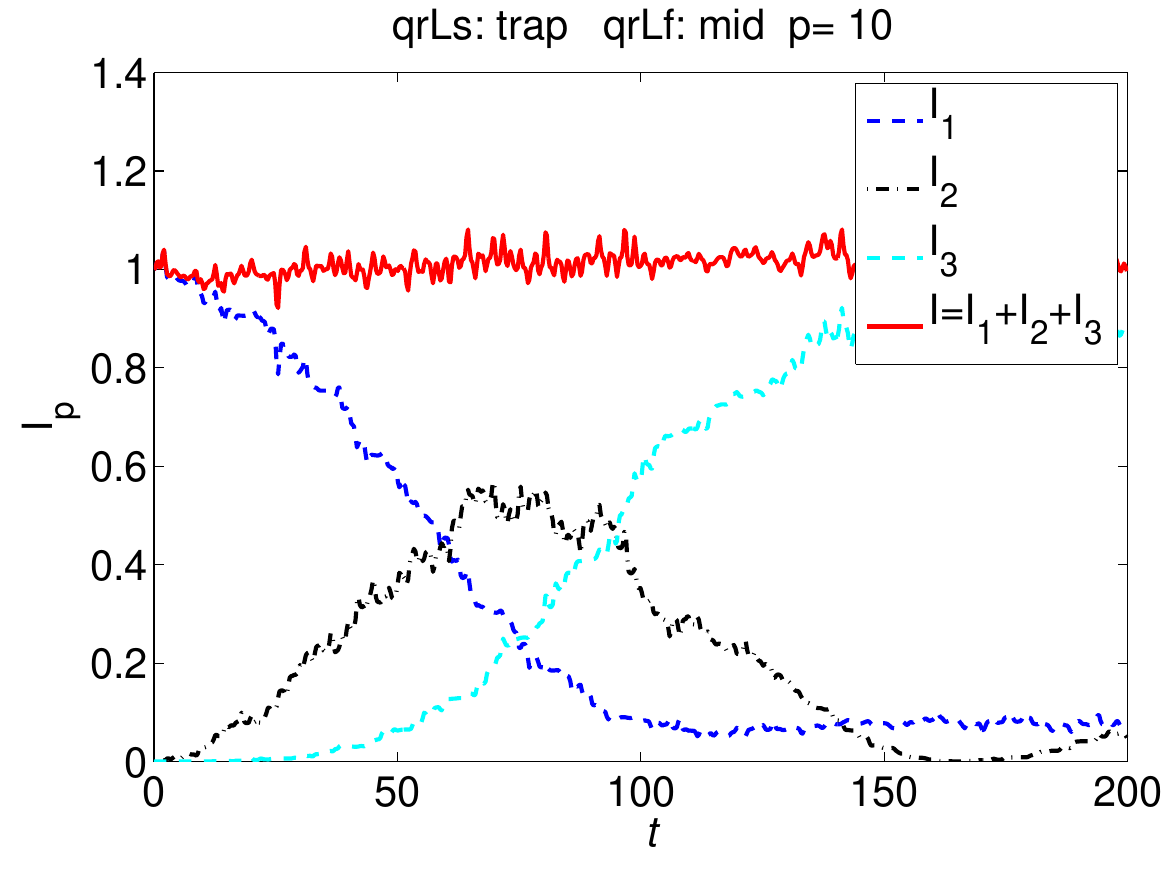}
	\caption{FPU: total energy in the stiff springs and sum of all stiff energies for \(p=1\) (left), \(p=5\) (centre), and \(p=10\) (right) with \(\Delta T = 0.3 \) for trapezoidal-midpoint scheme}
	\label{fig:fpu invariance exim}
\end{figure}

\begin{figure}[h]

	\centering
	\includegraphics[width=.32\textwidth]{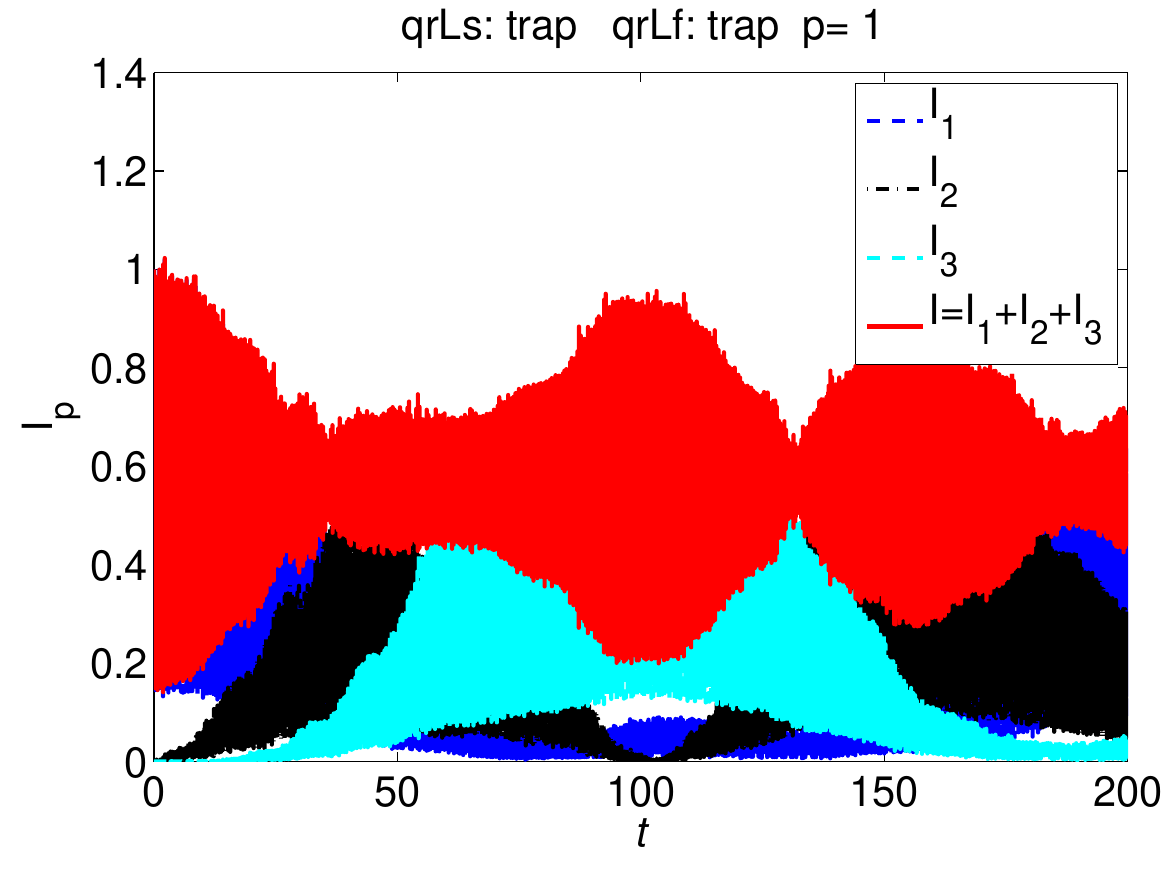}
	\includegraphics[width=.32\textwidth]{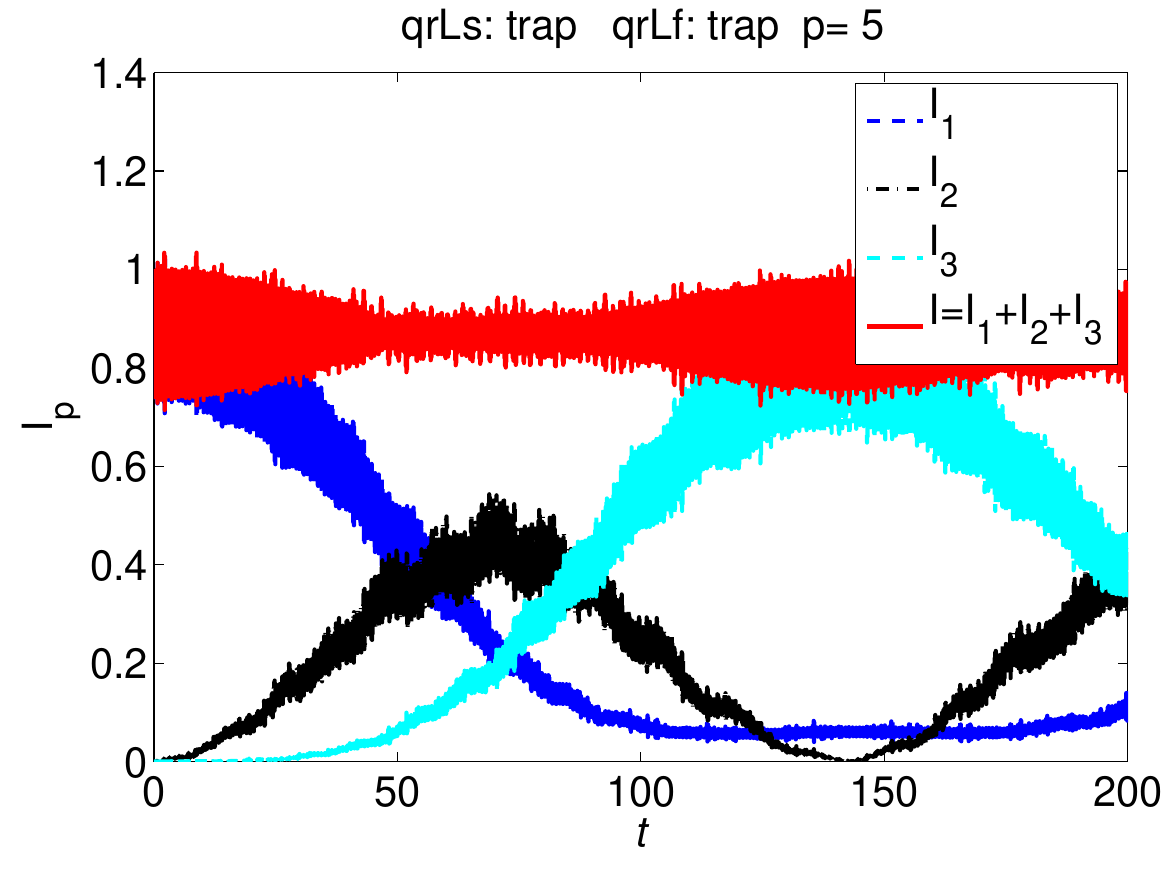}
	\includegraphics[width=.32\textwidth]{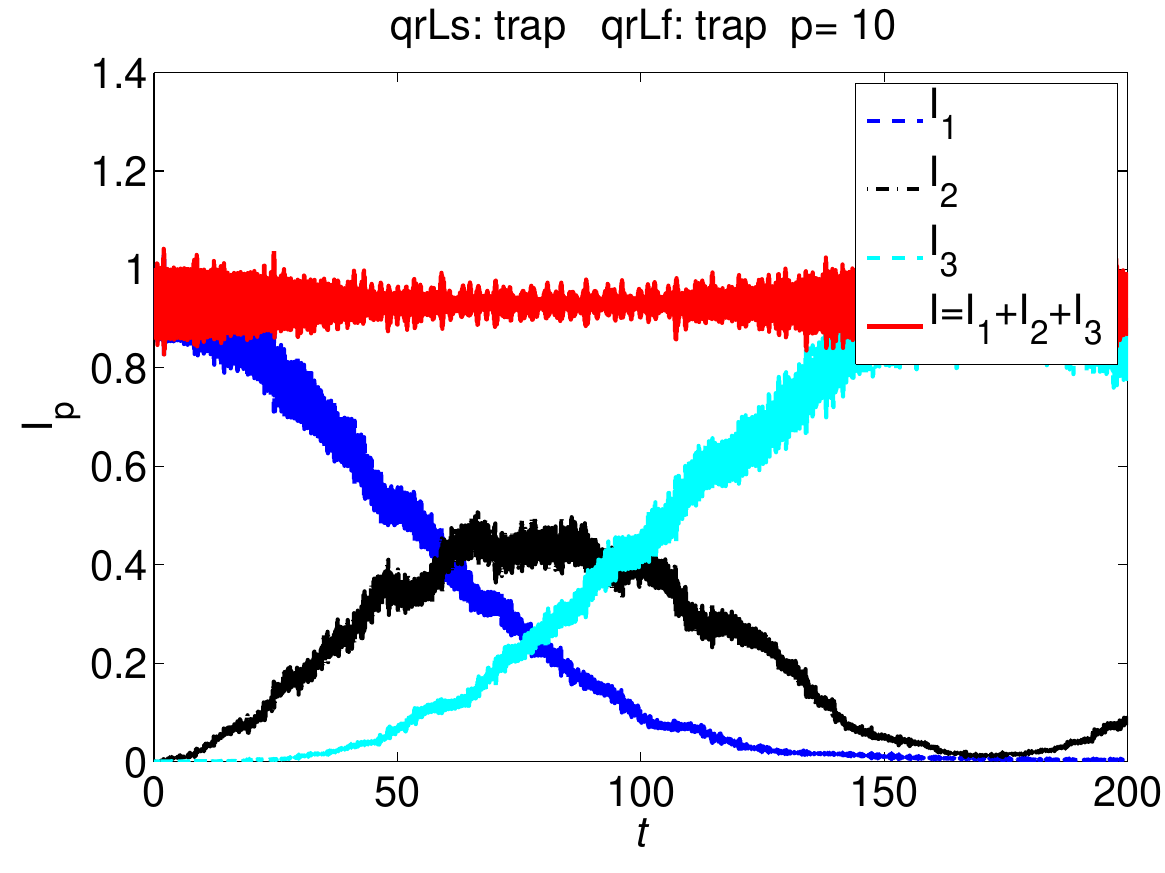}
	\caption{FPU: total energy in the stiff springs and sum of all stiff energies for \(p=1\) (left), \(p=5\) (centre), and \(p=10\) (right) with \(\Delta T = 0.03 \) for trapezoidal-trapezoidal scheme}
	\label{fig:fpu invariance exex}
\end{figure}

The FPU is simulated to an end time \(t_{N} = 200\) with \(\Delta T = 0.3\) for \(p\in\{1,5,10\}\). 
 The energy exchange is illustrated in the Fig.~\ref{fig:fpu invariance imim} -- \ref{fig:fpu invariance exex}.  
The energy $I_j$ stored in the $j$-th stiff spring, $j=1,2,3$, and the total stiff energy $I=\sum_j I_j$ is depicted.
One can see that the energy stored in the first stiff spring translates over to the second stiff spring and to the third stiff spring.	
While energy is exchanged between the stiff springs, 
the sum of energies over all stiff springs is an adiabatic invariant, oscillating around $I(0)=1$.
This illustrates that the integrators have good energy behaviour, with no numerical dissipation.
Furthermore, there are differences visible in the different quadrature schemes.
For the midpoint-midpoint scheme, the amplitude of the oscillations in the sum of all energies is not changing with \(p\).
However, the rate of the energy exchange is increasing with increasing \(p\).
This behaviour is also observed for the single rate midpoint-midpoint scheme when the time step size is decreased, see e.g.~\citep{hairer2006geometric}.
For the trapezoidal-midpoint scheme in Fig.~\ref{fig:fpu invariance exim} the time, at which the energy in the second spring (black line) is zero, is always around \( t =150 \) for all three values of \(p\).
In Fig.~\ref{fig:fpu invariance exex} the time at which the energy of the second spring is zero, is changing enormously with \(p\), from \( t =100 \) for \(p=1\) to approximately \( t=170 \) for \(p=10\). Note that the time step $\Delta T =0.03$ used for the simulations in Fig.~\ref{fig:fpu invariance exex} is near the instability region of the single rate integrator $p=1$ (left plot) as linear stability requires $\Delta T \omega < 2$, thus $\Delta T < 0.04$ here.
For both schemes, the trapezoidal-midpoint scheme and the trapezoidal-trapezoidal scheme, it is striking that the amplitude of the oscillations in the sum of the total oscillatory energies is decreasing with increasing \(p\) and therefore decreasing \(\Delta t\).
	
\begin{figure}[h]
	\centering
	\includegraphics[width=0.32\textwidth]{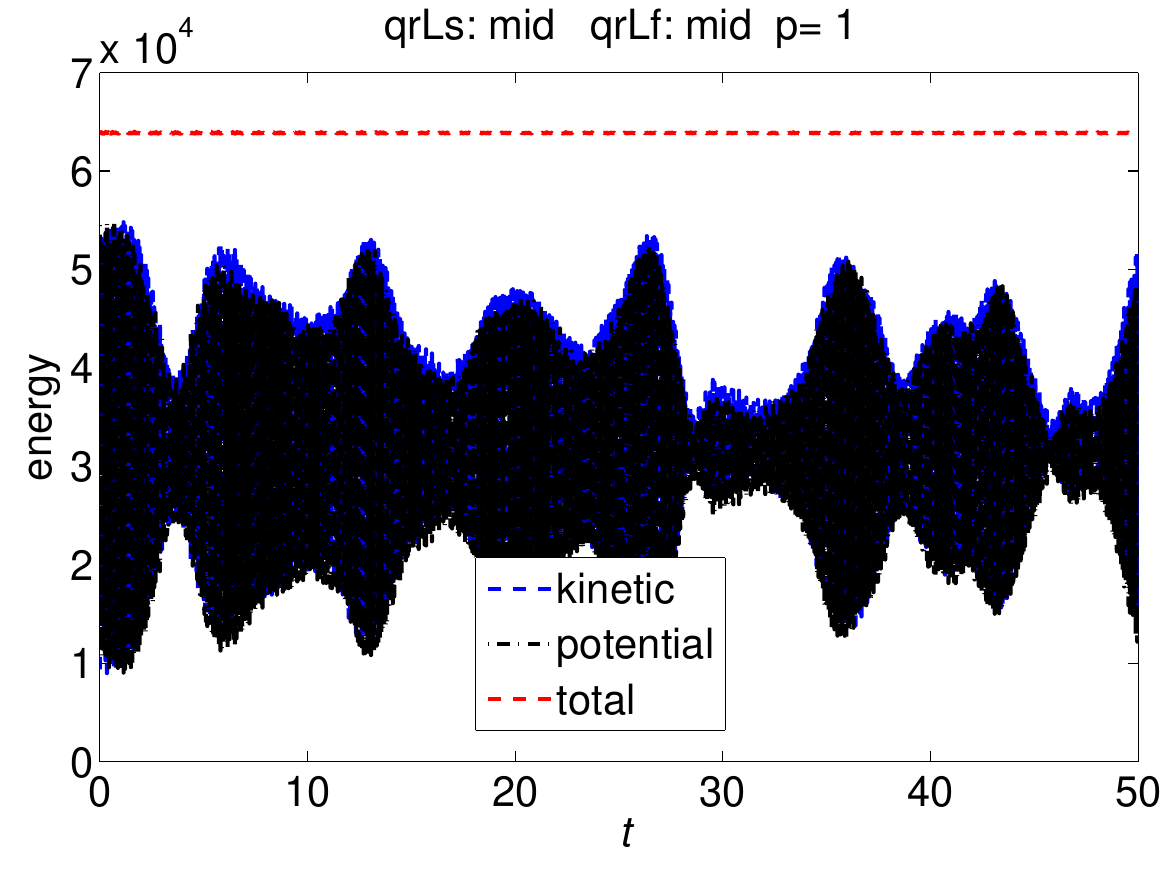}
	\includegraphics[width=0.32\textwidth]{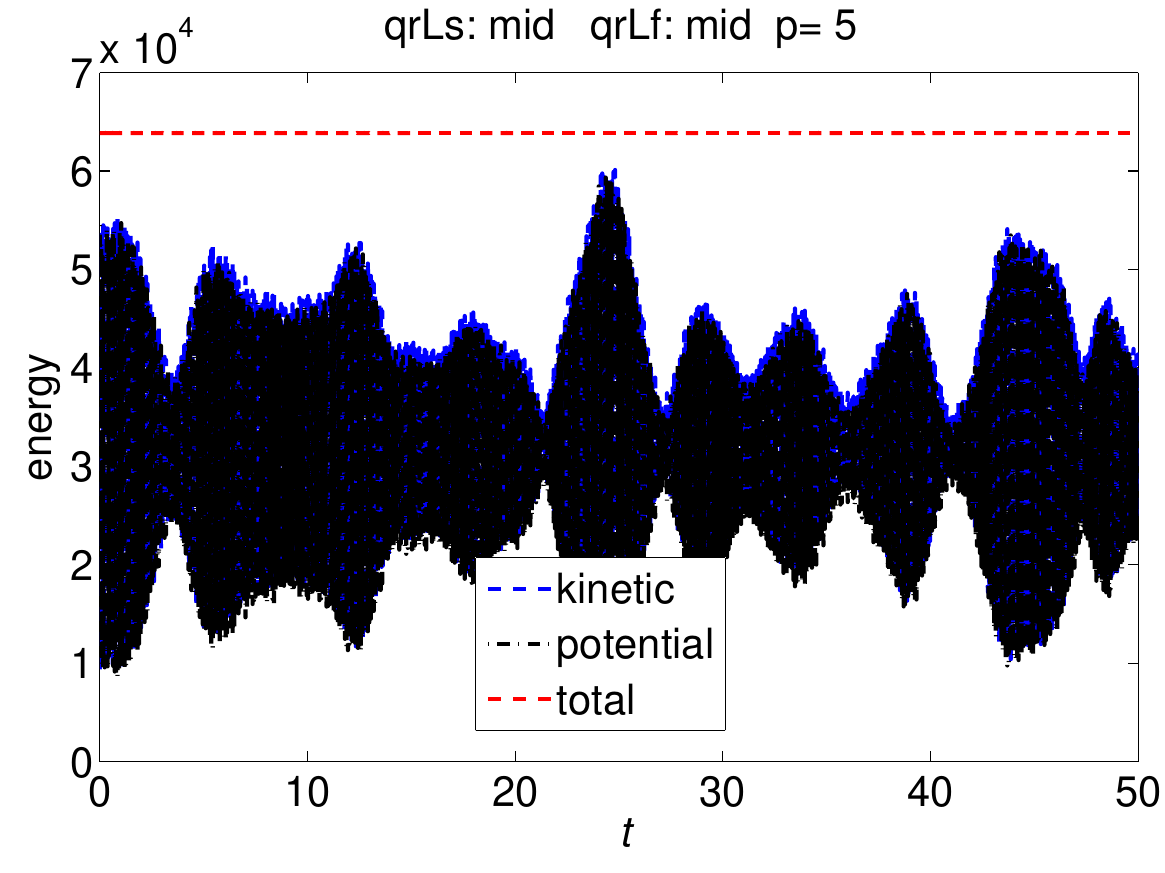}
	\includegraphics[width=0.32\textwidth]{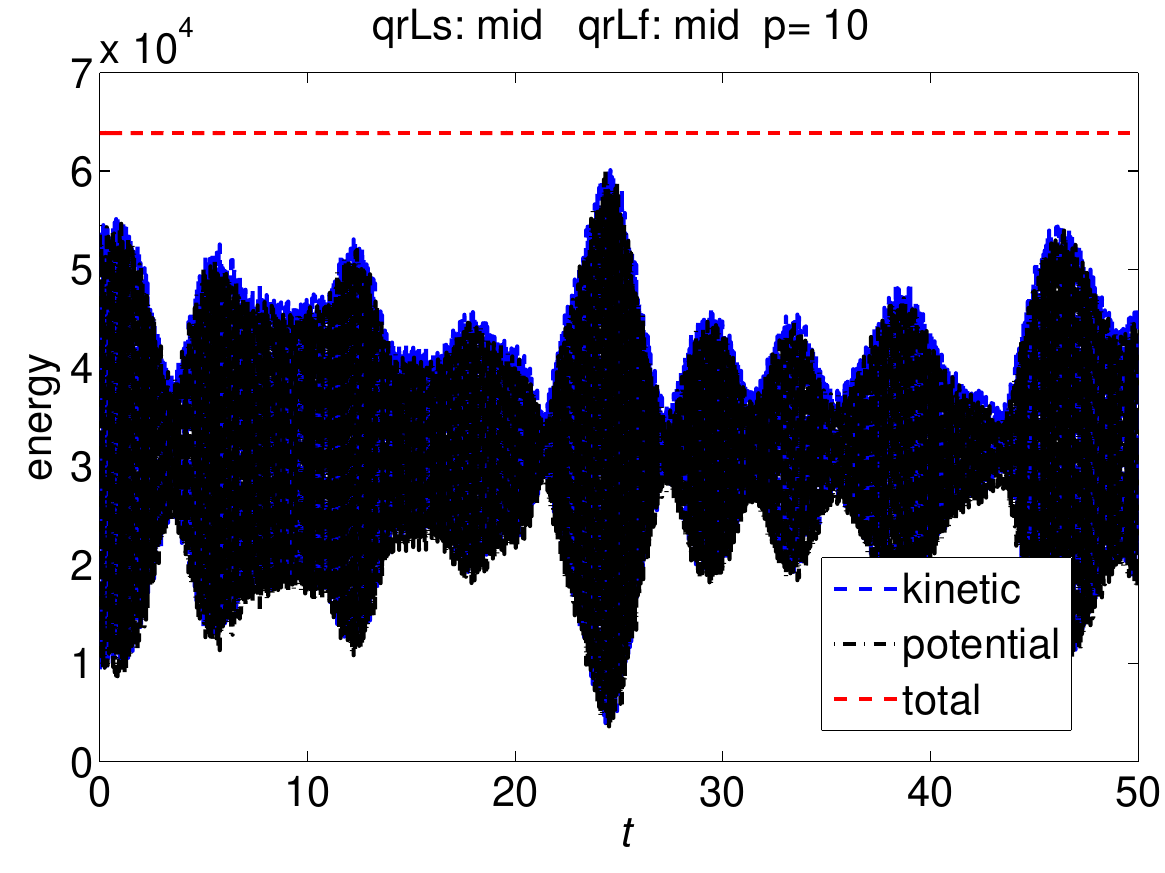}
	\caption{SR: evolution of energy for \(p=1\) (left), \(p=5\) (centre) and \(p=10\) (right) with \(\Delta T =0.01 \) for midpoint-midpoint scheme}
	\label{fig:sr energy imim}
\end{figure}

\begin{figure}[h]
	\centering\includegraphics[width=0.32\textwidth]{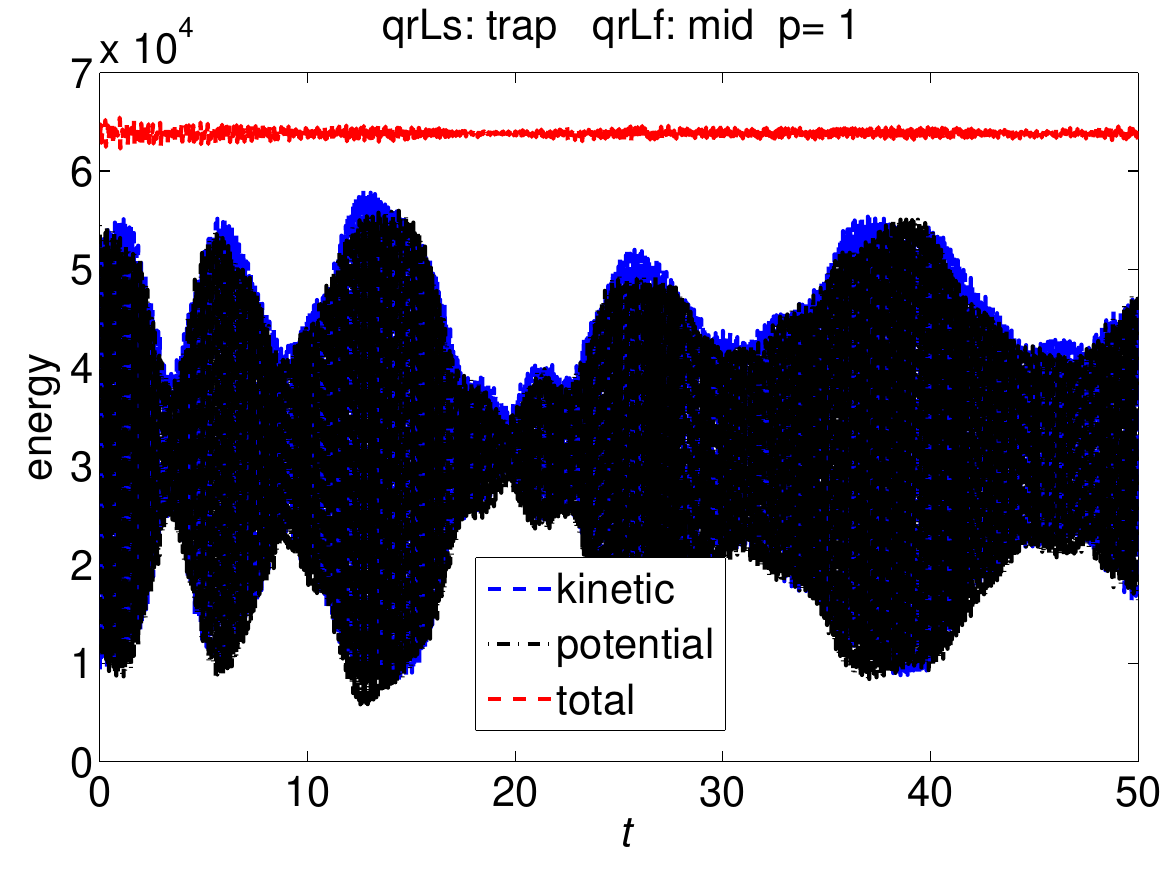}
	\includegraphics[width=0.32\textwidth]{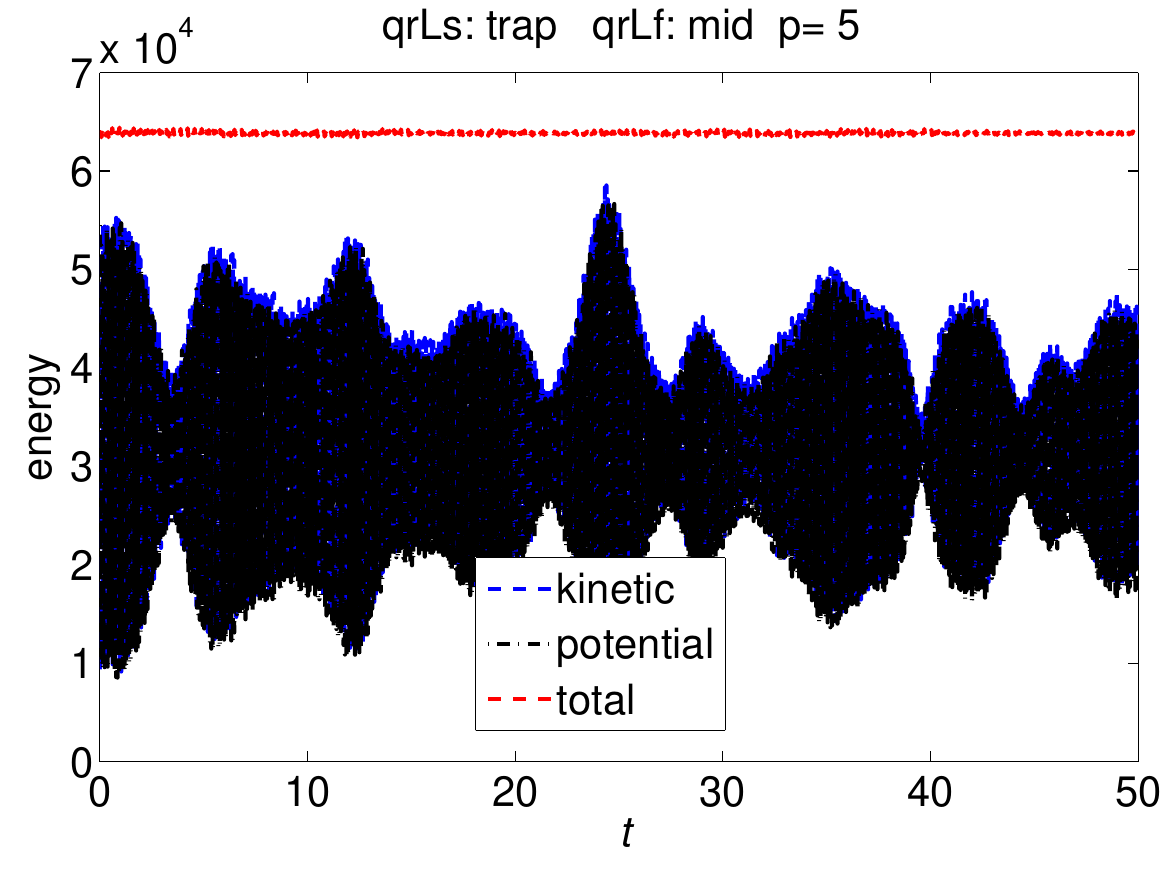}
	\includegraphics[width=0.32\textwidth]{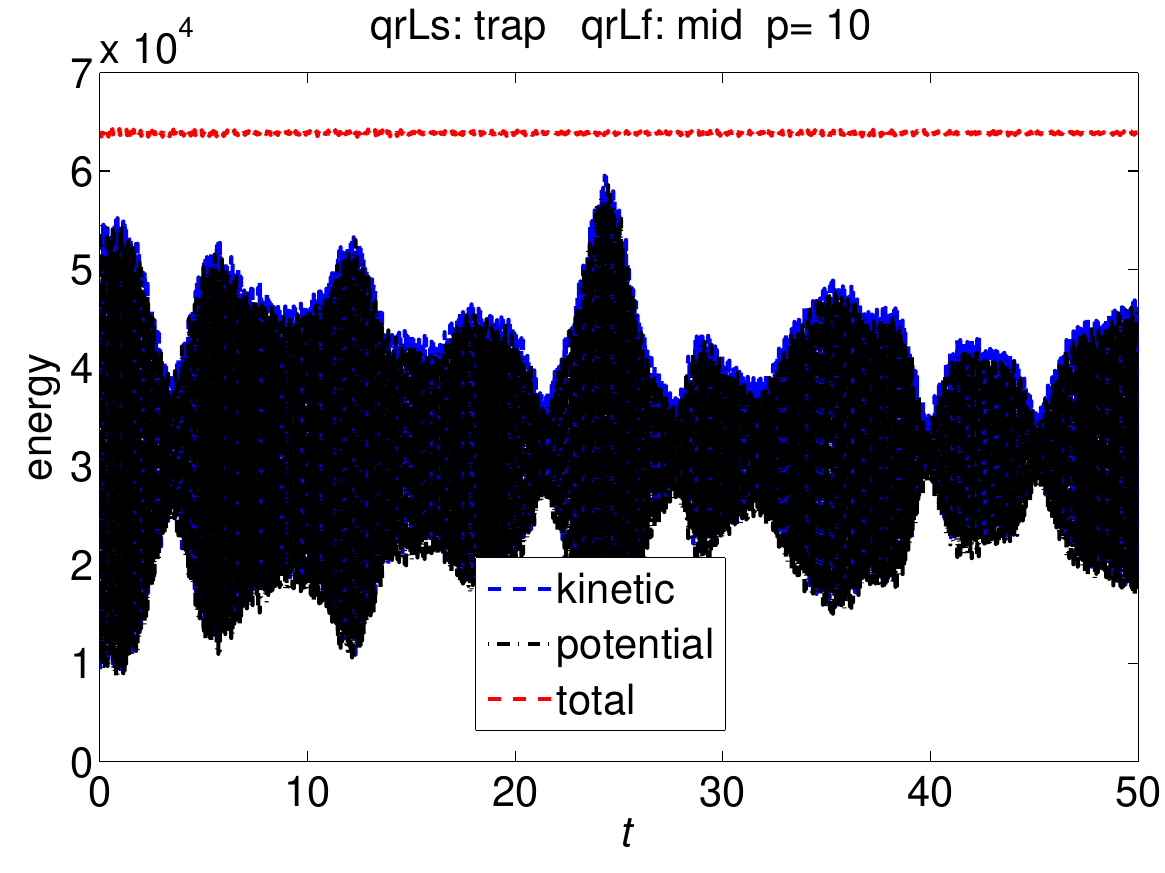}
	\caption{SR: evolution of energy for \(p=1\) (left), \(p=5\) (centre) and \(p=10\) (right) with \(\Delta T =0.01 \) for trapezoidal-midpoint scheme}
	\label{fig:sr energy exim}
\end{figure}

\begin{figure}[h]
	\centering\includegraphics[width=0.32\textwidth]{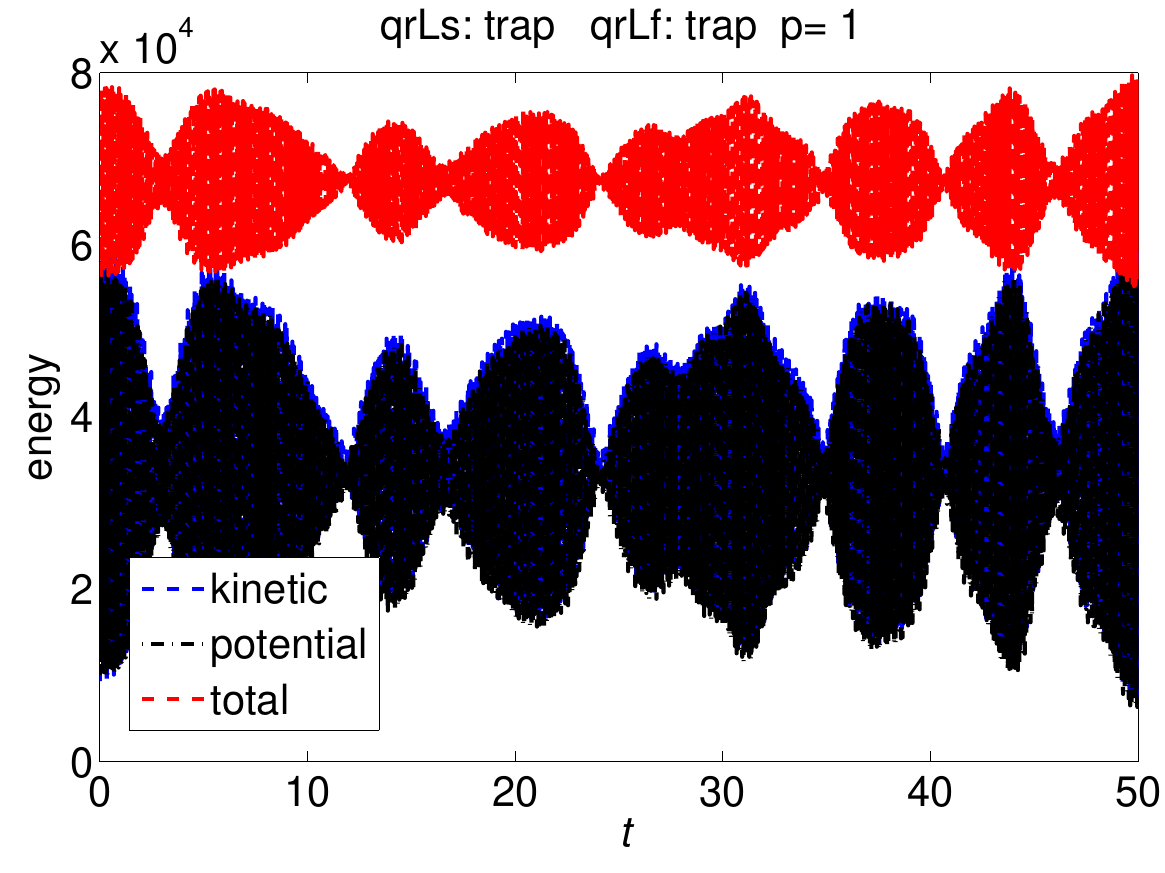}
	\includegraphics[width=0.32\textwidth]{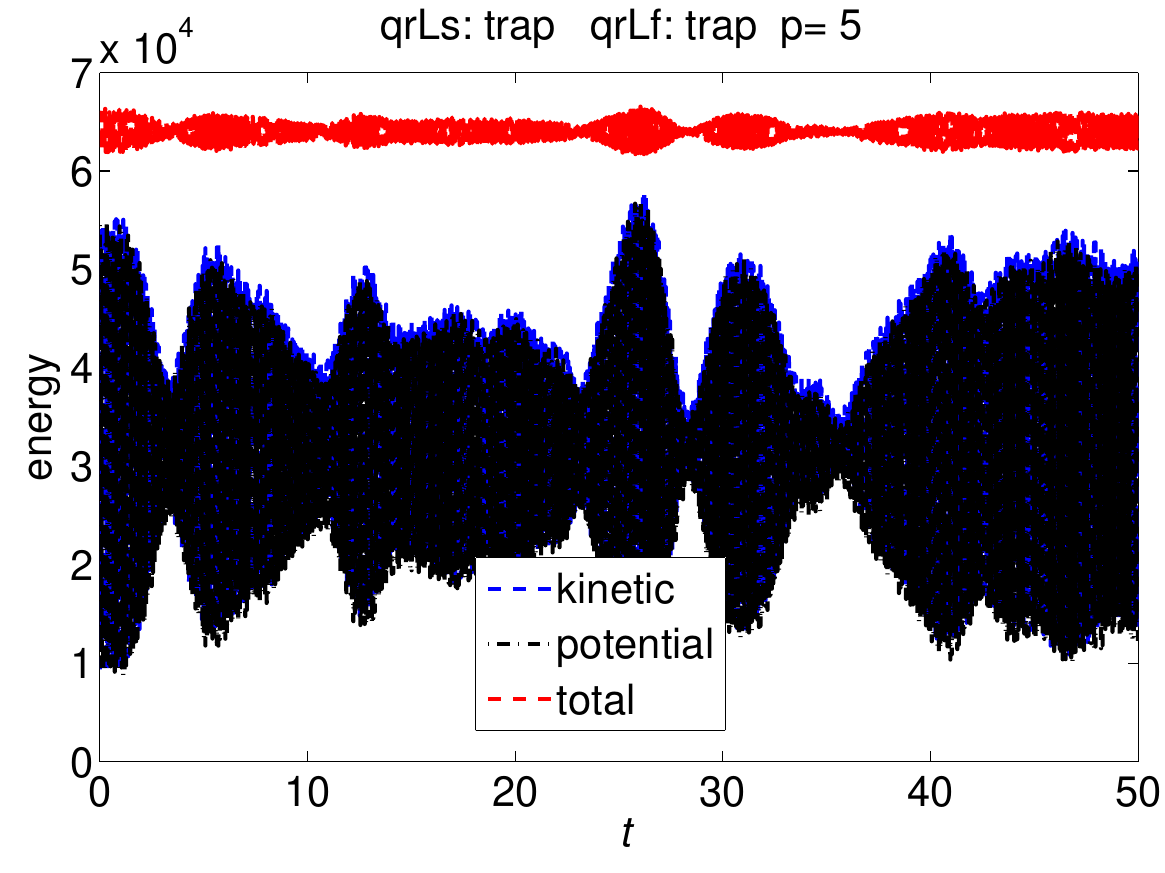}
	\includegraphics[width=0.32\textwidth]{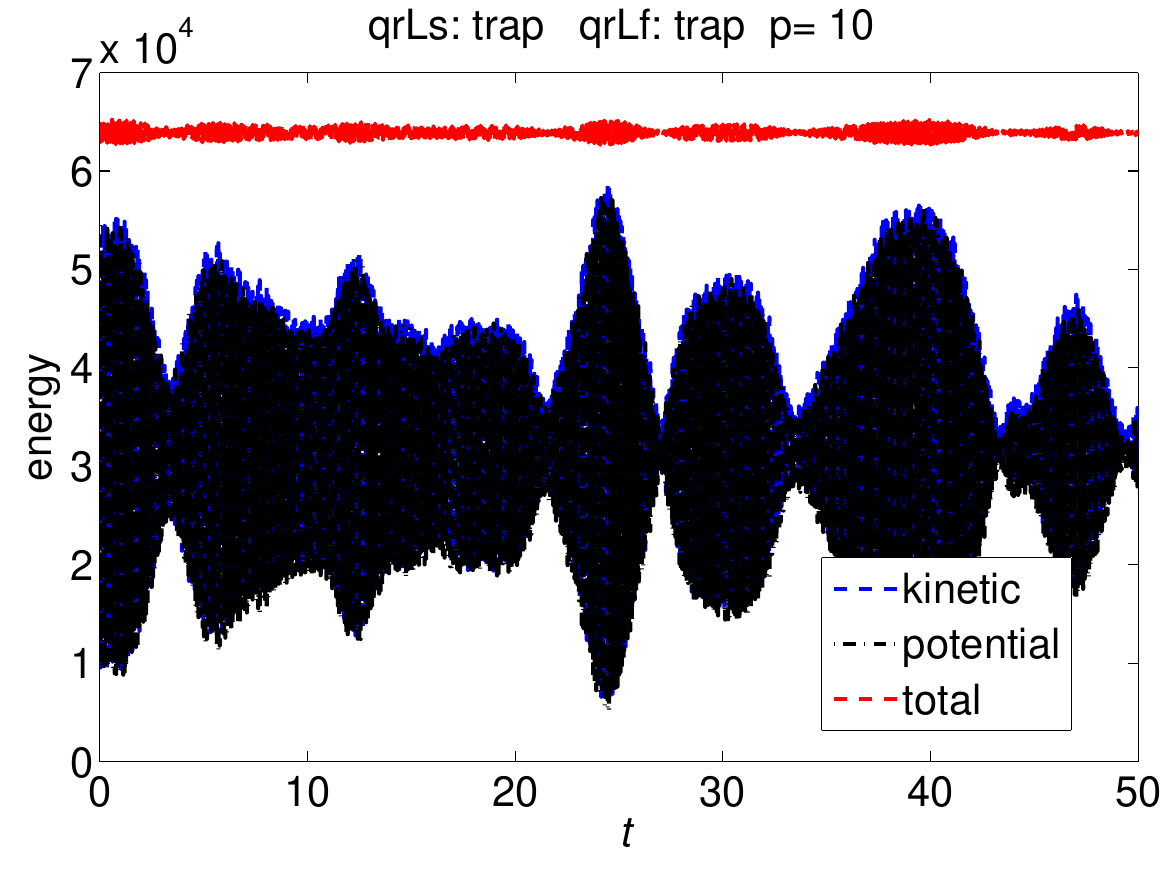}
	\caption{SR: evolution of energy for \(p=1\) (left), \(p=5\) (centre) and \(p=10\) (right) with \(\Delta T =0.01 \) for trapezoidal-trapezoidal scheme}
	\label{fig:sr energy exex}
\end{figure}

The simulation results of the spring ring are presented. First, we investigate the energy evolution of the spring ring. We simulate to \(t_{N}=50\) for all three quadrature schemes.
Fig.~\ref{fig:sr energy imim} shows the kinetic, potential, and total energy versus time \(t\) for the midpoint-midpoint scheme, Fig.~\ref{fig:sr energy exim} for the trapezoidal-midpoint scheme, and Fig.~\ref{fig:sr energy exex} for the trapezoidal-trapezoidal scheme.
In Fig.~\ref{fig:sr energy imim} the curve for the total energy looks like a straight line. However, in the curve of the total energy oscillations are present, which are not visible because their amplitude is several magnitudes smaller than the value of the total energy. The oscillations become visible if the relative change in the energy is depicted.
In the plots for the trapezoidal-midpoint scheme and the trapezoidal-trapezoidal scheme (Fig.~\ref{fig:sr energy exim} and Fig.~\ref{fig:sr energy exex}), 
these oscillations in the total energy are visible. The amplitudes of the oscillations decrease with an increase in \(p\), which means a smaller time step size $\Delta t$.
Furthermore, the error in energy is oscillating around a constant value for the whole simulation time.
No energy is gained nor lost, which shows the good energy behaviour of the presented integrators.
\begin{figure}[h]
	\centering
	\includegraphics[width=0.4\textwidth]{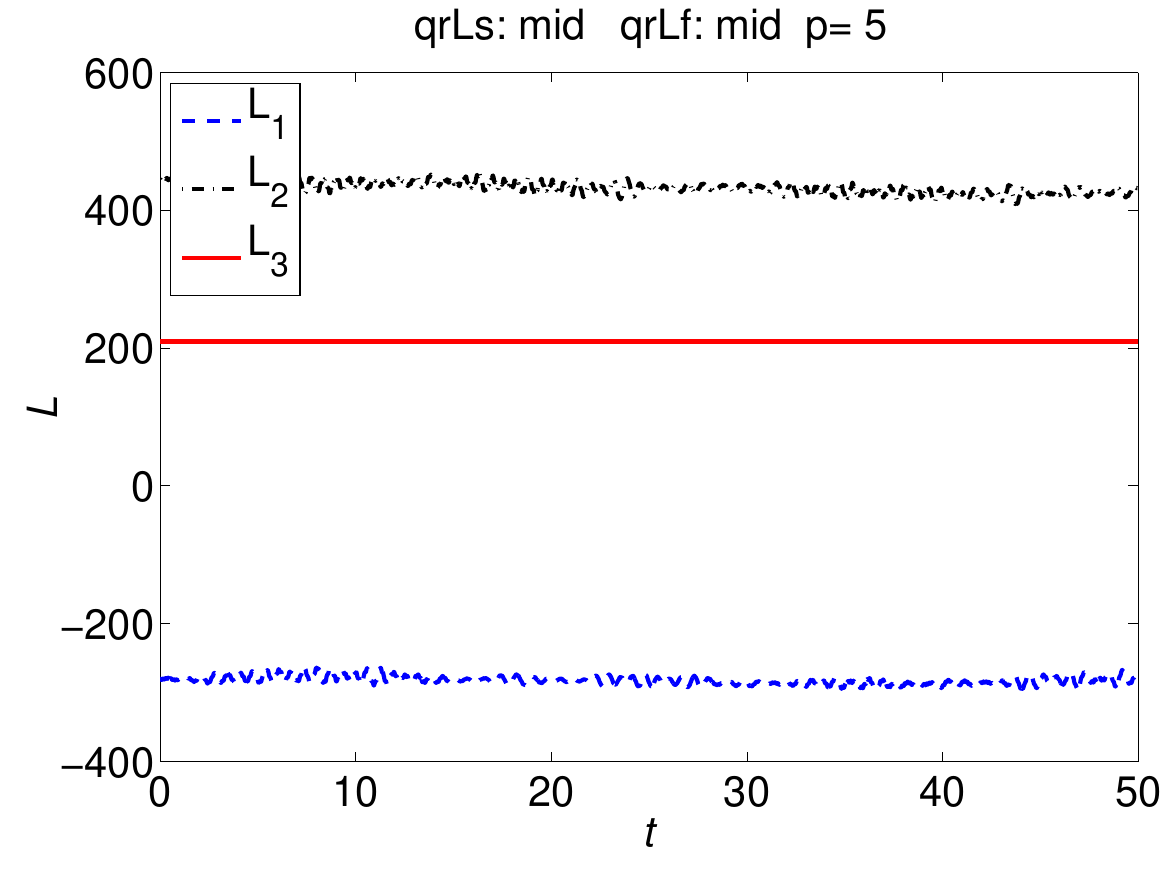}
	\includegraphics[width=0.4\textwidth]{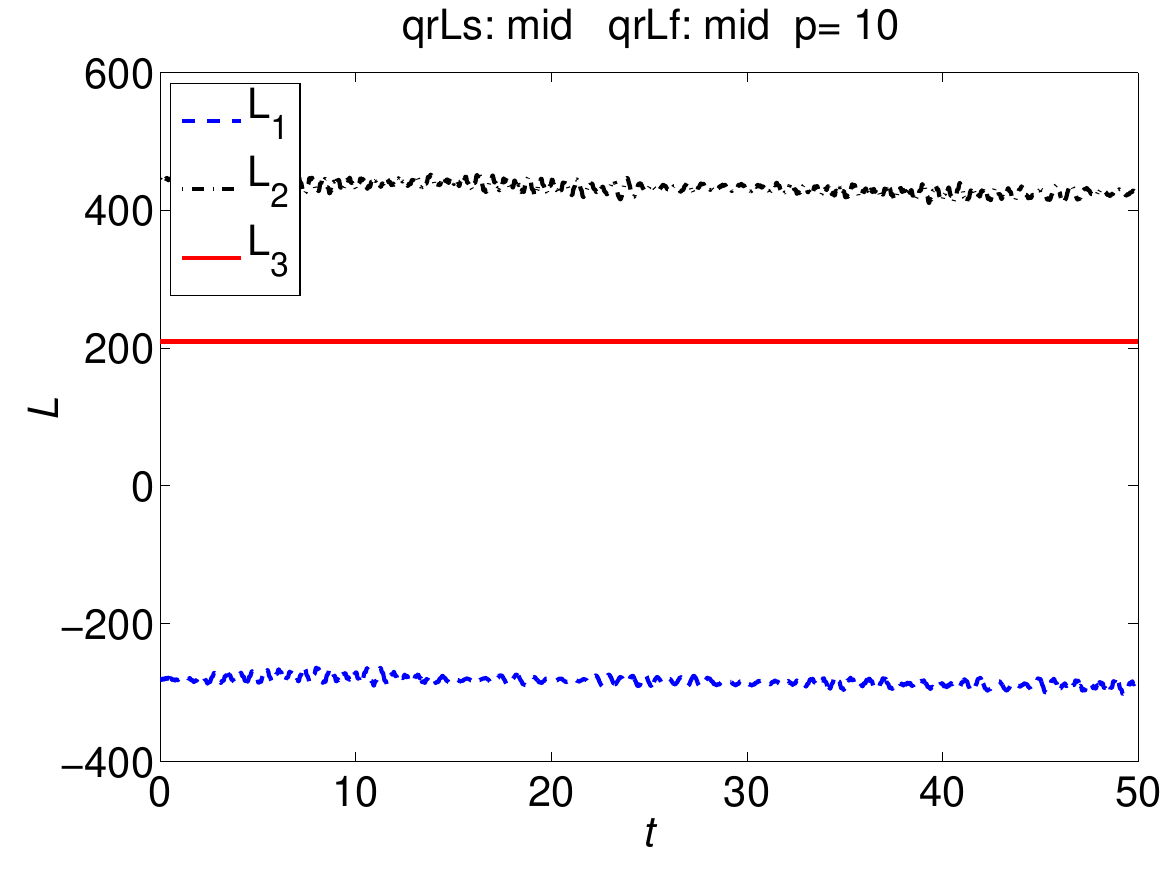}
	\caption{SR: angular momentum in all spatial directions for \(p=5\) (left) and \(p=10\) (right) with \(\Delta T = 0.01  \) for  midpoint-midpoint scheme}
	\label{fig:sr ang mom}
\end{figure}
Secondly, conservation of momentum maps is examined.
Because no external forces and torques act on the spring ring system, and only gravity is present, the angular momentum in gravity direction is invariant.
Fig.~\ref{fig:sr ang mom} shows the angular momentum evolution for a simulation with \(t_{N}=50\) for \(p=5\) and \(p=10\) using the midpoint-midpoint scheme. 
The angular momentum in $e_3$-direction is preserved up to the tolerance set in the Newton-Raphson iteration. As stated in Section \ref{sec:conservation_properties} the variational multirate integrators preserve momentum maps exactly up to numerical accuracy.

\subsection{Numerical convergence study} 
\label{sec:Numerical convergence study}

In this section, the convergence order of the multirate variational integrators is investigated numerically in order to illustrate Theorem \ref{th:approx}. 
	
We calculate the error in $q$ and $p$ with respect to a reference solution $q_{\text{ref}}$ and $p_{\text{ref}}$, respectively. As reference solution serves a single rate VI simulation with a tiny time step size.
	This reference solution is close to the real solution, because the convergence for the single rate VI is proven in \citep{MaWe01}.
	Therefore the convergence to the real solution is shown via the numerical error to the reference solution.
	Since there are two time grids, the error can be calculated on the macro and on the micro time grid.
	The error of $q$ on the macro nodes is defined as
	\begin{equation}
	\label{eq:err_mac}
	e_q^{mac} = \sup_{k=0,\ldots,N} \{\| q_k - q_{\text{ref}}(t_k)\|\} 
	\end{equation}
	and on the micro nodes as
	\begin{equation}
	\label{eq:err_mic}
	e_q^{mic}  = \sup_{\begin{array}{c}k=0,\ldots,N-1\\m=1,\ldots,p-1 \end{array}} \{\| q_{k}^{f,m} - q^f_{\text{ref}}(t_{k}^{m})\| \} 	
	\end{equation}
	Note that the macro nodes (\(m=0\), \(m=p\)) are excluded in \eqref{eq:err_mic}. Analogously, the error of $p$ is defined.
		
We investigate the convergence order of the variational multirate integrators with \(p=5\) and \(p=10\) numerically by means of the FPU and the spring ring example. For the reference solution of the FPU a single rate simulation with time step size \(\Delta T=5 \cdot 10^{-7} \) and simulation time \(t_{N}=0.5 \) is used. The reference solution of the spring ring is computed via a single rate integrator with time step size \(\Delta T = 1 \cdot 10^{-6} \) and \(t_{N} = 0.5  \). For each quadrature combination scheme, a corresponding reference solution is computed.

\begin{figure}[h]
	
		\centering
		\includegraphics[width=0.49\textwidth]{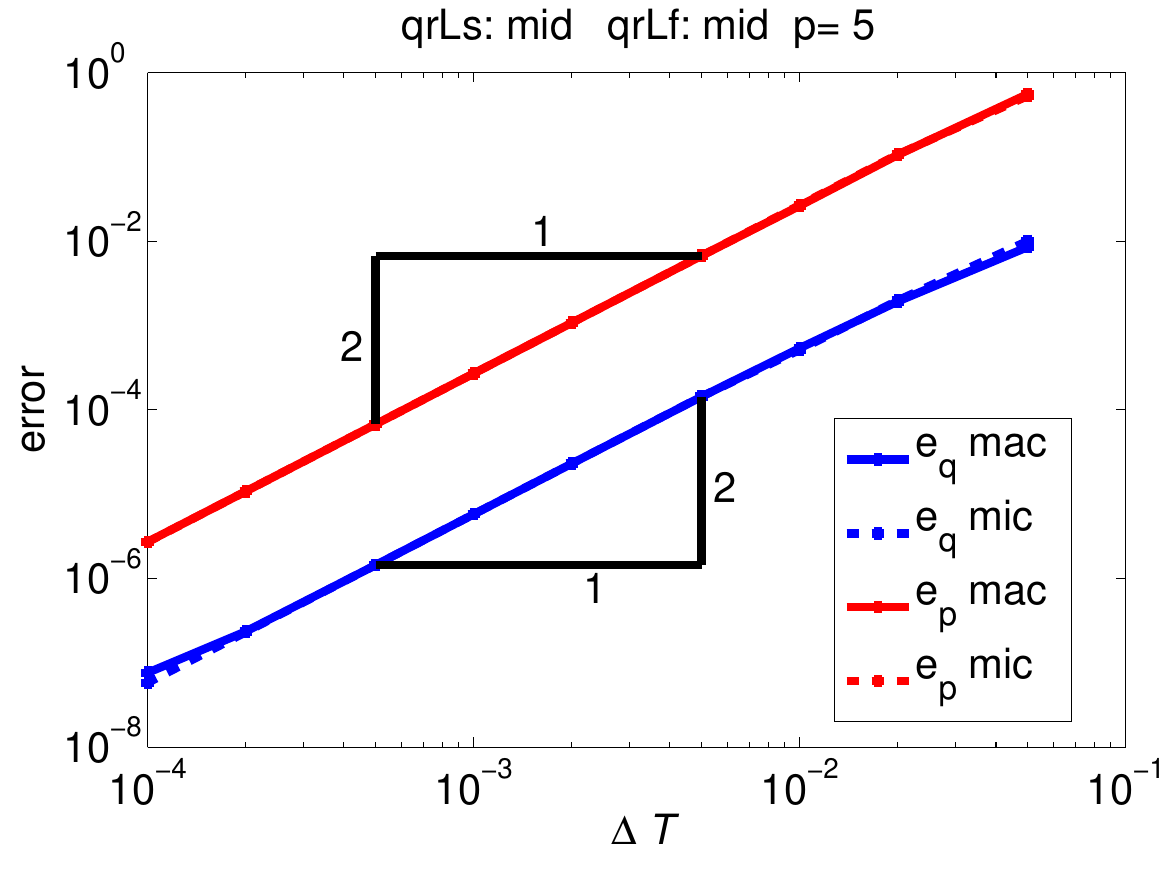}
		\includegraphics[width=0.49\textwidth]{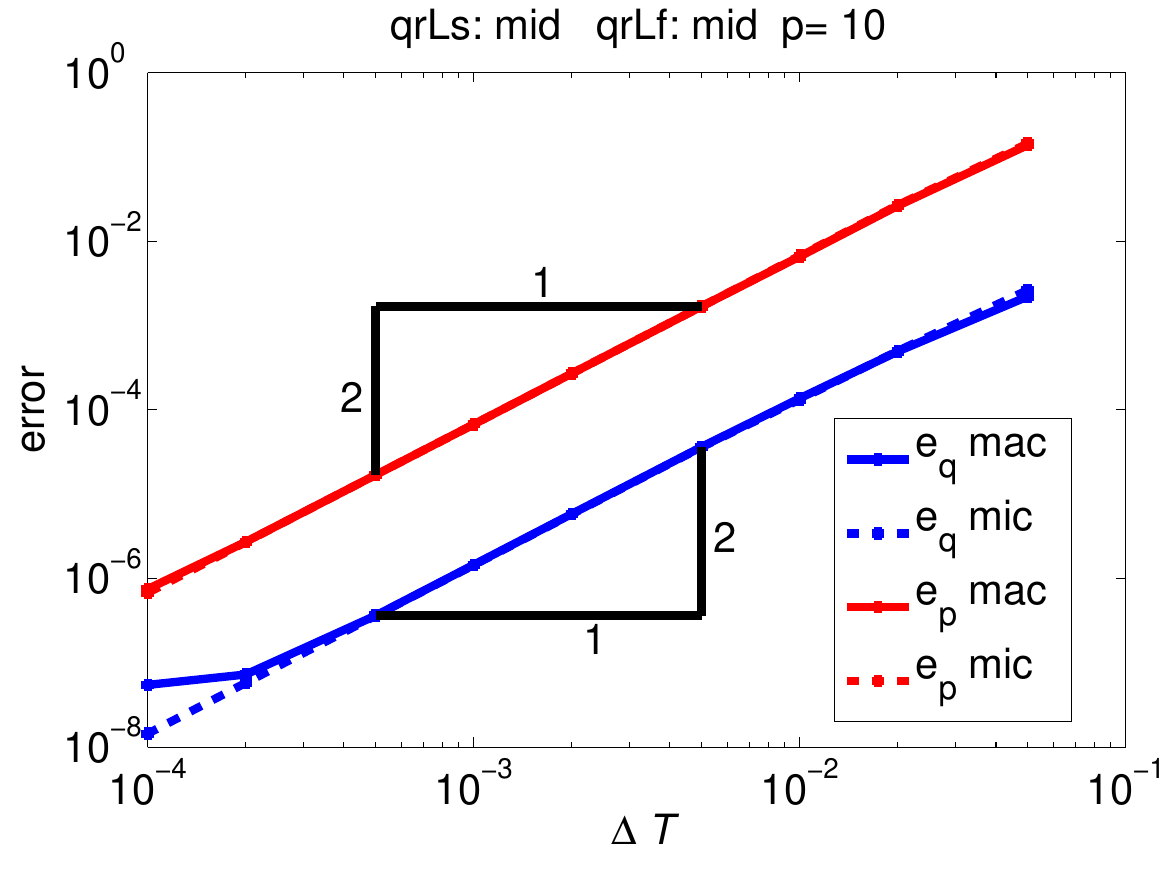}
		\caption{FPU: error in configuration \(q\) (blue) and conjugate momentum \(p\) (red) for macro nodes (solid line) and micro nodes (dashed line) for \(p=5\) (left) and \(p=10\) (right), midpoint-midpoint quadrature}
		\label{fig:fpu_imim_p_macmic}
	\end{figure}

\begin{figure}[h]
	
		\centering
		\includegraphics[width=0.49\textwidth]{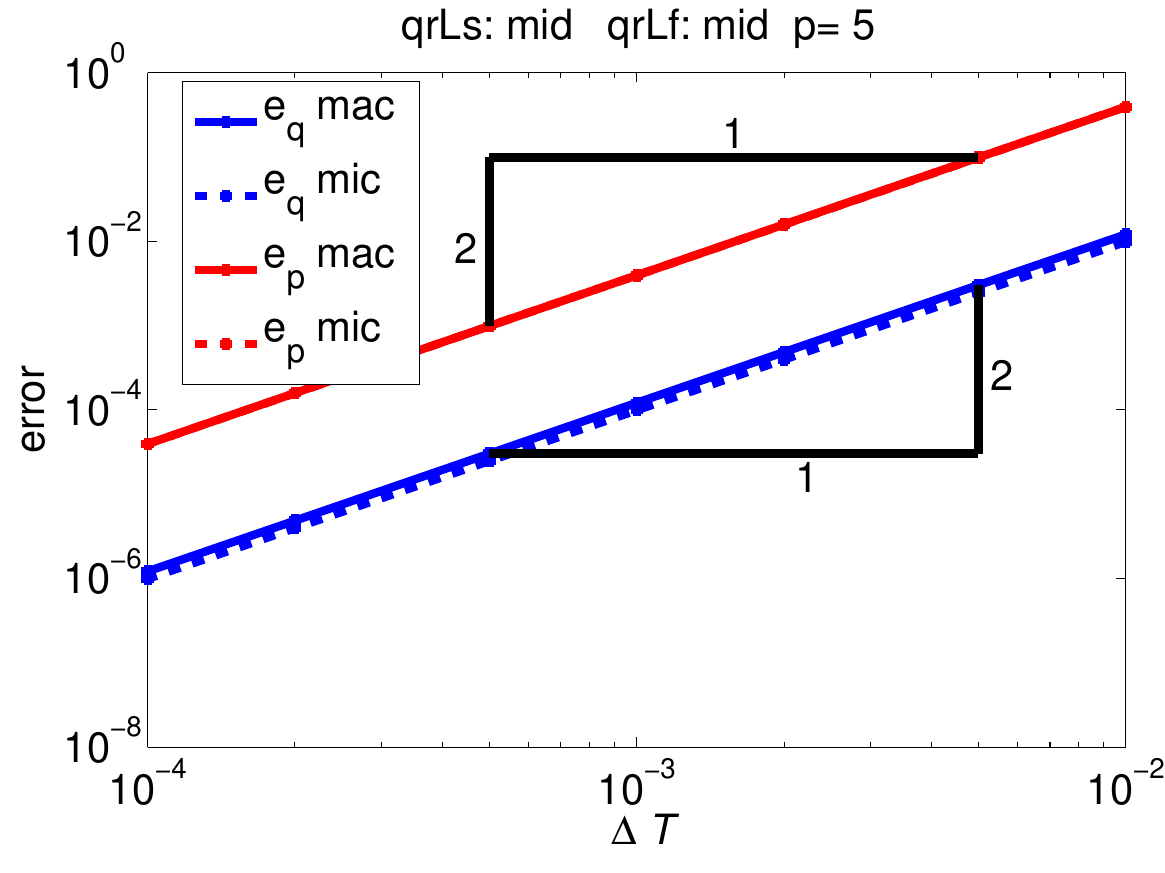}
		\includegraphics[width=0.49\textwidth]{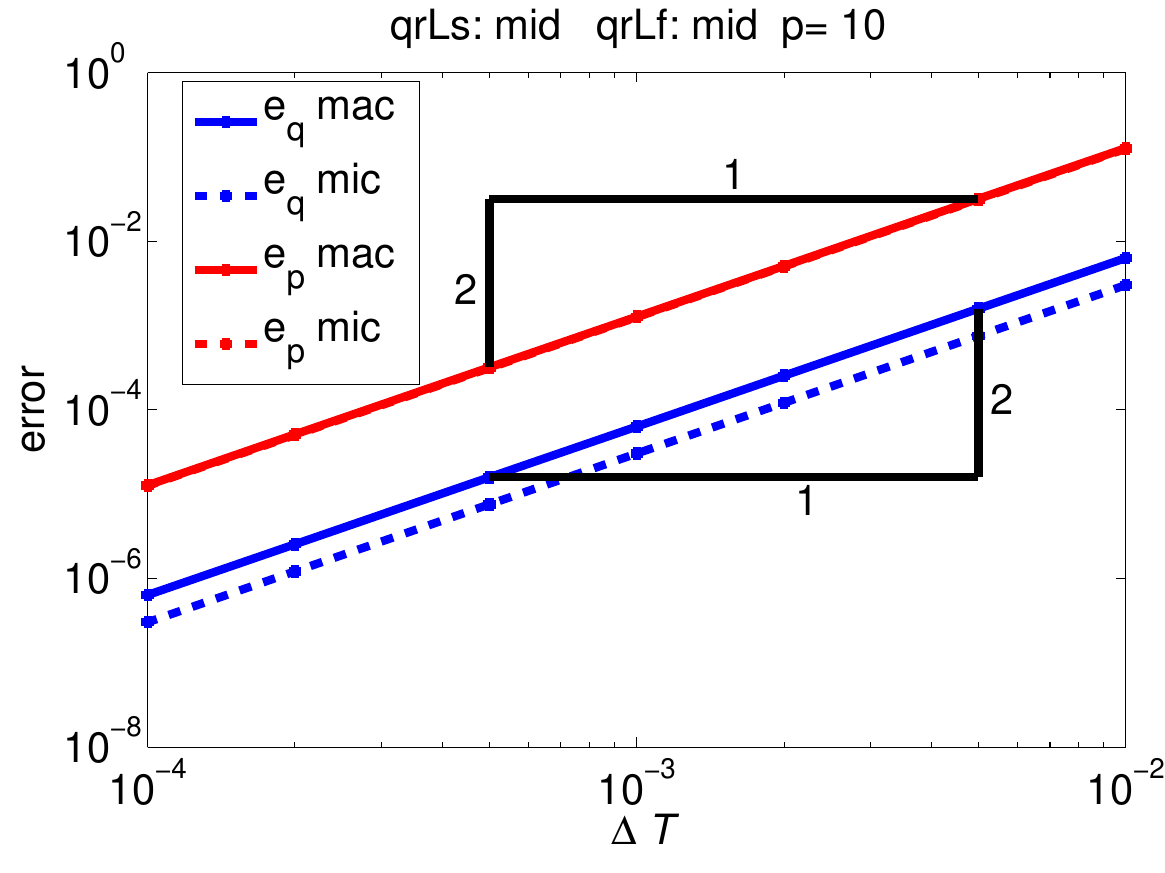}
		\caption{SR: error in configuration \(q\) (blue) and conjugate momentum \(p\) (red) for macro nodes (solid line) and micro nodes (dashed line) for \(p=5\) (left) and \(p=10\) (right), midpoint-midpoint quadrature}
		\label{fig:sr konv_imim_qp_macmic}
	\end{figure}

First, the convergence order of the midpoint-midpoint scheme is examined numerically. Fig.~\ref{fig:fpu_imim_p_macmic} and Fig.~\ref{fig:sr konv_imim_qp_macmic} show the error of configuration \(q\) and conjugate momentum \(p\) on macro nodes (solid) and of $q^f$ and $p^f$ on micro nodes (dashed) versus \(\Delta T \) with \(p=5\) and \(p=10\) for FPU and SR. On each plot the convergence of the error in configuration \(q\) resp.~$q^f$ and in the conjugate momentum \(p\) resp.~$p^f$ is of order 2. In Fig.~\ref{fig:fpu_imim_p_macmic} only solid lines are visible, because here the error on the macro nodes equals the error on the micro nodes. 

\begin{figure}[h]
		\centering
		\includegraphics[width=0.49\textwidth]{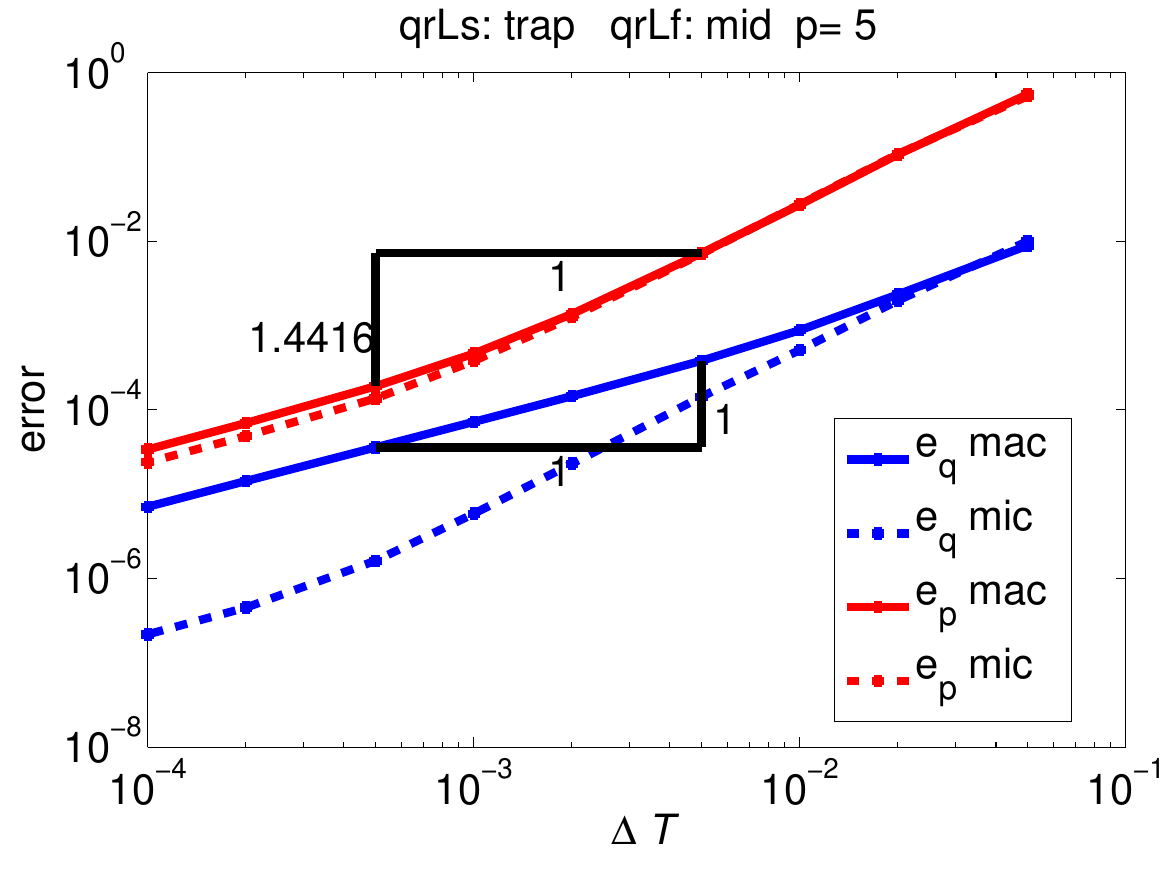}
		\hfill
		\includegraphics[width=0.49\textwidth]{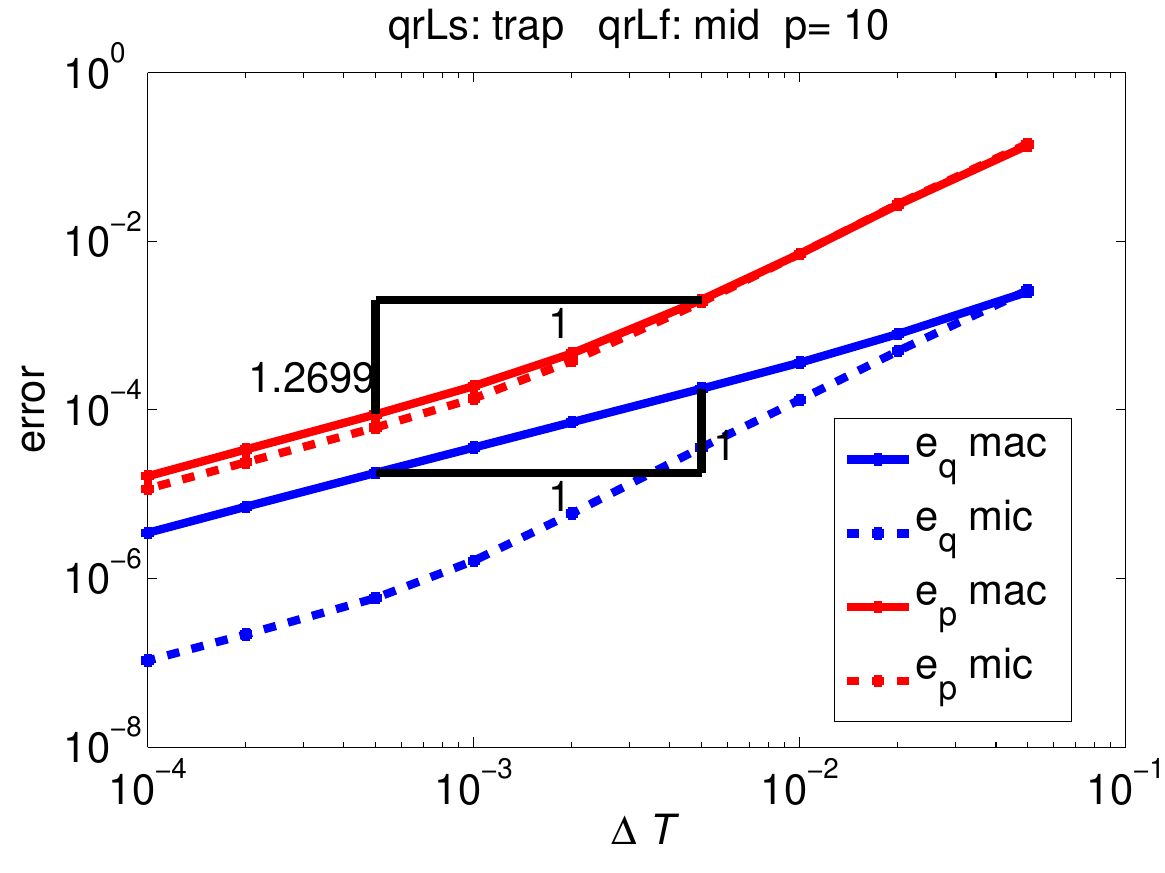}
		\caption{FPU: error in configuration \(q\) (blue) and conjugate momentum \(p\) (red) for macro nodes (solid line) and micro nodes (dashed line) for \(p=5\) (left) and \(p=10\) (right), trapezoidal-midpoint quadrature}
		\label{fig:fpu konv_exim}
	\end{figure}

\begin{figure}[h]
		\includegraphics[width=0.49\textwidth]{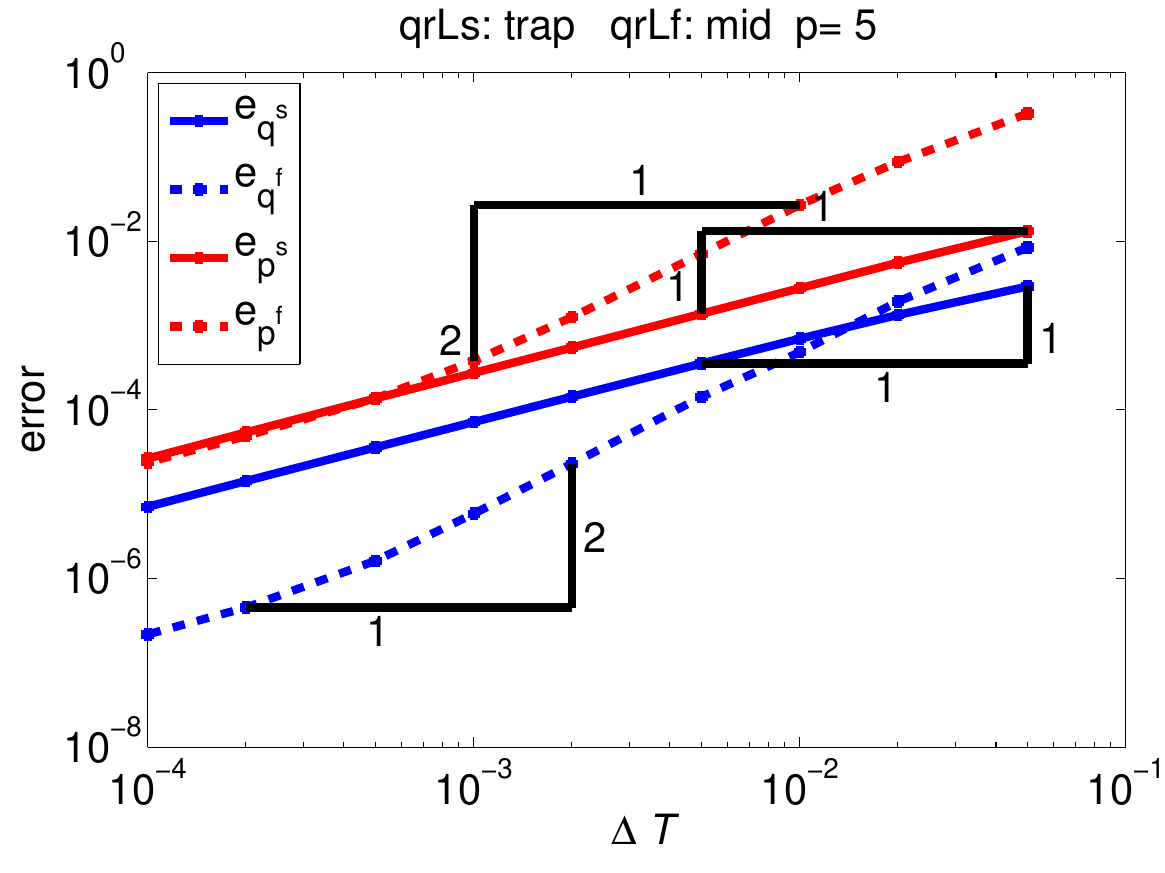}
		\includegraphics[width=0.49\textwidth]{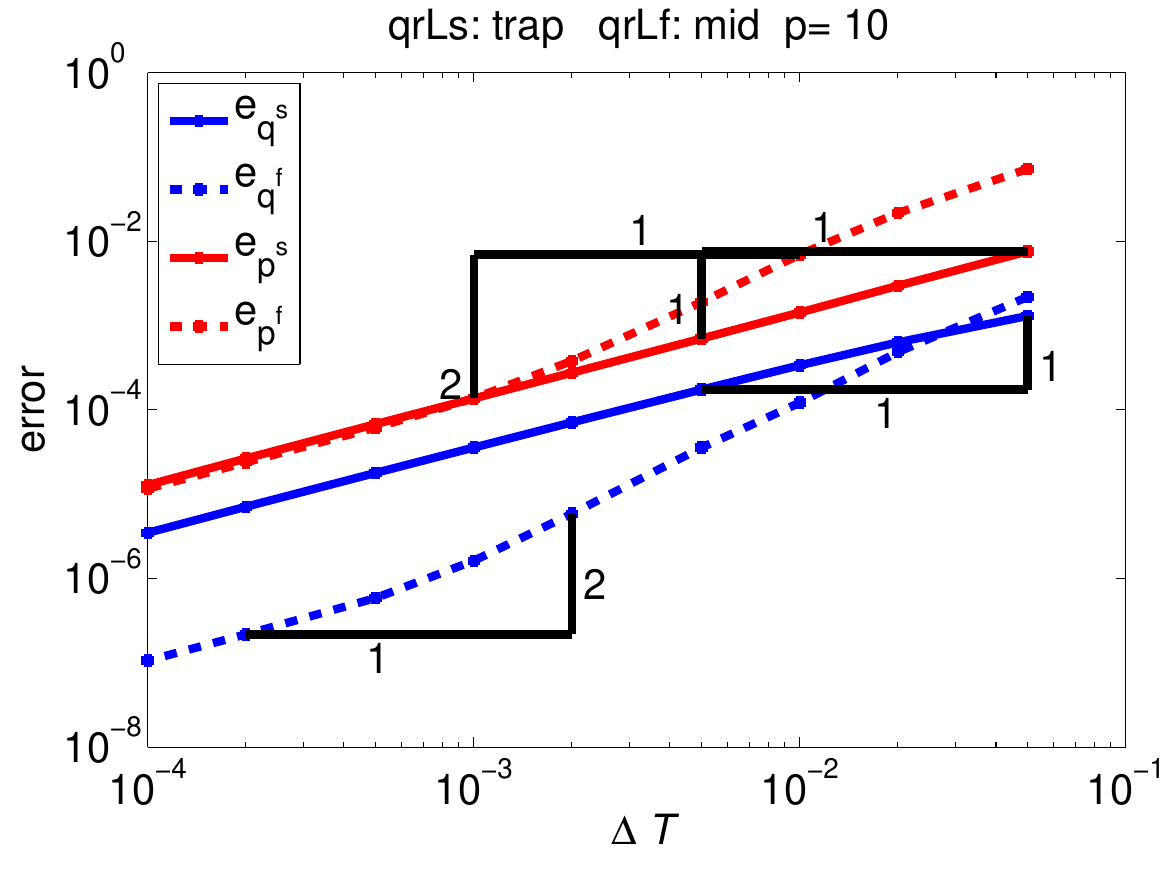}
		\caption{FPU: error on macro nodes split into slow (solid) and fast (dashed)  \(q\) (blue) and \(p\) (red) for \(p=5\) (left) and \(p=10\) (right), trapezoidal-midpoint quadrature}
		\label{fig:fpu err mac qsqf}
	\end{figure}

\begin{figure}[h]
		\centering
		\includegraphics[width=0.49\textwidth]{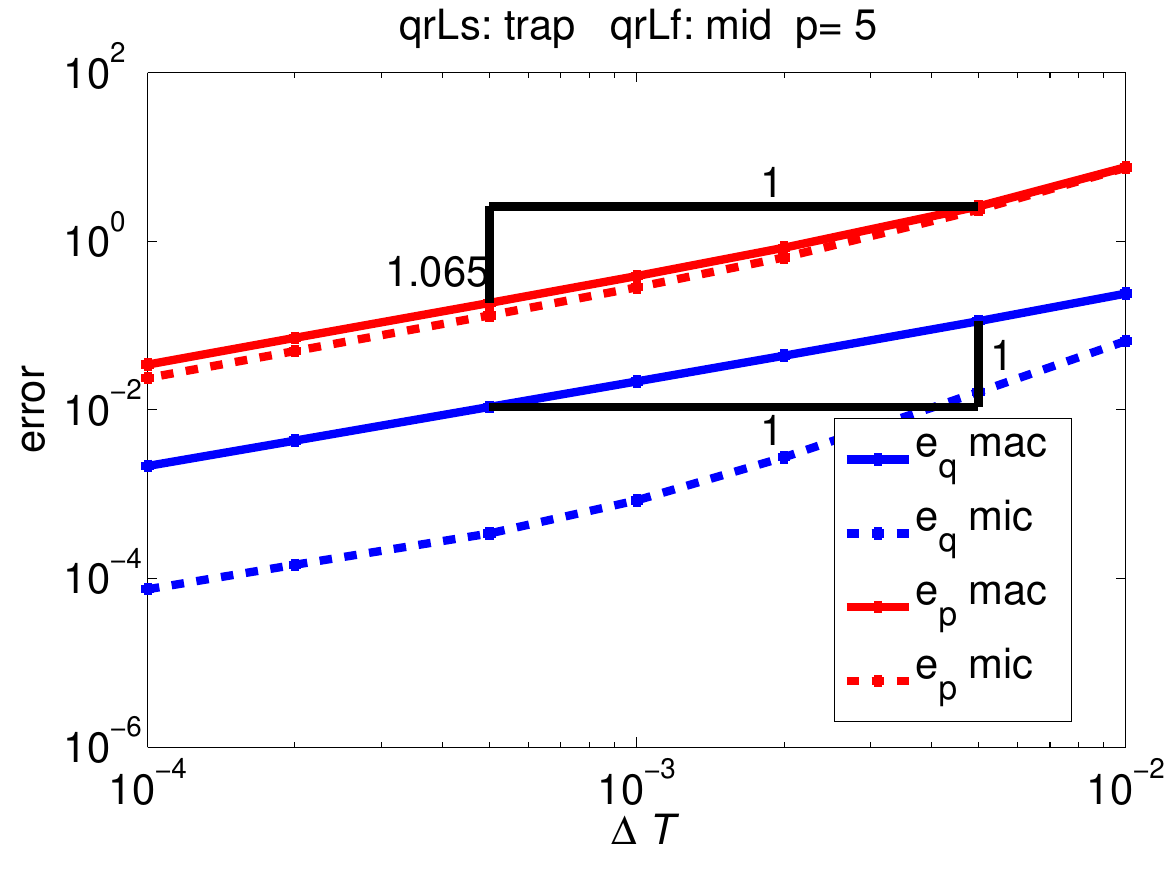}
		\hfill
		\includegraphics[width=0.49\textwidth]{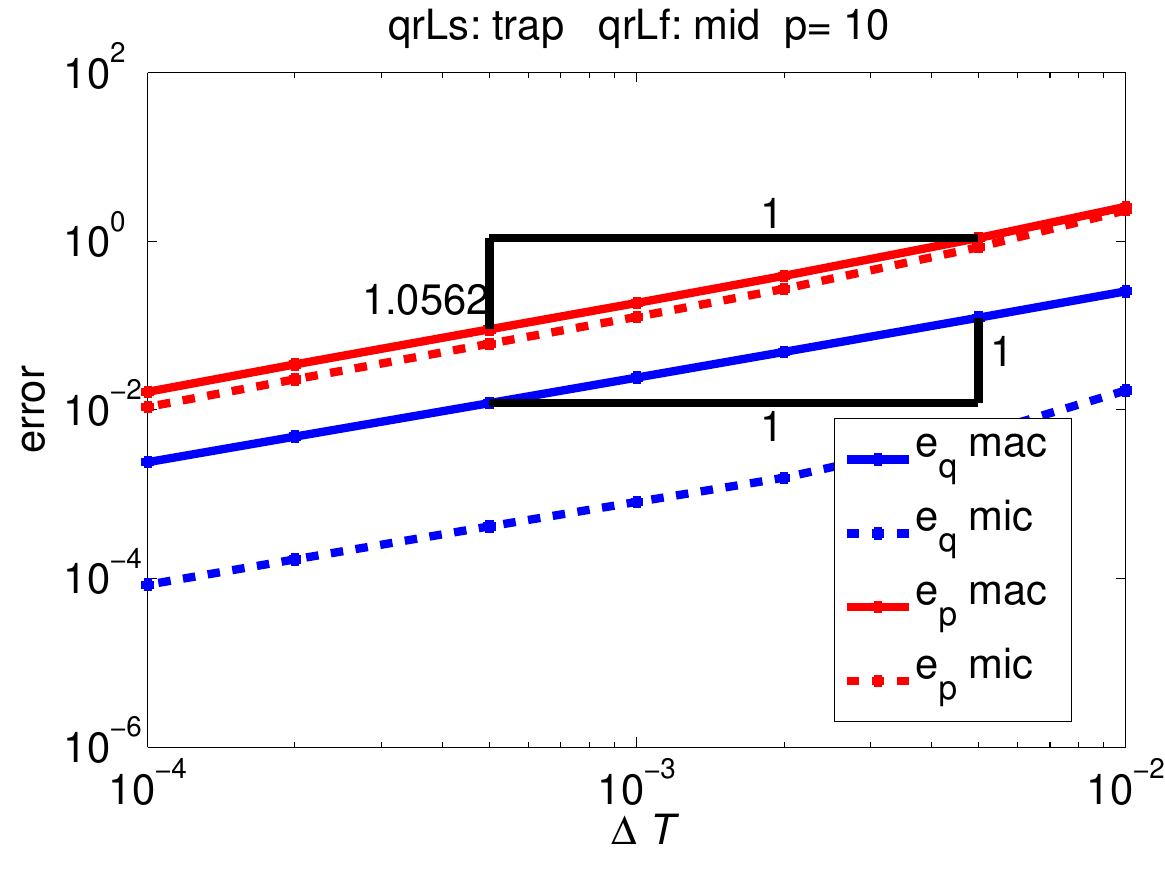}
		\caption{SR: error in configuration \(q\) (blue) and conjugate momentum \(p\) (red) for macro nodes (solid line) and micro nodes (dashed line) for \(p=5\) (left) and \(p=10\) (right), trapezoidal-midpoint quadrature}
		\label{fig:sr_konv_exim_qp_macmic}
	\end{figure}

\begin{figure}[h]	
		\includegraphics[width=0.49\textwidth]{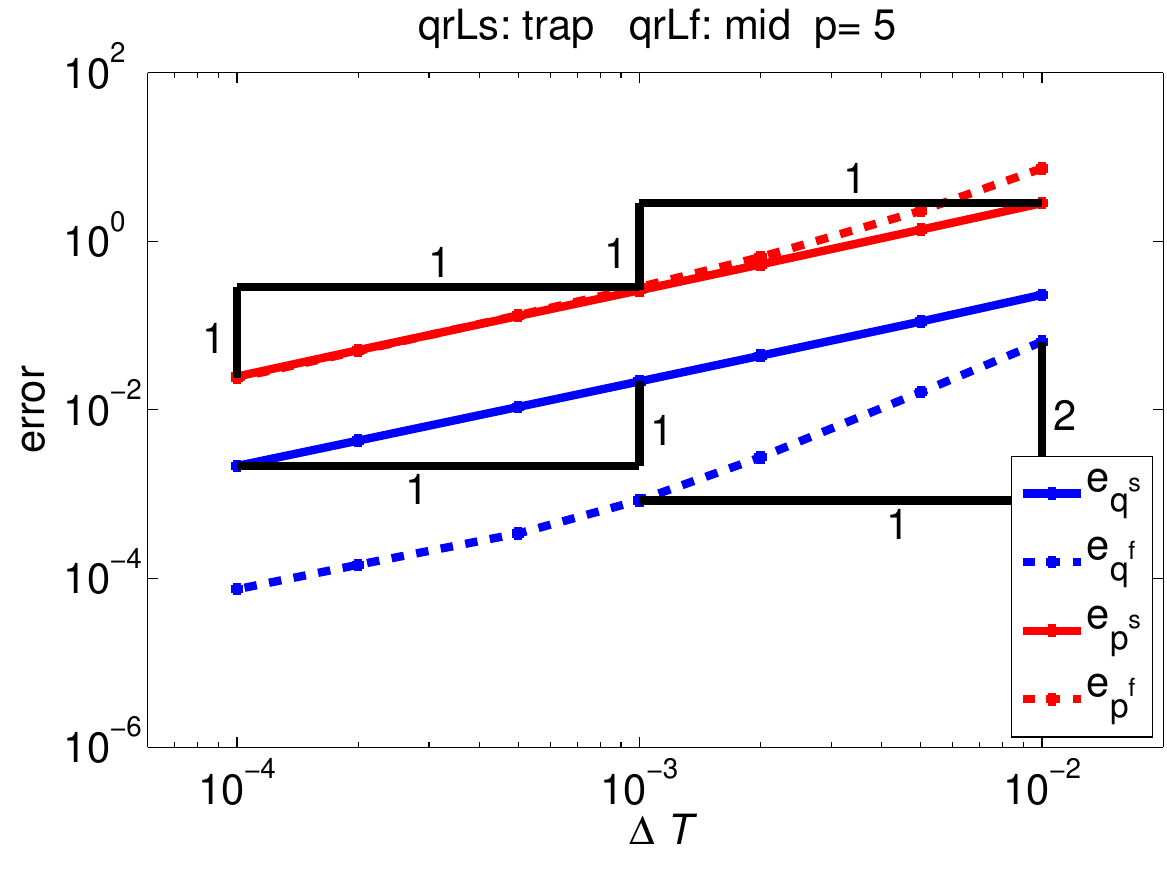}
		\includegraphics[width=0.49\textwidth]{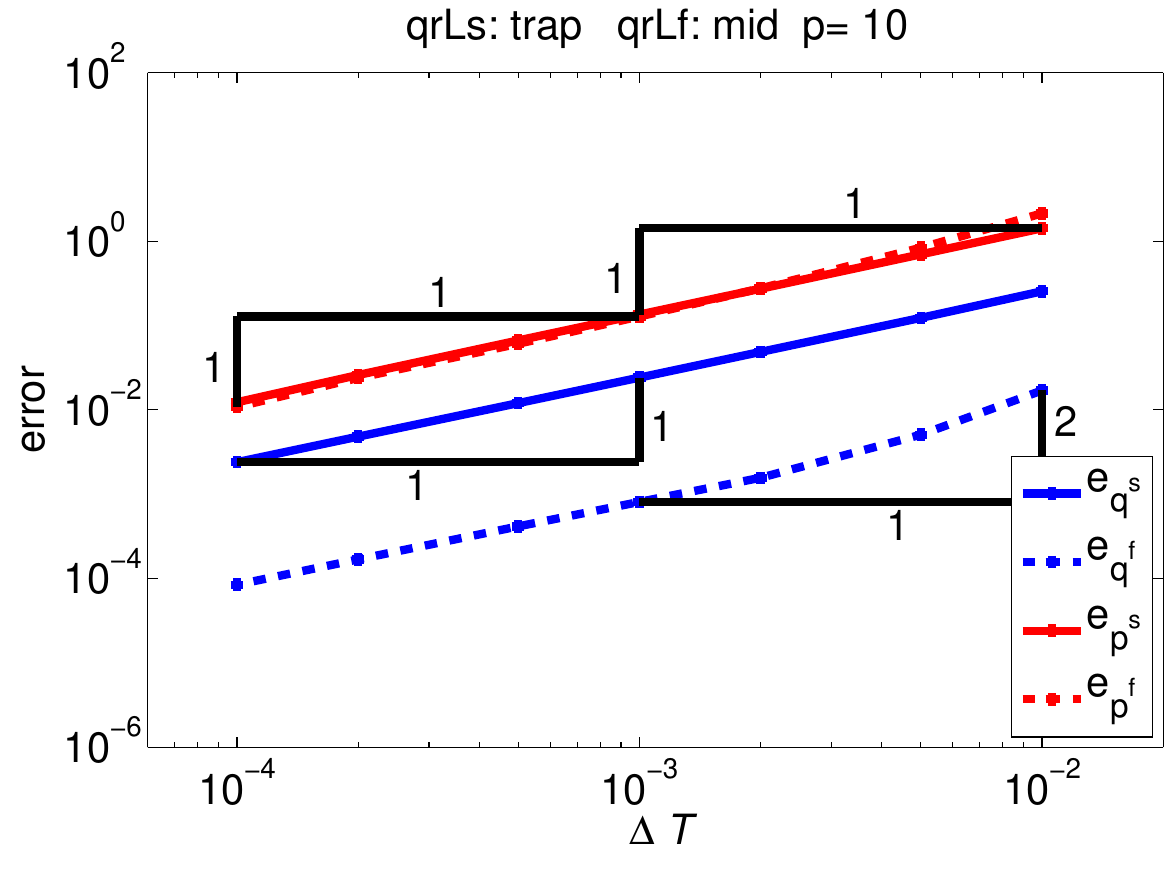}
		\caption{SR: error on macro nodes split into slow (solid) and fast (dashed) variables for \(p=5\) (left) and \(p=10\) (right) in trapezoidal-midpoint quadrature}
		\label{fig:sr err mac slow fast}
	\end{figure}

Next, we show the numerical convergence order of the trapezoidal-midpoint scheme. For the approximation of the integral of the slow potential the left rectangular rule is used ($\{ \alpha_V, \gamma_V \} \in \{0,1\}, \gamma_V=\alpha_V)$. 
Fig.~\ref{fig:fpu konv_exim} and \ref{fig:sr_konv_exim_qp_macmic} illustrate the error of the configuration \( q \) and the conjugate momentum \( p \) on macro nodes (solid) and of $q^f$ and $p^f$ on micro nodes (dashed) versus \(\Delta T \) for \(p=5\) (left plot) and \(p=10\) (right plot) for FPU (Fig.~\ref{fig:fpu konv_exim}) and SR (Fig.~\ref{fig:sr_konv_exim_qp_macmic}).
In Fig.~\ref{fig:fpu konv_exim}, we see that the error of \(q\) on the macro nodes, $e^{\text{mac}}_q$, converges with an order of one.
The error of $q^f$ on the micro nodes, $e^{\text{mic}}_q$, starts with a convergence order of two, but reduces for tiny time steps to one.
Similarly, the error of \(p\) on the macro nodes, $e^{\text{mac}}_p$, and of $p^f$ on the micro nodes, $e^{\text{mic}}_p$, starts with an order of two, but reduces to one for smaller time steps. 
An explanation is given in Fig.~\ref{fig:fpu err mac qsqf}, where the error on the macro nodes split into slow and fast variables for \(q\) and \(p\) versus \(\Delta T\) is presented. We see that the error in the slow variables $q^s$ and $p^s$ converges with order one. 
The error in the fast variables $q^f$ and $p^f$ firstly reduces quadratically until the errors in $q^s$ 
start to influence the approximations of $q^f$ and $p^f$ through the coupling potential $V$ and eventually dominate the error {and the error reduces linearly.}, Thus, the convergence order of $q^f$ and $p^f$ is also one. Fig.~\ref{fig:fpu konv_exim} can be interpreted in the same way. The error
$e^p_{\text{mac}}$ is first dominated by the error in $p^f$ on the macro nodes explaining the steeper slope for the larger time steps, while the error in $p^s$ dominates for small time steps resulting in convergence rates greater than one but lower than two. The order reduction of $e^{\text{mic}}_q$ can be explained again by the potential $V$ coupling the approximation of $q^f$ and $q^s$. The discrete fast conjugate momenta depend on the slow configuration what is the reason for the order reduction in $e^{\text{mic}}_p$.
In Fig.~\ref{fig:sr_konv_exim_qp_macmic} and \ref{fig:sr err mac slow fast}, the errors in the simulation of the spring ring are examined.  Both plots in Fig.~\ref{fig:sr_konv_exim_qp_macmic} show that the error of \(q\) on the macro nodes converges again with an order of one.
As seen before for the FPU, changes in the slope from two to one are visible for the errors $e^{\text{mac}}_p$, $e^{\text{mic}}_p$ and $e^{\text{mic}}_q$, when the time step $\Delta T$ is decreased. 
Fig.~\ref{fig:sr err mac slow fast} shows the error for \(q\) and \(p\) on macro nodes split into slow and fast variables. Again, 
the error in the slow variables $q^s$ and $p^s$ converges with order one. The convergence order of the error in the fast variables $q^f$ resp.~$p^f$ starts with two and changes to one. Reason for the order reductions of the error in $q^f$ resp.~$p^f$ on the macro nodes as well as on the micro nodes is again the influence of the coupling potential. Summarized, as the errors in \(q^s\) and \(p^s\) on the macro nodes are dominating, see Fig.~\ref{fig:sr err mac slow fast}, the order of convergence of the error in \(q\) and \(p\) on the macro nodes is one, respectively.

\begin{figure}[h]
\centering
		\includegraphics[width=0.49\textwidth]{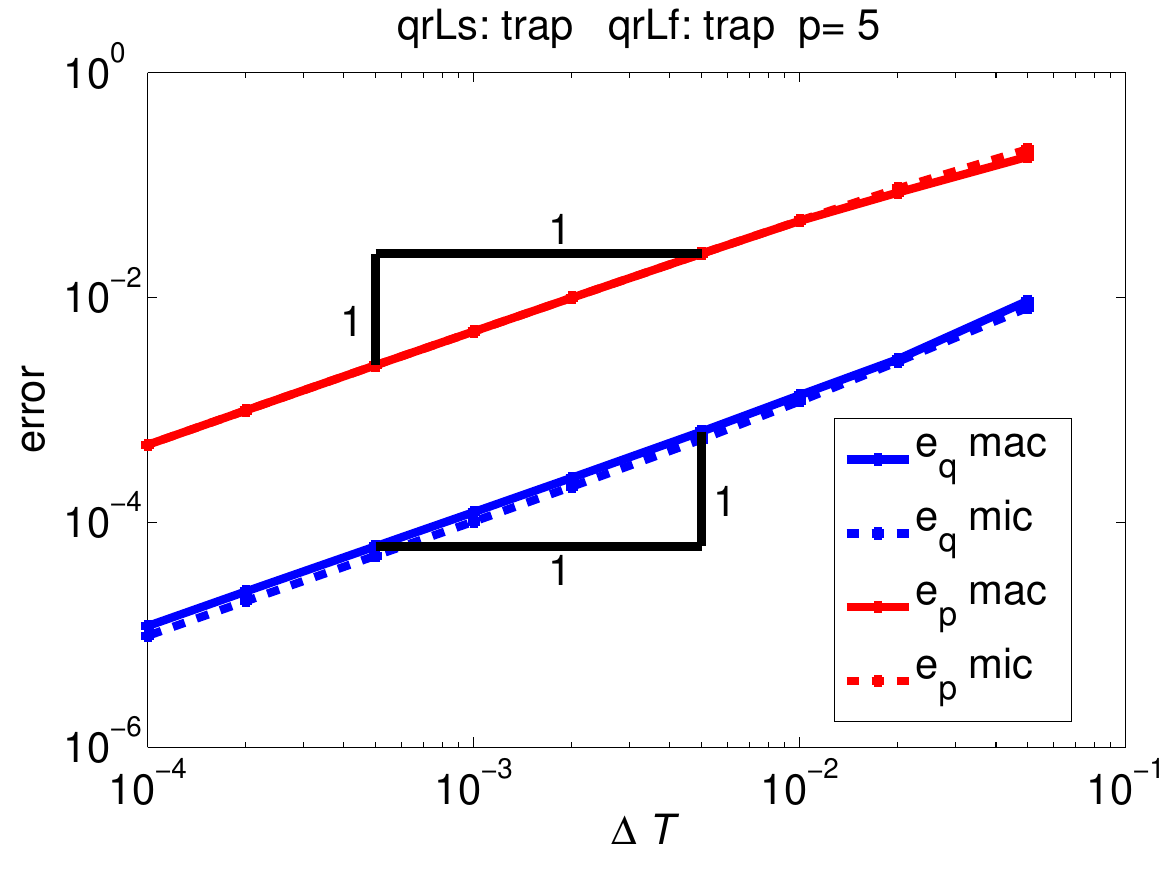}
		\hfill
		\includegraphics[width=0.49\textwidth]{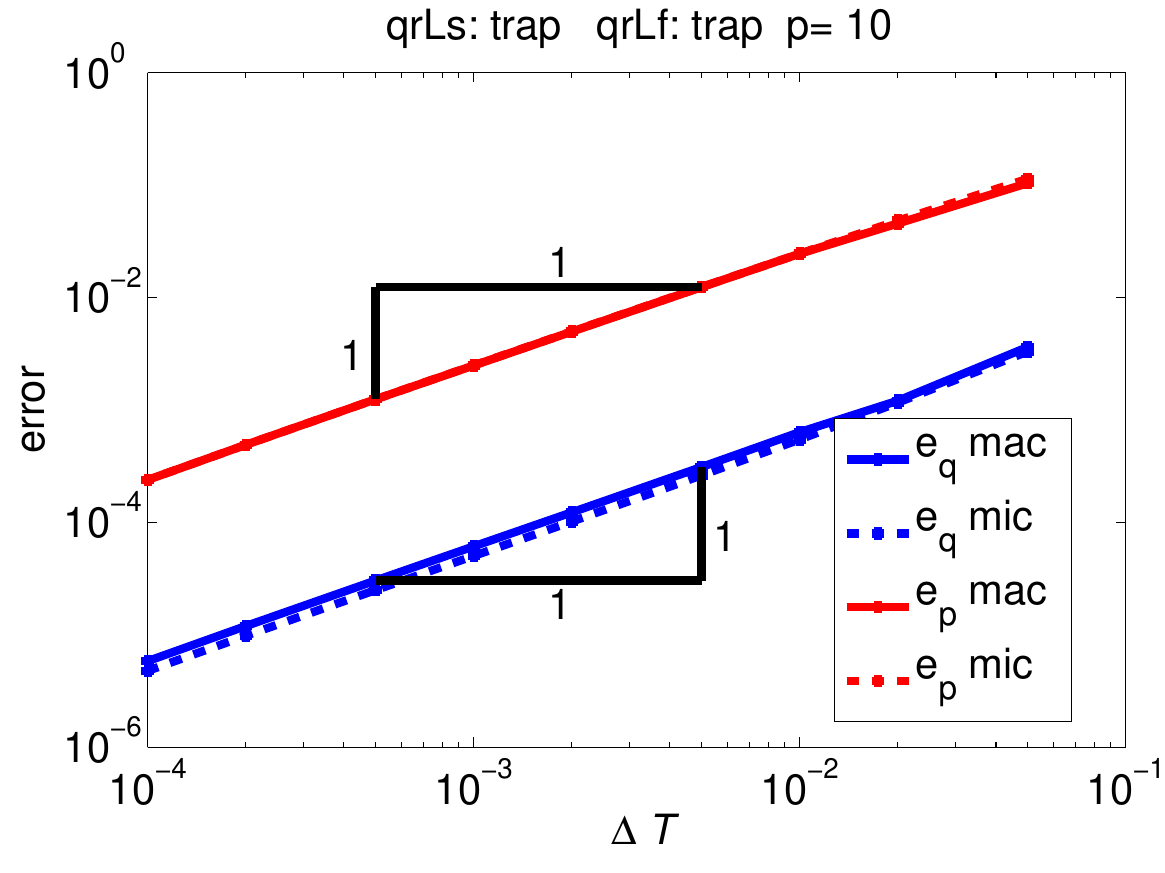}
		\caption{FPU: error in configuration q (blue) and conjugate momentum p (red) for macro nodes (solid line) and micro nodes (dashed line) for \(p=5\) (left) and \(p=10\) (right) micro steps, trapezoidal-trapezoidal quadrature}
		\label{fig:fpu konv_exex_p_macmic}
	\end{figure}

\begin{figure}[h]
	
		\centering
		\includegraphics[width=0.49\textwidth]{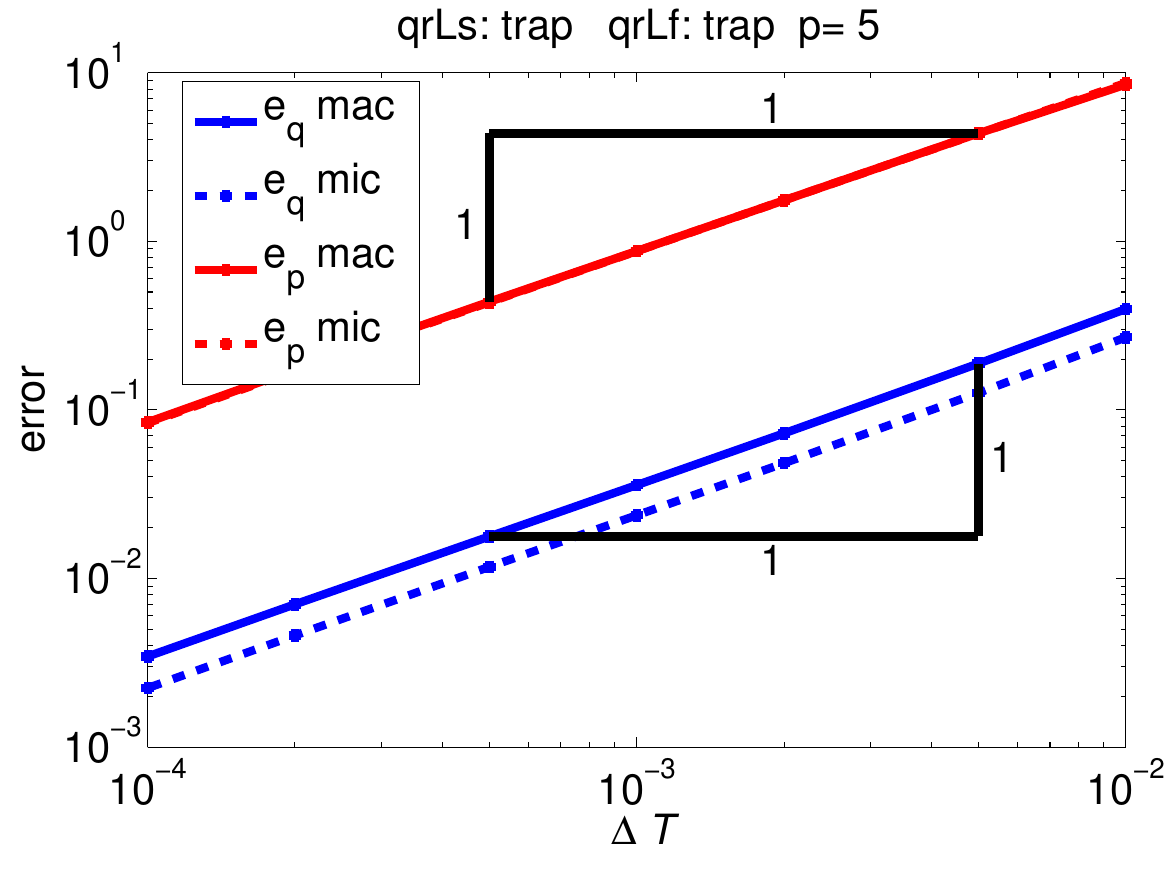}
		\hfill
		\includegraphics[width=0.49\textwidth]{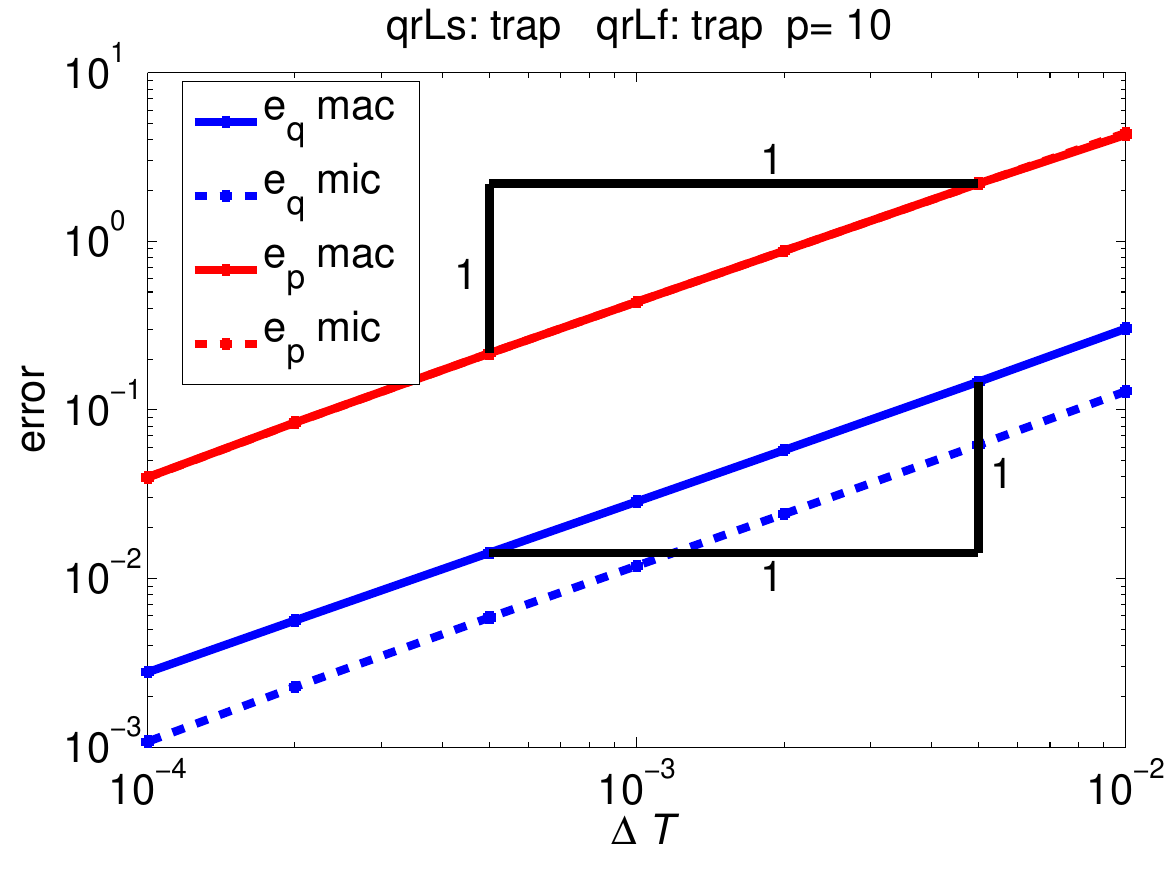}
		\caption{SR: error in configuration \(q\) (blue) and conjugate momentum \(p\) (red) for macro nodes (solid line) and micro nodes (dashed line) for \(p=5\) (left) and \(p=10\) (right), trapezoidal-trapezoidal quadrature}
		\label{fig:sr_konv_exex_p_macmic}
	\end{figure}

Finally, the convergence of the trapezoidal-trapezoidal scheme is discussed. Both, the integral of the slow and of the fast potential are approximated with the left rectangular rule ($\{ \alpha_V , \gamma_V \} \in \{0,1\}, \gamma_V=\alpha_V$ and $\{ \alpha_W , \gamma_W \} \in \{0,1\}, \gamma_W=\alpha_W)$.
Fig.~\ref{fig:fpu konv_exex_p_macmic} and \ref{fig:sr_konv_exex_p_macmic} show the error of the configuration \( q \) and the conjugate momentum \( p \) on macro nodes (solid) and of $q^f$ and $p^f$ on micro nodes (dashed) versus \(\Delta T \) for \(p=5\) and \(p=10\) for FPU and SR.
The errors in the configuration and the conjugate momenta convergence with order one.
	
A summary of the convergence orders of the configuration $q$ and conjugate momentum $p$ on macro nodes is given in Table~\ref{tab:convergence rates} for the different quadrature combinations. Regarding the midpoint-trapezoidal quadrature we have seen a convergence order in $p$ lower than two but higher than one for FPU. However this is a result that holds for FPU specifically but not generally. In general we can expect a convergence order of at least 1. The results given in Table \ref{tab:convergence rates} confirm the statements in Theorem \ref{th:approx}.
\begin{table}
\caption{Convergence orders of configuration $q$ and conjugate momentum $p$ on macro nodes }
\label{tab:convergence rates}     
\begin{tabular}{lccc}
\hline\noalign{\smallskip}
& midpoint-midpoint & trapezoidal-midpoint & trapezoidal-trapezoidal \\
&  &  \footnotesize{ $\{ \alpha_V , \gamma_V \} \in \{0,1\}$} &  \footnotesize{ $\{\alpha_i , \gamma_i \} \in \{0,1\}$} \\
&  &  $\gamma_V=\alpha_V$  &  \footnotesize{ $ \gamma_i=\alpha_i, i=\{V,W\}$} \\
\noalign{\smallskip}\hline\noalign{\smallskip}
		 p & $ q  \qquad \qquad  p$          &  $ q  \qquad \qquad   p$        &   $ q  \qquad \qquad   p$  \\
\noalign{\smallskip}\hline\noalign{\smallskip}
	 5 & $ 2  \qquad \qquad  2$         &   $ 1  \qquad \qquad  1$        &      $ 1  \qquad \qquad   1$    \\
		 10& $ 2  \qquad \qquad  2$         &   $ 1  \qquad \qquad  1$         &     $ 1  \qquad \qquad  1$   \\
\noalign{\smallskip}\hline
\end{tabular}
\end{table}

\subsection{Computing time investigation}

This section shows that the variational multirate integrator does save computing time compared to the single rate variational integrator by means of the FPU.

	For the measurements to be comparable, we simulate to the same \(t_{N} \) for all simulations.
	Furthermore, \(\Delta t = 0.001 \) is constant for all simulations and is the time step for the single rate simulations, i.e.~for $p=1$.
	From \(p=1\) we increase \(p\) in the range from \(p=5\) to \(p=1000\).
	Then, with \(\Delta t\) being constant, \(\Delta T\) increases, since \(\Delta T = p\Delta t\) holds.
	With the constant \(t_{N}\) the number of macro steps \(N\) decreases, because 
	\(N = t_{N}/\Delta T\) holds. With an increase in \(p\), the number of equations for one macro time step is increasing and a higher dimensional system of equations needs to be solved, which is more costly. This is an opposing effect to the effect of less needed macro time steps. How these two effects affect the computing time is what we look at now.

	The computing time is measured using the Matlab build in functions \emph{tic} and \emph{toc}.
	The simulations are carried out on an ``Mac $21,5"$ 2011 Edition'' with the ``Intel\textsuperscript{\textregistered} Core\textsuperscript{\texttrademark} i5 2400s 2,5GHz" processor with the operating system Mac OS X 10.6.8 (Kernel: Darwin 10.8.0) 64 bit running Matlab 2011a R7.12.0 64 bit version.
	
\begin{figure}[h]
	
		\includegraphics[width=0.49\textwidth]{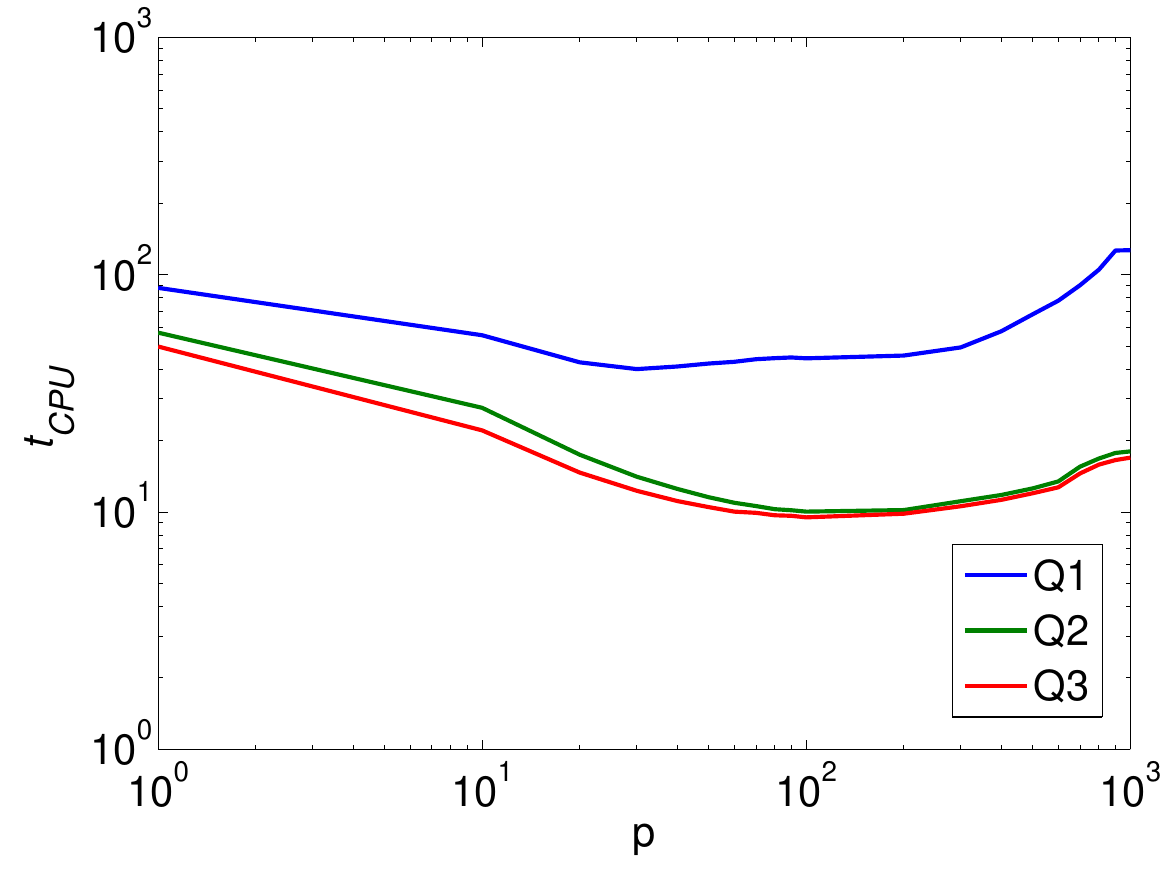}
		\includegraphics[width=0.49\textwidth]{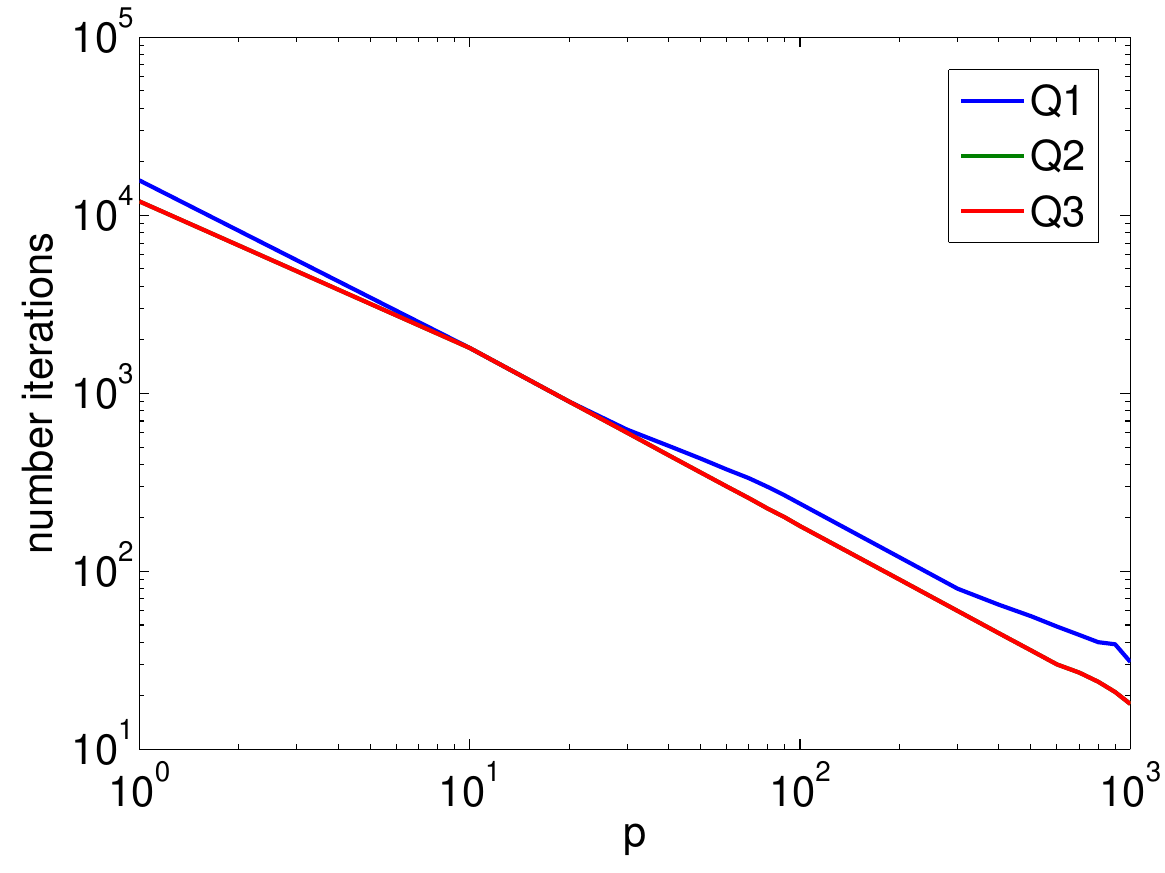}
		\caption{FPU: total computing time (left) and number of iterations (right) versus \(p\)  for midpoint-midpoint (Q1), trapezoidal-midpoint (Q2) and trapezoidal-trapezoidal (Q3) quadrature scheme with \(\Delta t =0.001\mathrm{s}\) and \(p\in \{1,5,\ldots,1000\}\)}
		\label{fig:fpu t cpu}
	\end{figure}

\begin{figure}[h]
	
		\includegraphics[width=.49\textwidth]{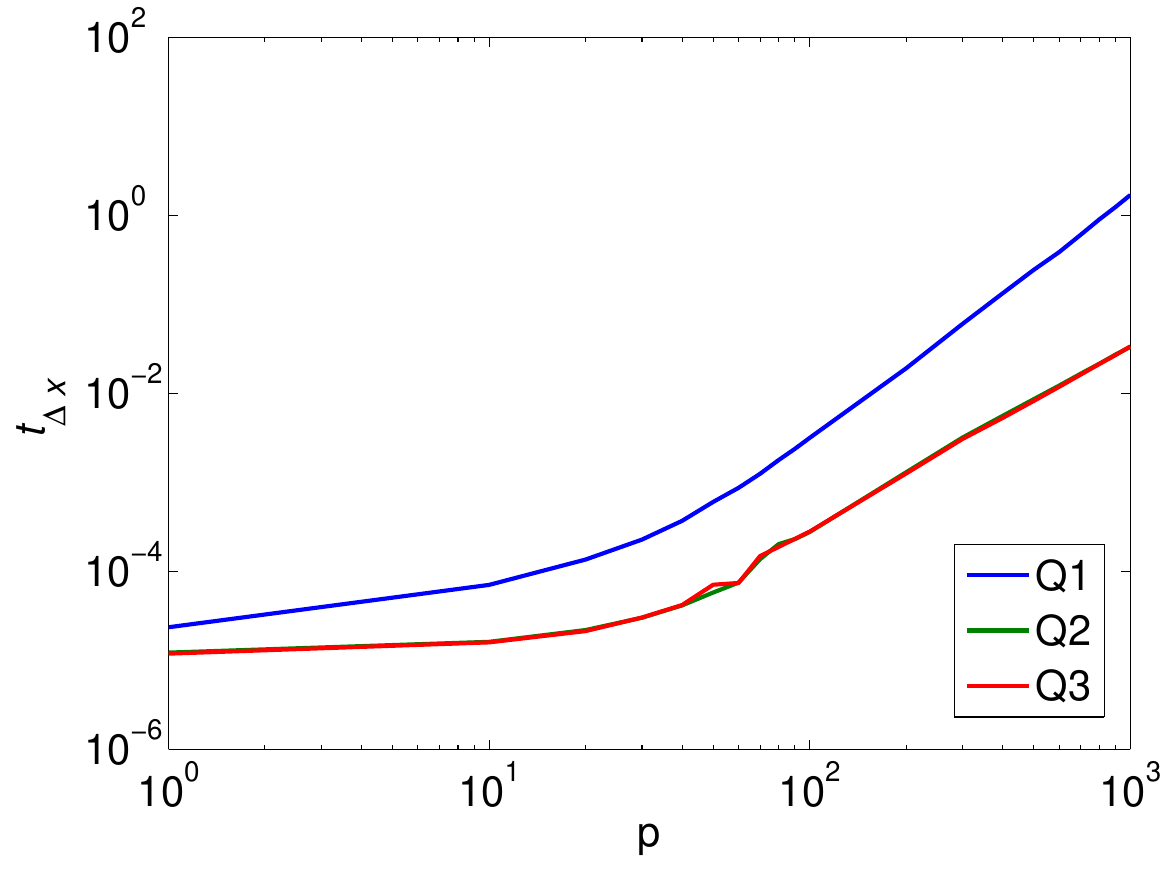}
		\includegraphics[width=.49\textwidth]{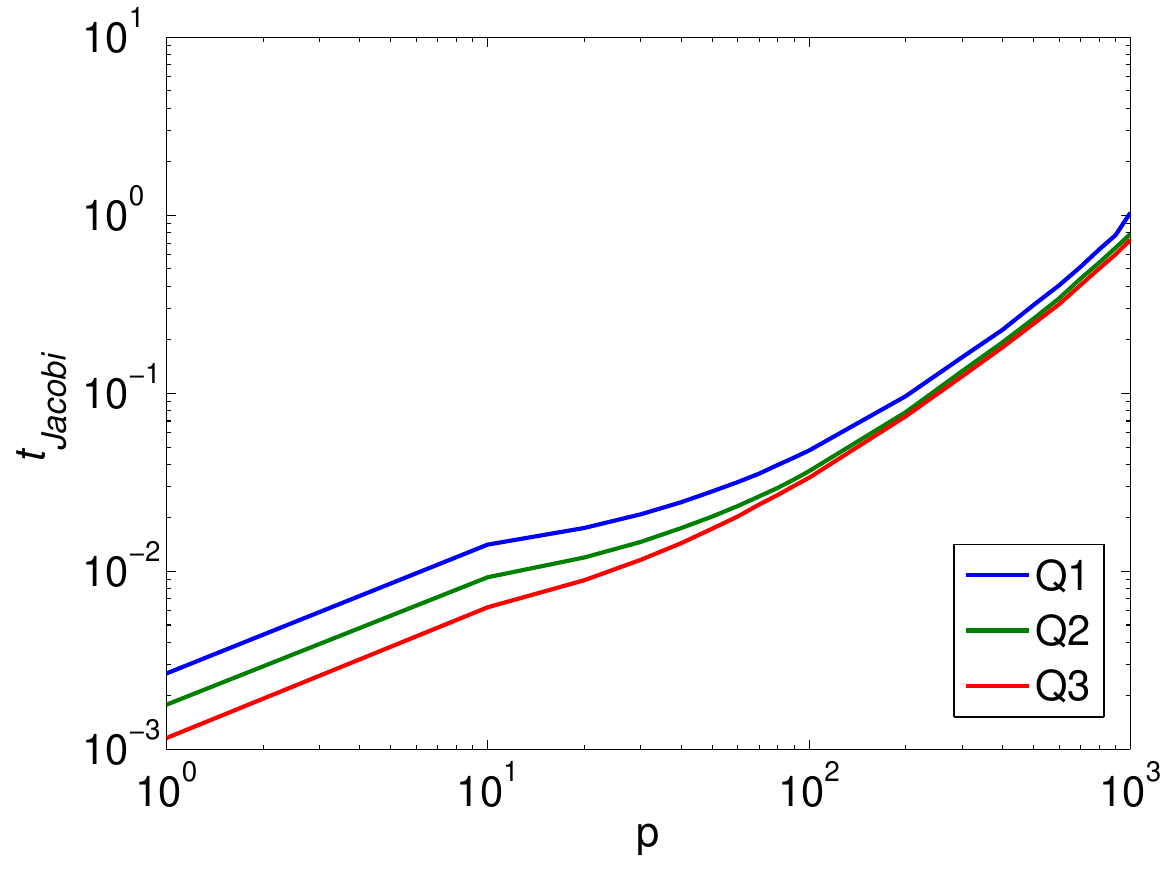}
		\caption{FPU: \(t_{\Delta x}\) (left) and \(t_{Jacobi}\) (right) versus \(p\) for all three {quadrature} schemes for \(\Delta t =0.001\mathrm{s}\) with \(p\in\{1,5,\ldots,1000\}\)}
		\label{fig:fpu t deldta x}
	\end{figure}

Fig.~\ref{fig:fpu t cpu} presents the overall computing time for all macro steps \(t_{cpu}\) (left plot) and the total number of Newton iterations (right plot) versus \(p\) for all three quadrature combinations, the midpoint-midpoint (Q1), the trapezoidal-midpoint (Q2) and the trapezoidal-trapezoidal (Q3) scheme.
We see that the variational multirate integrator compared to the single rate variational integrator $(p=1)$ saves computing time (left plot). The right plot in Fig.~\ref{fig:fpu t cpu} shows that because less macro time steps are needed, the total number of Newton iterations needed is decreasing, resulting in the seen computing time savings. However, the plot on the left also shows that with higher numbers of \(p\), there are no more computing time savings. 
To investigate this behaviour, we look at \(t_{\Delta x}\), the time to solve the nonlinear system of equations per macro step, and \(t_{Jacobi}\), the time to evaluate the Jacobian of the equations per macro step. 
Fig.~\ref{fig:fpu t deldta x} shows \(t_{\Delta x}\) (left plot) and \(t_{Jacobi}\) (right plot) versus \(p\) for all three quadrature schemes.
	Both plots show that the computing time \(t_{\Delta x}\) and \(t_{Jacobi}\) is increasing with increasing \(p\), because the dimension of the system of equations is increasing.
	We see that \(t_{\Delta x}\) is increasing super linearly, while \(t_{Jacobi}\) increases linearly.
	From these plots, we firstly see that the trapezoidal-trapezoidal scheme (Q3) needs the least computing time, while
	the trapezoidal-midpoint scheme (Q2) needs slightly more computing time. The difference in both schemes is the time \(t_{Jacobi}\), where the trapezoidal-trapezoidal scheme needs slightly less time to evaluate the Jacobian.
	Secondly, these plots explain that there are no more computing time savings for higher \(p\).
	With an increase in \(p\), less macro time steps are needed at the cost of an increase in the dimension of the system of equations.
	This increase in the dimension yields an increase in \(t_{\Delta x}\).

\section{Conclusion and future work}

In this work, integrators to efficiently approximate the dynamics of multirate systems are presented. Paralleling continuous and discrete variational multirate dynamics provides the theoretical framework to derive the integrators on a time grid consisting of macro and micro time nodes and to show the inheritance of the system's intrinsic geometric structure to the approximate solution. 
Good energy behaviour due to symplecticity and momentum map conservation are proven analytically and illustrated by numerical examples.
The choice of the quadrature rule determines the scheme and its convergence order what is also proven and demonstrated numerically. The integrators are derived and formulated as $(p,q)$-schemes making them easy to implement. A linear stability study is carried out for the case, where the non-split variable is linearly interpolated between the macro nodes together with action approximation on the micro grid, giving the range of micro steps that can be taken.
The presented computing time investigations show that a high number of micro steps and therefore large macro time steps, around 1 sec, can be chosen for the simulation of a nontrivial highly oscillatory system with nonlinear slow potential without loosing stability.
Increasing the size of the macro time step with fixed micro time step chosen appropriately to resolve the fast motion, the number of discrete slow configurations $q^s_k$ that have to be calculated over the simulation time horizon reduces and therefore computing time can be saved showing the efficiency of the presented approach.
The proposed variational multirate integrators are computationally advantageous, when the dimension of the slow variables equals the dimension of the fast variables as in the presented numerical examples, see also \citep{lishkova2020multirate}. In a situation, where the fast variables are low dimensional and the slow ones are high dimensional, the method might be even more beneficial as reducing the number of slow discrete variables throughout the simulation time horizon has even more impact.

Variational multirate integrators have successfully been extended to the framework of optimal control simulations for multirate systems. In \citep{gail_variational_2017}, multirate discrete mechanics and optimal control MDMOC is introduced. Another class of systems for which variational multirate integration is beneficial is the class of constrained systems, where in particular, slow and fast parts of the dynamics are coupled by constraints. Holonomic constraints can be enforced using Lagrange multipliers and in our experience, the numerical results are promising. However, the elimination of the constraint forces from the discrete system via the discrete null space method is not easily feasible to eliminate constraint forces that couple the slow and the fast parts. A further idea is combining time step adaption methods and multirate variational integrators.




\section*{Funding}
Tobias Gail was supported by the German Research Foundation (DFG) within the project \textit{Simulation und optimale Steuerung der Dynamik von Mehrkörpersystemen in der Biomechanik und Robotik}, GZ:LE 1841/2.

\bibliographystyle{abbrvnat}
\bibliography{references}

\end{document}